\documentclass[11pt]{article}

\usepackage{amsmath}
\usepackage{amsfonts}
\usepackage{amssymb}
\usepackage{amsthm}
\usepackage{graphics}
\usepackage{amscd}
\usepackage{graphicx,epsfig}
\usepackage{color}

\DeclareMathOperator{\Div}{div}

\DeclareMathOperator{\argch}{argch}

\renewcommand{\epsilon}{\varepsilon}
\newcommand{\boQ}{\mathcal{Q}}
\newcommand{\boV}{\mathcal{V}}
\newcommand{\boY}{\mathcal{Y}}
\newcommand{\boW}{\mathcal{W}}
\newcommand{\boF}{\mathcal{F}}
\newcommand{\boC}{\mathcal{C}}

\newcommand{\boA}{\mathcal{A}}

\newcommand{\boB}{\mathcal{B}}
\newcommand{\boO}{\mathcal{O}}
\newcommand{\boL}{\mathcal{L}}

\newcommand{\boP}{\mathcal{P}}

\newcommand{\boD}{\mathcal{D}}
\newcommand{\boN}{\mathcal{N}}
\newcommand{\Ome}{\Omega}

\newcommand{\R}{\mathbb{R}}

\newcommand{\Z}{\mathbb{Z}}
\newcommand{\C}{\mathbb{C}}

\renewcommand{\S}{\mathbb{S}}
\newcommand{\N}{\mathbb{N}}
\newcommand{\id}{\text{id}}
\newcommand{\der}[2]{\dfrac{\partial #1}{\partial #2}}
\newcommand{\dd}{\mathrm{d}}

\newtheorem{defn}{Definition}
\newtheorem{thm}{Theorem}
\newtheorem{cor}{Corollary}
\newtheorem{prop}{Proposition}
\newtheorem{lem}{Lemma}
\newcommand{\inter}[1]{\overset{\circ}{#1}}

\renewcommand{\phi}{\varphi}
\newcommand{\dis}{\displaystyle}
\newtheorem{Thm}{Theorem}[section]
\newtheorem{Prop}{Proposition}[section]
\newtheorem{Lem}{Lemma}[section]

\newcounter{remark}

\newenvironment{rem}[1]{\refstepcounter{remark}\label{#1}
  \noindent\textbf{Remark \ref{#1}.}}{}

\newcounter{case}

\newenvironment{case}[1]{\refstepcounter{case}\label{#1}
  \noindent\textbf{Case \ref{#1}.}}{}

\newcounter{construction}

\newenvironment{cons}[1]{\refstepcounter{construction}\label{#1}
  \noindent\textbf{Construction \ref{#1}.}}{}

\title{The Plateau problem at infinity for horizontal ends and genus $1$}  
\author{Laurent Mazet}
\date{}

\begin{document}

\maketitle

\begin{abstract}
In this paper, we study Alexandrov-embedded $r$-noids with genus $1$ and
horizontal ends. Such minimal surfaces are of two types and we build several
examples of the first one. We prove that if a polygon bounds an immersed
polygonal disk, it is the flux polygon of an $r$-noid with genus $1$ of the
first type. We also study the case of polygons which are invariant under a
rotation. The construction of these surfaces is based on the resolution of the
Dirichlet problem for the minimal surface equation on an unbounded domain.
\end{abstract}



\section*{Introduction}
The classical Plateau problem consists in finding a surface of least
area bounded by a given closed curve in $\R^3$, such a surface satisfies that
its mean curvature vanishes. A surface in $\R^3$ with zero mean curvature will
be called a minimal surface. The Plateau problem has been
solved by T. Rad\'o in 1930. A generalization of this problem is: finding a
minimal surface for a given asymptotic behaviour. We first give sense to this
question. 

We know that, if a complete minimal surface $M$ has finite
total curvature and $r$ embedded ends (such a surface is called an $r$-noid),
each end of this minimal surface 
is asymptotic to a plane or to a half-catenoid; besides, we can
associate to each end a vector in $\R^3$, this vector is called the
flux vector of the end. These vectors satisfy the following condition:
the sum of the flux vectors over all ends is zero. So the generalization of
the Plateau problem is: given a
finite number of vectors such that their sum is zero, can we find an
$r$-noid which has these vectors as flux vectors? This problem is
called the Plateau problem at infinity. Besides, we know that $M$ is
conformally equivalent to a compact Riemann surface $\overline{M}$ minus $r$
points, the punctures of $M$, and what we call the genus of $M$ is in fact the
genus of $\overline{M}$.  

In this paper, we shall study a particular case of this problem. In \cite{CR},
C. Cos\'\i n and A. Ros give a description of the space of 
solutions of the Plateau problem at infinity with an asymptotic
behaviour which 
is symmetric with respect to an horizontal plane (\emph{i.e.} all the
flux vectors are horizontal) for genus $0$ (the Riemann surface $\overline{M}$
is the Riemann sphere $\S^2$). In the genus $0$ case, there is a
natural order on the ends. Then, since
the flux vectors are horizontal 
and their sum is zero, the flux vectors draw a polygon in
$\R^2$; this polygon is called the flux polygon of $M$. C.~Cos\'\i n and
A.~Ros give a necessary and 
sufficient condition on this polygon for having a solution to the Plateau
problem at infinity. In our work, we
study the case where $\overline{M}$ is of genus $1$ (\emph{i.e.} $\overline{M}$
is a torus).

For the case of genus $1$, we need to distinguish two types of $r$-noid with
horizontal ends; this classification depends on the place of the punctures on
the torus: when $M$ is of the first type, there is a natural order on the
punctures and for the second type, there is not. If
$M$ is an $r$-noid of genus $1$ and horizontal ends of the first type, since
there is 
a natural order on the punctures, we can define as in the genus $0$ case the
flux polygon associated to $M$. Then our main result can be state as follow
(see Theorem \ref{main}).

\bigskip

\emph{Let $M$ be an $r$-noid with genus $0$ and horizontal ends, then there
  exists $\Sigma$ an $r$-noid with genus $1$ and horizontal ends of the first
  type which  have the same flux polygon as $M$.} 

\bigskip

In Corollary \ref{regpoly}, we give  examples of
polygons which are the flux polygons of $r$-noid with genus $1$ but not the
flux polygons of $r$-noid with genus $0$.

\medskip

The $r$-noids $M$, we consider, are symmetric with respect to a horizontal
plane that we can normalized to be $\{z=0\}$. Then to build them, it is
sufficient to build the part $M^+$ of $M$ included in $\R^2\times\R_+$. The
proof of 
our result is then based on the fact that $M^+$ is the conjugate of a minimal
surface that can be seen as the graph of a function $u$ over a planar "domain"
which 
depends on the flux polygon of $M$. This "domain" is in fact a multi-domain
(see definitions in Section \ref{defff}).

If $\Ome$ is a domain in $\R^2$ and $u$ is a function on $\Ome$, the graph of
$u$ is a minimal surface in $\R^3$ iff $u$ satisfies the elliptic partial
equation called the minimal surface equation: 

\begin{equation}{\label{MSE}}
\Div\left(\frac{\nabla u}{\sqrt{1+|\nabla u|^2}}\right)=0
\tag{MSE}
\end{equation}

Like every partial differential equation, we can associate to \eqref{MSE} the
Dirichlet problem that consists in finding a function $u$ on $\Ome$ which is a
solution of the minimal surface equation and takes on assignated values on the
boundary of $\Ome$.

The first step to build a $r$-noid with genus $1$ having a polygon $V$ as flux
polygon is then to consider a multi-domain $\Ome$ associated to this polygon
and 
solve a Dirichlet problem on $\Ome$. The boundary data is such that the graph
of the solution can be the conjugate of $M^+$. As an example if $V$ is a
triangle $ABC$, we can glue along the three edges of $V$ three half-strips
($[A,B]\times \R_+$ along $[A,B]$, $[B,C]\times \R_+$ along $[B,C]$...), we
get an unbounded domain $D$. Then the 
multi-domain $\Ome$ is the universal cover of $D\backslash \{Q\}$ where $Q$ is 
a point in the triangle $ABC$. The boundary value we take on $\Ome$ is $\pm
\infty$ alternating the sign such that for every half-strip in $\Ome$ one side
has $+\infty$ and the other has $-\infty$.

If $\Ome$ is a multi-domain associated to the polygon $V$, the conjugate
surface 
to the graph of the solution of the Dirichlet problem is a minimal surface
which is invariant under a translation by a horizontal vector. If this vector
is non zero, this minimal surface is not the piece $M^+$ of an $r$-noid $M$,
then we must 
choose a multi-domain associated to $V$ such that the corresponding vector is
zero. In our example of the triangle, the only choice we have is the position
of the point $Q$ in $ABC$, then we need to prove that there exists $Q\in ABC$
such 
that the corresponding vector vanishes. This problem is called the period
problem. 

In the case of the triangle, the idea to solve the period problem is the
following. To each point $Q$ in the interior of the triangle $ABC$, we can
associate an 
horizontal vector, in fact, this defines a continuous map from the triangle
to $\R^2$ and we want to show that this map vanishes; this map is the period
map. Then we whall compute 
the degree of the period map along the boundary of the triangle and find a
non-zero number, then the period map must vanishes at one point. The proof of
our main result is based on a generalization of this arguement.

\bigskip

The paper is organized as follows; in the first section, we define all the
notions of multi-domain we use in the following and the objects associated to
a function on such domains. We also precise what kind of $r$-noid we study in
this paper and give the first results concerning them.

Section \ref{exist} is devoted to the resolution of the Dirichlet problem on
a multi-domain $\Ome$ associated to a polygon $V$. In this section, we make
use of tools developped in \cite{Ma1} and recalled in Appendix
\ref{lineofdiv}. 

In the third section, we study the regularity of the graph of the solution of
the Dirichlet problem of the preceding section near the singularity point of
the multi-domain $\Ome$.

In Section \ref{peri}, we explain what is the period problem for the
construction of a $r$-noid. We also generalize the notion of the period map and
give the proof of our main result (Theorem \ref{main}) in using Theorem
\ref{degree} which is proved in Section \ref{extention}.

Section \ref{extention} is devoted to an extension of the period map and the
computation of its degree.

In section \ref{rotation}, we give examples of $r$-noids with genus $1$ that
are not given by Theorem \ref{main}. In particular, we consider the case where
the polygon $V$ is a regular polygon.

Let us fix some notations. In the following, when $u$ is a function on a
domain of $\R^2$ we shall note $W=\sqrt{1+|\nabla u|^2}$. We shall also use
the classical following notations for partial derivatives: $p=\der{u}{x}$ and 
$q=\der{u}{y}$. Besides, for the graph of $u$, we shall always chose the
downward pointing normal to give an orientation to the graph.


\section{Preliminaries{\label{defff}}}
\subsection{Graph on multi-domains}
In this section,we give several generalizations of the notion of domain in
$\R^2$, since our aim is to describe some minimal surfaces as the graph of a
function which is a solution of \eqref{MSE} (see the example of the half
helicoid). We use notions introduced in \cite{CR} and \cite{Ma1}.

Let us consider a pair $(D,\psi)$ where $D$ is $2$-dimensionnal flat
manifold with piecewise smooth boundary and $\psi:D\rightarrow \R^2$
is a local isometry. The map $\psi$ is called the \emph{developing map}
and the points where the boundary $\partial D$ is not smooth are
called vertices. If a part of the boundary of $D$ is linear then this
part will be called an \emph{edge} of the boundary. 

\begin{defn}
A pair $(D,\psi)$, where $D$ is a simply-connected $2$-dimensionnal
complete flat manifold with piecewise smooth boundary and $\psi:
D\rightarrow \R^2$ is a local isometry, is a \emph{multi-domain} if
each connected component of the smooth part of $\partial D$ is a
convex arc.
\end{defn}

Let $D$ be a complete metric space and $Q$ a point of $D$ we say that
$D$ admits a \emph{cone singularity} at $Q$ of angle $\alpha$ if
$D\backslash\{Q\}$ is a $2$-dimensional flat manifold and if there
exist $\rho_0>0$ such that $\{M\in D|\ d(M,Q)<\rho_0\}$ is isometric to
$\{(\rho,\theta)|\,0\le\rho<\rho_0,\,0\le\theta \le \alpha\}$ with the
polar metric $\dd s^2=\dd \rho^2+\rho^2\dd \theta^2$ where all the
points $(0,\theta)$ are identified and where for every $\rho$ we
identified $(\rho,0)$ with $(\rho,\alpha)$ (the isometry sends $Q$ to
$(0,0)$). 

\begin{defn}
A triplet $(D,Q,\psi)$ is a \emph{multi-domain with a cone singularity
  at $Q$} if 
\begin{enumerate}
\item $D$ is a simply-connected complete metric space,
\item there exists $q\in\N$ such that $D$ admits a cone singularity at
  $Q$ of angle $2q\pi$,
\item $D$ has piecewise smooth convex boundary and 
\item $\psi: D\rightarrow \R^2$ is a local isometry outside $Q$.
\end{enumerate}
\end{defn}

We can remark that a multi-domain $(D,\psi)$ can be seen as a multi-domain
with a 
cone singularity at $Q$ if $Q$ is some point of $D$. The angle at the
singularity is $2\pi$.

Let $(D,\psi)$ be a compact multi-domain such that its boundary is only
composed of edges. The developing map allows us to see $\partial D$ included in
$\R^2$ since there are only edges $\psi(\partial D)$ is a polygon in
$\R^2$. The same thing can be done for $(D,Q,\psi)$ a multi-domain with cone
singularity. We then say that a polygon $V$ bounds a multi-domain (with
perhaps a cone singularity) if there exists $(D,\psi)$ or $(D,Q,\psi)$ such
that $V=\psi(\partial D)$. When $V$ bounds a multi-domain $(D,\psi)$ we shall
also say that $v$ bounds an immersed polygonal disk as in \cite{CR}.

The last generalization we need is to give a sense to a cone singularity with
infinite angle. 

Let us consider $\boD=\{(\rho,\theta)|\ \rho\in\R_+,\ \theta\in\R\}$
with the polar metric and where all
the points $(0,\theta)$ are identified, this point will be called the
vertex of $\boD$ and noted $\boO$. The space $\boD$ is a
simply-connected complete metric space and $\boD\backslash \boO$ is a
2-dimensional flat manifold. 

\begin{defn}{\label{log}}
A triplet $(\Omega,\boQ,\phi)$ is a \emph{multi-domain with a logarithmic
singularity at $\boQ$} if 
\begin{enumerate}
\item $\Omega$ is a simply-connected complete metric space,
\item $\boQ\in\Omega$,
\item $\Omega\backslash \boQ$ is a $2$-dimensional flat manifold with
piecewise smooth convex boundary,
\item $\phi: \Omega\rightarrow \boD$ is a local isometry such that
$\phi(\boQ)=\boO$ and 
\item there exist a neighborhood $\boN$ of $\boQ$ in $\Omega$ and $\rho>0$ such
that $\phi|_\boN$ is an isometry into $\{M\in\boD\ |\ d(\boO,M)<\rho\}$. 
\end{enumerate} 
\end{defn}

Let us define $R_\alpha:\boD\rightarrow\boD$ by $R_\alpha(r,\theta)=
(r,\theta+\alpha)$, $R_\alpha$ is an isometry of $\boD$.

\begin{defn}
A multi-domain with a logarithmic singularity $(\Omega,\boQ,\phi)$ is
\emph{periodic} if there exists $f:\Omega\rightarrow\Omega$ an isometry
and $n\in\N^*$ such that 
\begin{equation*}
\phi\circ f=R_{2n\pi}\circ\phi \tag{$\ast$}\label{au}
\end{equation*}
The period of $\Omega$ is then $2\pi q$ where $q$ is the smallest $n$
such that there exists $f$ making \eqref{au} true.
\end{defn}

The first example of a multi-domain with a logarithmic singularity, we
can give, is $(\boD,\boO,\id)$. This multi-domain is periodic of period
$2\pi$.

\begin{cons}{C1}
Let us consider $(D,\psi)$ a multidomain and $Q$ a point in $D$. We then
note $\Omega \xrightarrow{\pi}D\backslash Q$ a universal cover of
$D\backslash Q$. We can pull back to $\Omega$ the flat metric of
$D$. The metric
completion of $\Omega$ is just $\Omega\cup \{\boQ\}=\overline{\Ome}$
where $\boQ$ is a point 
``above'' $Q$ (\emph{ie} if $\boA_n\rightarrow \boQ$, we have
$\pi(\boA_n)\rightarrow Q$). If $(\rho,\theta)$ are the polar coordinates
on $\R^2$ of center $\psi(Q)$ then, on $\Ome$ the $1$-forms
$(\pi\circ\psi)^*\dd \rho$ and $(\pi\circ \psi)^*\dd\theta$ are 
exact and by integration we can define a map $\phi:\Omega\cup\{\boQ\}
\rightarrow\boD$ such that $(\overline{\Ome},\boQ, \phi)$ is a
multi-domain with a logarithmic singularity. The multi-domain, we
have just build, is a periodic one of period $2\pi$. 

In fact, we can do the same work for $(D,Q,\psi)$ a multi-domain with a
cone singularity at $Q$ of angle $2q\pi$. We get 
$(\overline{\Omega},\boQ,\phi)$ a periodic multi-domain with a
logarithmic singularity (the period is less than $2q\pi$) and a
covering map $\pi:\overline{\Ome} \rightarrow D$ with $\pi(\boQ)=Q$. 
\end{cons}

\begin{cons}{quotient}
The inverse construction is also possible. Let $(\Ome,\boQ,\phi)$ be a
periodic multi-domain with a logarithmic singularity of period $2q\pi$
and isometry $f$. Then by taking the quotient of $\Ome$ by the group
$\{f^n\}_{n\in \Z}$, we
build a multi-domain with a cone singularity at $Q$, the image of
$\boQ$ in the quotient, and angle $2q\pi$.
\end{cons}

\begin{rem}{angle}
We make a remark about these two constructions. Let $(\Ome,\boQ,\phi)$
be a periodic multi-domain with a logarithmic singularity of period
$2q\pi$ and isometry $f$. If we considere the quotient of $\Ome$ by
the group $\{f^{an}\}_{n\in\Z}$ for $a\in \N^*$, we get a multi-domain
with a cone singularity $(D,Q,\psi)$. The cone singularity at $Q$ of
$D$ has $2qa\pi$ as angle. But if we apply Construction \ref{C1} to
$D$ we get $\Ome$ which have a period less than $2qa\pi$ if $a>1$.

Let $V=(v_1,\dots,v_r)$ be a polygon which, for example, bounds an immersed
polygonal disk $(D,\psi)$. If $Q\in D$, we make Construction \ref{C1} and we
get 
a multi-domain with logarithmic singularity $(\overline{\Ome}, \boQ,
\phi)$. The quotient $(D',Q',\phi')$ of $\overline{\Ome}$ by
$\{f^{2n}\}_{n\in\Z}$ is a multi-domain with cone singularity of angle $4\pi$;
besides this multi-domain bounds the polygon $(v_1,\dots,v_r,v_1,\cdots,v_r)$.
But since Construction \ref{C1} gives $\overline{\Ome}$ for $D$ and $D'$, the
two polygons $V$ and $(v_1,\dots,v_r,v_1,\cdots,v_r)$ will not be
distinguished in the following. 
\end{rem}
\smallskip

Let $(\Ome,\boQ,\phi)$ be a multi-domain with logarithmic
singularity and $A$ be a point in $\R^2$. We then can define the map
$\phi_A: \Ome \rightarrow \R^2$ by $\phi_A=G\circ\phi$ with
$G(\rho,\theta)=A+(\rho\cos\theta, \rho\sin\theta)$. If $\Ome$ is
given by Construction \ref{C1} we always choose $A=\psi(Q)$.

\medskip

Let $\Ome$ such that $(\Ome,\phi)$ or $(\Ome,Q,\phi)$ corresponds to
one of the three definitions of multi-domain given above. Let $u$ be a
function defined on $\Ome$ or $\Ome$ minus its singularity. The
\emph{graph} of $u$ is then the surface in $\R^3$ defined by
$\{\phi(x),u(x)\}_{x\in\Ome}$ or $\{\phi_A(x),u(x)
\}_{x\in\Ome\backslash\{\boQ\}}$ with $A\in\R^2$. In the following
the function $u$ will often be a solution of the minimal surface
equation, in this case the graph of $u$ becomes a minimal surface of
$\R^3$. The fact that $u$ is a solution of \eqref{MSE} allows us to
define a closed $1$-form $\dd\Psi_u$ on $\Ome$, $\dd\Psi_u$ is the inner
product $\frac{\nabla u}{W}\lrcorner \dd V$ where $\dd V$ is the volume form
on $\Ome$. Since $\dd\Psi_u$ is
closed we can define locally a function $\Psi_u$ (obviously $\Psi_u$ is well
defined only if we fix its value at one point). $\Psi_u$ is locally defined in
the interior of $\Ome$ minus the possible singularity, but $\Psi_u$ is
$1$-Lipschitz countinuous then it can be countinuously extended to the
singularity and the boundary, then since $\Ome$ is simply connected $\Psi_u$
is then globally defined on $\Ome$. In fact $\Psi_u$ correponds to the third
coordinates of the conjugate surface to the graph of $u$; $\Psi_u$ is called
the conjugate function to $u$. For other properties on $\Psi_u$ we refer to
\cite{JS} and Appendix \ref{lineofdiv}.

In \cite{JS} and \cite{Ma1}, we can find the most general answer to the
Dirichlet problem on compact multi-domain $(D,\psi)$: the Dirichlet
problem consists in finding a solution $u$ on $D$ of \eqref{MSE}
knowing its value on the boundary.

On multi-domain with cone or logarithmic singularity there is no
general answer. To give an exemple of solution of \eqref{MSE} on a
multi-domain with logarithmic singularity, let us consider the function
$u$ defined on $(\boD,\boO,\id)$ by $u(\rho,\theta)=\theta$, it is
obvious that the graph of $u$ is the half of an helicoid; more
precisely it is the surface given in isothermal coordinate by
$(a,b)\mapsto (\sinh a\cos b,\sinh a\sin b,b)$ for
$(a,b)\in\R_+^*\times \R$. The function $u$ is then a solution of the
minimal surface equation.

\subsection{The $r$-noids of genus $1$}

In this section we give precise definitions of $r$-noids, flux at one end and
other objects linked to the Plateau problem at infinity. We also give some
results concerning this problem and explain how we can build solutions in the
genus $1$ case.

Let M be a complete minimal surface with finite total curvature in
$\R^3$; we know that $M$ is isometric to a compact Riemann surface
$\overline{M}$ 
minus a finite number of points (we can refer to \cite{Os}). Then $M$ has a
finite number of annular ends; when these ends are embedded they are
asymptotic either 
to a half-catenoid or to a plane. A properly immersed minimal surface
with $r$ embedded ends will be called  a \emph{$r$-noid}. We can
associate to each end a vector which caracterizes the direction and the 
growth of the asymptotic half-catenoid (when the end is asymptotic to
a plane this vector is zero); this vector is called the \emph{flux} of
the end (for a precise definition of the flux see \cite{HK}). If
$v_1,\dots,v_r$ are the fluxes at each end, we have the following balancing
condition: 
\begin{equation}{\label{balcond}}
v_1+\cdots+v_r=0
\end{equation}
This condition tells us that the total flux of the system vanishes. If 
$v_1,\dots,v_r$ are vectors in $\R^3$ such that \eqref{balcond} is
verified and $g$ is a non-negative integer, the Plateau problem at
infinity for these data is to find an $r$-noid of genus $g$ which has
$v_1,\dots,v_r$ as fluxes at its ends (the genus $g$ is the genus of
$\overline{M}$). 

Let $X:M\longrightarrow\R^3$  be an $r$-noid. $M$ is conformally
equivalent to a compact surface $\overline{M}$ minus $r$
points $p_1,\dots,p_r$. We will say that $M$ is
\emph{Alexandrov-embedded} if $\overline{M}$ bounds a compact $3$-manifold
$\overline{\Omega}$ and the immersion $X$ extends to a proper local
diffeomorphism $f:
\overline{\Omega}\backslash\{p_1,\dots,p_r\}\longrightarrow \R^3$. An
Alexandrov-embedded surface has a canonical orientation; we choose the 
Gauss map to be the outward pointing normal. An Alexandrov-embedded
$r$-noid can not have a planar end (see \cite{CR}).

We are interested in the case where $X:M\longrightarrow \R^3$ is an
Alexandrov-embedded $r$-noid of genus $g$ and $r$ horizontal ends
(\emph{ie} the flux at each end is an horizontal vector).

Let $X: M\longrightarrow\R^3$ be a nonflat immersion of a connected 
orientable surface $M$ and $\Pi$ be a plane in $\R^3$, normalized to
be $\{x_3=0\}$. We note by $S$ the Euclidiean symmetry with respect to 
$\Pi$ and consider the subsets:
\begin{gather*}
M^+=\{p\in M | x_3(p)>0\}\\
M^-=\{p\in M | x_3(p)<0\}\\
M^0=\{p\in M | x_3(p)=0\}
\end{gather*}
With these notation we have:
\begin{defn}\label{stsym}
We shall say that $M$ is strongly symmetric with respect to $\Pi$ if
\begin{itemize}
\item There exists an isometric involution $s:M\longrightarrow M$
  such that $\psi\circ s=S\circ \psi$.
\item $\{p\in M|s(p)=p\}=M^0$.
\item The third coordinate $N_3$ of the Gauss map of $M$ takes
  positive (resp. negative) values on $M^+$ (resp. $M^-$).
\end{itemize}
\end{defn}

In \cite{CR}, C. Cos\'\i n and A. Ros prove

\begin{prop}{\label{strongsym}}
Let $M$ be an $r$-noid with horizontal ends. Then $M$ is strongly
symmetric with respect to an horizontal plane if and only if $M$ is
Alexandrov-embedded.
\end{prop}

The notion of strong symmetry is then important for the study of the
Alexandrov-embedded $r$-noid.

The case of genus $0$ was studied by C.~Cos\'\i n and A.~Ros in
\cite{CR}; they show that in this case there is a natural order on the
ends. Let $M$ be an Alexandrov-embedded $r$-noid of genus $0$ if
$2v_1,\dots,2v_r$ are the fluxes of $M$ ordered as the 
ends, the polygon $(v_1,\cdots,v_r)$ is called the \emph{flux polygon}
of $M$ and is noted $F(M)$. We then have
\begin{thm}
Let $v_1,\dots,v_r$ be horizontal vectors such that
$v_1+\cdots+v_r=0$ and $V$ the associated polygon, then there exists
$M$ an Alexandrov-embedded $r$-noid of genus $0$ such that $F(M)=V$
if, and only if, $V$ bounds an immersed polygonal disk. Besides there
is a bijection from the set of $M$ such that $F(M)=V$ and the set of
the immersed polygonal bounded by $V$.
\end{thm}
In the following, when $(\boP,\psi)$ is an immersed polygonal disk we
shall call $\Sigma(\boP)$ the $r$-noid of genus $0$ associated to
$\boP$ by this bijection; we refer to \cite{CR} and \cite{Ma1} for more
explanations on this theorem.

In this paper we are interested in the case of genus $1$. Let $M$ be
an Alexandrov-embedded $r$-noif of genus $1$ with horizontal ends. $M$
is conformally a torus $\overline{M}$ minus $r$ points
$p_1,\dots,p_r$. By Proposition \ref{strongsym}, $M$ is strongly
symmetric with respect to the plane $\{z=0\}$ (in fact $M$ is strongly
symmetric with respect to an horizontal plane that we can normalize to
be $\{z=0\}$). As in \cite{CR}, the punctures $p_i$ are fixed by the
antiholomorphic involution $s$ given by Definition \ref{stsym}; we
have the following lemma.

\begin{lem}
Let $\overline{M}$ be a conformal torus and $s$ an antiholomorphic
involution on $\overline{M}$. We suppose that $s$ has a fixed point
then the set of the fixed points of $s$ is two separated circles.
\end{lem}

We can then give the following definition.

\begin{defn}
Let $M$ be an Alexandrov-embedded $r$-noid with genus one and
horizontal ends and $p_i$ the punctures of $M$. We note $s$ the
antiholomorphic involution associated to $M$, we then have $C_1$ and
$C_2$ the two circles of fixed points of $s$. With this notations, we
shall say that: 
\begin{itemize}
\item $M$ is of \emph{type I} if all the $p_i$ are in one of the two
  circles $C_1$ and $C_2$,
\item $M$ is of \emph{type II} if not.
\end{itemize}
\end{defn}

We suppose now that $M$ is of type I, the circles of fixed points are the
boundary of $M^+\cup M^0$, this minimal surface is oriented by the outward
pointing normal and then its boundary has a natural orientation. Then the
circle that contains the points $p_i$ has a
natural orientation and we suppose that the points $p_i$ are numbered with
respect to this orientation. We note $2v_i$ the flux vector associated to
$p_i$. We know that we have $\dis\sum_{i=1}^rv_i=0$, so, as in the case $g=0$,
the list $(v_1,\dots,v_r)$ defines a polygon which we call the \emph{flux
  polygon} of $M$ and we note it $F(M)$. We can remark that if $M$ is 
of type II we do not have such an easy definition.

In the following, we use some arguements developped by C.~Cos\'\i n
and A.~Ros in \cite{CR} to the study of the genus $0$ case.
Let $X:\overline{M}\backslash \{p_1,\dots,p_r\}\rightarrow M$ be an
Alexandrov-embedded $r$-noid with genus one and horizontal 
ends of type I. The surface $M^+\cup M^0$ is topologically an annulus, so
by passing to the universal cover, we can define its conjugate
surface  $M_{0,+}^*$. The surface $M_{0,+}^*$ is a periodic minimal
surface and its period vector is $\dis\int_\gamma\dd X^*$ where
$\gamma$ is a path in $M^+\cup M^0$ that generates $\pi_1(M^+\cup
M^0)$. Since $M$ is of type I all the $p_i$ are in one circle of fixed
points, then the other circle of fixed points generates $\pi_1(M^+\cup M^0)$
and since 
the normal along this circle is horizontal the first two
coordinates of the period vector are zero. This proves that the period
vector is a vertical vector.

The symmetry curves $M^0$ consists of $r$ complete strictly convex curves
in the plane $\{z=0\}$ (they are the images by $X$ of the arcs
$p_ip_{i+1}$) and  the image by $X$ of the circle of fixed points that
contains no $p_i$. In the conjugate surface, these curves are
transformed in vertical lines; the image of the arc $p_ip_{i+1}$ in
the conjugate surface are vertical straight-lines over the vertices of
the polygon $F(M)$

Since on $M^+$ the third coordinate of the normal is positive, the
projection map $\Pi$ from $M_+^*$, the conjugate of $M^+$, into the
plane 
$\{z=0\}$ is a local diffeomorphism. So we can pull back to $M_+^*$
the flat metric of this plane.
Besides $M_+^*$ is stable by the translation $t$ of vector
$\dis\int_\gamma\dd X^*$. Since this vector is vertical, $\Pi\circ t=\Pi$; $t$
is then an isometry 
of $M_+^*$ with the flat metric $\Pi^*(\dd s_{\R^2}^2)$.

Using arguements of C.~Cos\'in and A.~Ros, we can prove that there
exists $(\Ome,\boQ,\phi)$ a multi-domain with logarithmic singularity
such that $M_+^*$ with the flat metric can be seen as the interior of
$\Ome$ minus the singularity point $\boQ$. Besides $\Ome$ is periodic
because of the existence of the isometry $t$. The boundary of $\Ome$
is composed of half-lines. $M_+^*$ is then a graph over $\Ome$ such
that the line which is the conjugate of the circle of fixed points
that contains no $p_i$ is the part of the boundary of the graph which
is above $\boQ$ and the conjugates of the arcs $p_ip_{i+1}$ are lines
over the vertices of $\Ome$.

The quotient of $\Ome$ by the group $(t^n)_{n\in\Z}$ is then a
multi-domain with cone singularity $(D',Q,\psi)$. $D'$ has $r$
vertices $P_1,\dots,P_r$ and is bounded by $2r$ half-lines; more
precisely, $D'$ is a multi-domain with cone singularity $D$ which is
bounded by the flux polygon $F(M)=(v_1,\dots,v_r)$ (we have $v_i=
\overrightarrow{\psi(P_i)\psi(P_{i+1})}$) to which we have glued $r$
half-strips (along the edge $[P_i,P_{i+1}]$, we glue a half-strip
isometric to $[P_i,P_{i+1}]\times \R_+$). This proves that the flux
polygon $F(M)$ bounds a multi-domain with cone singularity.

\medskip

In Section \ref{exist} and Section \ref{regularite}, we shall
prove many results that explains, when a 
polygon $V$ bounds a multi-domain with cone singularity, how we can
build a candidate for a surface $M_+^*$ such that $F(M)=V$. 

The main result we prove concerning the Plateau problem at infinity
for genus one is the following.

\begin{thm}{\label{main}}
Let $v_1,\dots,v_r$ be $r$ non zero vectors of $\R^2$ such that
$(v_1,\dots,v_r)$ is a polygon that bounds an immersed polygonal disk. 
Then there exists $M$ an
Alexandrov-embedded $r$-noid with genus one and horizontal ends of
type I such that $F(M)=(v_1,\dots,v_r)$. 
\end{thm}

This implies that every polygon that can be realized as the flux
polygon of an $r$-noid of genus $0$ is also the flux polygon of an
$r$-noid of genus $1$.


\section{A Dirichlet problem{\label{exist}}}
Let $v_1,\dots,v_r$ be $r$ non zero vectors of $\R^2$ such
that $V=(v_1,\dots,v_r)$ bounds a multi-domain with cone
singularity $(D,Q,\psi)$. Let us note $P_1,\dots,P_r$ the vertices of the
polygon, we put $P_{r+1}=P_1$; in the following, we suppose that the
orientation on the polygon $V$ is the one given as boundary of
$D$ (We remark that, in the following, the vertices of $V$ will be
identified with the vertices of $D$). Following Construction \ref{C1},
we get a multi-domain with logarithmic singularity $(\boW,\boQ,\phi)$
with a projection map $\pi:\boW \rightarrow D$. Because of Remark
\ref{angle}, we suppose that the period $2q\pi$ of $\boW$ is equal to the
angle of the cone singularity of $D$, in fact we suppose that $D$ is
the quotient of $\boW$ by the isometry $f$; In the following we shall say that
$D$ satisfies the hypothesis H.

\smallskip

\begin{cons}{C2}
Let $i\in\{1,\dots,r\}$ and consider $E$ an edge of $\boW$ which is
send by $\pi$ to the edge $[P_i,P_{i+1}]$ of $D$; we can glue to
$\boW$ along $E$ a 
half strip $S_i$ isometric to $[P_i,P_{i+1}] \times \R_+^*$ (if
$A\in E$ the point $A$ is identified with $(\pi(A),0)$). Making
this for every $i$ and every edges $E$, we get a
new multi-domain with a logarithmic singularity at $\boQ$ that we
note $(\Ome,\boQ,\phi)$ (we keep the same notation for the developping
map since it is an extension of the original developping map). Since
$\boW$ is periodic, $\Ome$ is periodic and have the same period, we
note $f$ the corresponding isometry.

Let us note, for $i\in\{1,\dots,r\}$, $\boL_i^+$ (resp. $\boL_i^-$)
the union of the straight-lines corresponding to $\{P_i\}\times
\R_+^*$ in the half-strips glued along the edges $E$ such that
$\pi(E)=[P_i,P_{i+1}]$
(resp. $\pi(E)=[P_{i-1},P_i]$). $\boL_i^+$ and $\boL_i^-$ are a
countable union of half straight-lines. We note $\boV$ the set of the
vertices of $\Ome$ this set is $\pi^{-1}\{P_1,\dots,P_r\}$.
\end{cons}

We then have the following existence and uniqueness result.

\begin{thm}{\label{graph}}
Let $(D,Q,\psi)$ be a multi-domain with cone singularity that bounds a polygon
$V$, Constructions \ref{C1} and \ref{C2} give us a periodic multi-domain with
logarithmic singularity $(\Ome, \boQ,\phi)$; we suppose that the period of
$\Ome$ is the cone angle of $D$ at $Q$. Then there exists a solution $u$ of
the minimal surface equation on $\Omega$ such that
\begin{enumerate}
\item $u$ tends to $+\infty$ along $\boL_i^+$ and $-\infty$ along
$\boL_i^-$ and
\item $\Psi_u(\boQ)=0=\Psi_u(\boV)$,
\end{enumerate}
where $\boL_i^+$, $\boL_i^-$ and $\boV$ are the notations given in
Construction \ref{C2}. Besides, the solution is unique up to an additive
constant. 
\end{thm}

The conditions imposed to a function $u$ solution of this Plateau problem make
that the graph of $u$ is a good candidate for giving a solution to the Plateau
problem at infinity with $V$ as flux polygon.

\begin{cor}{\label{period}}
If $u$ is a solution on $\Ome$ of the Dirichlet problem asked in Theorem
\ref{graph}, there exists a constant $c\in \R$ such that $u\circ f=u+c$. 
\end{cor}

\begin{proof}
Let $u$ be a solution of the Dirichlet problem asked in Theorem
\ref{graph}, then $u\circ f$ is also a solution of this Dirichlet
problem. This proves that $u\circ f - u$ is constant.
\end{proof}

The following of this section is devoted to the proof of Theorem \ref{graph}

\subsection{Notations}

First, the period of $\Ome$ is $2q\pi$.
We can suppose that $P_1$ is such that $d(Q,P_1)=\min_id(Q,P_i)$. Let
us consider a vertex in $\boV$ which is a lift of $P_1$ and denote it
$\boP_1(0)$. Since $d(M,P_1)$ is minimal, the geodesic joining $\boQ$
to 
$\boP_1(0)$ is embedded, this implies that the first coordinate of
$\phi(\boP_1(0))$ is $d(\boQ,\boP_1(0))=d(Q,P_1)$. By considering
$R_\alpha\circ \phi$ instead of $\phi$, we can suppose that
$\phi(\boP_1(0))=(d(\boQ,\boP_1(0),0)$. Let us note
$\boP_1(k)=f^k(\boP_1(0))$ for $k\in\Z$; we then have
$\phi(\boP_1(k))=(d(M,P_0),2kq\pi)$. Then $\{\boP_1(i),
i\in\Z\}$ is the set of the vertex of $\Ome$ corresponding to $P_1$
(\emph{i.e.} $\{\boP_1(i),i\in\Z\}=\pi^{-1}(P_1)$). 

If we remove the geodesics $[\boQ,\boP_1(k)]$ and $[\boQ,\boP_1(l)]$
($k<l$) from $\Ome$, we get a space which have three connected
components; one of them is such that its intersection with $\boN$, the
neighborhood of $\boQ$ introduced in Definition \ref{log}, is
isometric with $]0,r[\times ]2kq\pi,2lq\pi[$ by $\phi$; we shall note 
this part $\Ome_k^l$. For $n\in \Z$ we have $f^n(\Ome_k^l)=
\Ome_{k+n}^{l+n}$. In $\Ome_0^1$, there is exactly one lift of every
$P_i$ ($2\le i\le r$), we note $\boP_i(0)$ the lift of $P_i$ that is in
$\Ome_0^1$. We then note for $k\in\Z$ $\boP_i(k)=f^k(\boP_i(0))$; we
then have that $\boP_i(k)$ is the only lift of $P_i$ in
$\Ome_k^{k+1}$. A part $\Ome_k^{k+1}$ of $\Ome$ is called a
\emph{period} of $\Ome$.

For $i\in\{1, \dots,r\}$ and $ k\in\Z$, let us note $\boL_i^+(k)$
(resp. $\boL_i^-(k)$) the half straight-line included in $\boL_i^+$
(resp. $\boL_i^-$) and having $\boP_i(k)$ as end-point. 

Let us consider $k$ and $l$ in $\Z$  such that $k<l$, then $\Ome_k^l$
is a multi-domain. Its vertices are $\boQ$, the $\boP_1(m)$ for $k\le
m \le l$ and the $\boP_i(m)$ for $2\le i\le r$ and $k\le m<l$. Its
boundary is composed of two segments $[\boQ,\boP_1(k)]$ and
$[\boQ,\boP_1(l)]$ and $2r(l-k)$ half straight-lines which are the
$\boL_1^+(m)$ for $k\le m <l$, the $\boL_1^-(m)$ for $k< m \le l$ and
the $\boL_i^\pm(m)$ for $2\le i\le r$ and $k\le m<l$. From
$\Ome_k^l$, we now define a new multi-domain: we can glue to
$\Ome_k^l$ two half-strips, one is isometric to $[\boQ,\boP_1(k)]
\times \R_+$ and is glued along $[\boQ,\boP_1(k)]$, the second is
isometric to $[\boP_1(l),\boQ] \times \R_+$ and is glued along
$[\boP_1(l),\boQ]$. We note $\widetilde{\Ome}_k^l$ this new
multi-domain. We note $\widetilde{\boL}_k^-$
(resp. $\widetilde{\boL}_l^+$) the new half straight-line in the
boundary of $\widetilde{\Ome}_k^l$ with $\boP_1(k)$
(resp. $\boP_1(l)$) as end-point. We also note $\widetilde{\boL}^+$
and $\widetilde{\boL}^-$ the two half straight-lines in the boundary
with $\boQ$ as end-point such that $\widetilde{\boL}^+$ is in the
same half-strip than $\widetilde{\boL}_k^-$.

\subsection{Proof of the existence}

We shall now prove the existence part of Theorem \ref{graph}. First,
for $n$ in $\N^*$, we prove that there exists a function $u_n$ on
$\widetilde{\Ome}_{-n}^n$ which:
\begin{enumerate}
\item is a solution of \eqref{MSE},
\item tends to $+\infty$ on $\widetilde{\boL}^+$,
  $\widetilde{\boL}_n^+$ and all the $\boL_i^+(k)$ that are in the
  boundary of $\widetilde{\Ome}_{-n}^n$ and
\item tends to $-\infty$ on $\widetilde{\boL}^-$,
  $\widetilde{\boL}_{-n}^-$ and all the $\boL_i^-(k)$ that are in the
  boundary of $\widetilde{\Ome}_{-n}^n$.
\end{enumerate}

Let us consider the following polygon: 
$$(\overrightarrow{QP_1}, \underbrace{ v_1,\dots,v_r,\dots
  \dots,v_1,\dots,v_r}_{2n\ \textrm{times}},
\overrightarrow{P_1Q})$$  
If we remove to $\widetilde{\Ome}_{-n}^n$ all the half-strips, we get a
multi-domain which is bounded by the above polygon. Then, the
existence of the solution $u_n$ is ensured by Theorem $7$ in \cite{Ma1}. 

Let us now only consider the restriction of $u_n$ to $\Ome_{-n}^n$. We
then have a sequence of solutions of \eqref{MSE} defined in an
increasing sequence of domains $\Ome_{-n}^n$ and
$\dis\bigcup_{n\in\N^*} \Ome_{-n}^n=\Ome$. We can then consider that
$(u_n)$ is a sequence of functions on $\Ome$. We want $(u_n)$ to converge,
so we shall prove that there is no line of divergence (see Appendix
\ref{lineofdiv}). 

First we give some preliminary results on $u_n$.

\begin{lem}{\label{jenkins}}
Let $u$ be a solution of \eqref{MSE} on the half-strip
$[0,a]\times\R_+$ such that $u$ tends to $-\infty$ on
$\{0\}\times\R_+$ and $+\infty$ on $\{a\}\times \R_+$. Then, for $y\ge
4a$, we have:
\begin{gather*}
\frac{|p|}{W}(x,y)\ge1-\frac{a^2}{y^2}\\
\frac{|q|}{W}(x,y)\le \sqrt{2}\frac{a}{y}
\end{gather*} 
\end{lem}

\begin{proof}
We note $A$ the point of coordinates $(x,y)$ and $B$ a point in the
boundary of the half-strip which realizes the distance from $A$ to this
boundary. Since $y\ge 4a$ the coordinates of $B$ are $(0,y)$ or
$(a,y)$; besides the distance $|AB|$ is less than $a/2$. Because of
the value of $u$ on the boundary, the distance along the graph from
the point above $A$ to the boundary of the graph is bigger than
$4a$. The ratio of this two distances is less than $1/8$, then we can
apply Lemma $1$ in \cite{JS}; this gives the lemma.
\end{proof}

\begin{cor}{\label{P=Q}}
Let $u$ be a solution of \eqref{MSE} on the half-strip
$[0,a]\times\R_+$ such that $u$ tends to $-\infty$ on
$\{0\}\times\R_+$ and $+\infty$ on $\{a\}\times \R_+$. We have
$\Psi_u(0,y) =\Psi_u(a,y)$ and if $\Psi_u(0,0)=0$ then $\Psi_u\ge 0$
in the half-strip.
\end{cor}

\begin{proof}
Because of the value of $u$ on the boundary we have
$\Psi_u(0,y)=\Psi_u(0,0)+y$ and $\Psi_u(a,y)=\Psi_u(a,0)+y$. This
implies that for $y\ge 4a$ we have
\begin{equation*}
\begin{split}
|\Psi_u(0,0)-\Psi_u(a,0)|&=|\Psi_u(0,y)-\Psi_u(a,y)|\\
&=|\int_0^a-\frac{q}{w}(x,y)\dd x|\\
&\le\int_0^a \sqrt{2}\frac{a}{y}\dd x\le \sqrt{2}\frac{a^2}{y}
\end{split}
\end{equation*}
Then by letting $y$ goes to $+\infty$ we have
$\Psi_u(0,0)=\Psi_u(a,0)$ and then $\Psi_u(0,y)=\Psi_u(a,y)$.

If $\Psi_u(0,0)=0$ then $\Psi_u(0,y)=y$. If $A$ is in the interior of the
half-strip, there exists $y$ such that $A$ is at a distance less than
$y$ from the point $(0,y)$; then, since $\Psi_u$ is $1$-Lipschitz
continuous $\Psi_u(A)\ge 0$. 
\end{proof}

This corollary implies that, if we choose $\Psi_{u_n}$ such that
$\Psi_{u_n}(\boQ)=0$, we have $\Psi_{u_n}(\boP_i(k))=0$ for all
$\boP_i(k)\in \Omega_{-n}^n$ and $\Psi_{u_n}\ge 0$ in the
half-trips contained in $\widetilde{\Ome}_{-n}^n$. Since $\Psi_{u_n}$ satisfies
a maximum principle, we have 
$\Psi_{u_n}\ge 0$ in $\Ome_{-n}^n$.

Assume that there exists a line of divergence $L$. By Lemma
\ref{lemdiv1}, $L$ 
can not have an end point in the interior of one of the half
straight-lines that compose the boundary of $\Omega$; this prove that
if it has end-points it must be $\boQ$ or one $\boP_i(k)$.

If $L$ has
no end point, let $\boA$ be a point of $L$ and let us note $d$ the
distance between $\boA$ and $\boQ$, since $\Psi_{u_n}(\boQ)=0$ and
$\Psi_{u_n}$ is $1$-Lipschitz, we have $0\le \Psi_{u_n}(\boA) \le
d$. Let us fix an orientation to $L$ such that the limit normal along
this line of divergence is the right-hand unit normal. We then note
$s$ the arc-length along $L$ with $\boA$ as origin. Let us consider
$\boA'$ the point on $L$ of arc-length $s=-2d$. Then we know that, for
the subsequence that makes $L$ appear, we have $\Psi_{u_{n'}}(\boA)-
\Psi_{u_{n'}}(\boA') \rightarrow 2d$ but, since $\Psi_{u_n}(\boQ)\le d$
this implies that $\Psi_{u_{n'}}(\boA')<0$ for big $n'$; this is a
contradiction.

If $L$ is a segment ($L$ has two end-points noted $\boA_1$ and
$\boA_2$), we know that for each $n$ we have $\Psi_{u_n}(\boA_1)=0=
\Psi_{u_n}(\boA_2)$, but for the subsequence that make $L$ appear we
have $|\Psi_{u_{n'}}(\boA_1)-\Psi_{u_{n'}}(\boA_2)| \rightarrow
d(\boA_1,\boA_2)>0$; this gives us a contradiction.

We then can suppose that $L$ has only one end-point that we note
$\boF$ and goes to infinity in one half-strip that we can parametrized
isometricaly by $[0,a]\times\R_+$. We remark that for one half-strip
in $\Ome$ the number of such lines $L$ is finite. We have, for each
$n$, 
$\Psi_{u_n}(\boF)=0$. There exists $b\in]0,a[$ such that the part of
$L$ in the half-strip is $\{b\}\times\R_+$. By changing $L$ if
necessary we can suppose that the part $]b,a[\times\R_+$ is included
in $\boB(u_n)$. Let $\boA=(b,0)$, $\boB=(b,2(a-b))$, $\boC=(a,2(a-b))$
and $\boD=(a,0)$ be four points in the half-strip. Since
$\dd\Psi_{u_n}$ is closed we have:

\begin{equation*}
\begin{split}
\int_{[\boA,\boB]}\dd\Psi_{u_n}&= -\int_{[\boB,\boC]}\dd\Psi_{u_n}
-\int_{[\boC,\boD]}\dd\Psi_{u_n}- \int_{[\boD,\boA]}\dd\Psi_{u_n}\\
&=2(a-b)-\int_{[\boB,\boC]}\dd\Psi_{u_n}-
\int_{[\boD,\boA]}\dd\Psi_{u_n}\\
&\ge 2(a-b)-(a-b)-(a-b)=0 
\end{split}
\end{equation*}

This implies that we have only one possibility for the limit
normal . Since $]b,a[\times\R_+\subset \boB(u_n)$ we can
suppose that a subsequence $(u_{n'})$ converges to a function $v$ on
$]b,a[\times\R_+$ ($n'$ is chosen such that the line of divergence $L$
appears). $v$ is a solution of \eqref{MSE} and by Lemma 
\ref{lemdiv2}, we know 
that $v$ tends to $+\infty$ along $\{b\}\times \R_+^*$ and $-\infty$
along $\{a\}\times \R_+^*$. Then by Lemma \ref{P=Q}, $\Psi_v(\boA)=
\Psi_v(\boD)= 
\lim \Psi_{u_n}(\boD)=0$. But, for the subsequence, $\Psi_{u_{n'}}(\boA)-
\Psi_{u_{n'}}(\boF) \rightarrow d(\boF,\boA)>0$, this contradicts the fact
that $\Psi_{u_n}(\boF)=0$ and $\lim \Psi_{u_{n'}}(\boA)=0$. 

We then have prove that $\boB(u_n)=\Omega$, and, by taking a subsequence, we
can suppose that $u_n$ converges to a solution $u$ of \eqref{MSE} on
$\Omega$. By Lemma \ref{lemdiv2}, $u$ is such that $u$ tends to
$+\infty$ along the 
$\boL_i^+(k)$ and $-\infty$ along the $\boL_i^-(k)$. By construction,
we have also that $\Psi_u(\boQ)=0=\Psi_u(\boP_i(k))$ and $\Psi_u\ge 0$.

\subsection{Property of the solution 1}

Before proving the uniqueness, we need to give some properties of a
solution of the Dirichlet problem asked in Theorem \ref{graph}. We
recall that there exist $\rho_0$ and a neighborhood $\boN$ of $\boQ$
such that $\phi$ is an isometry from $\boN$ into $\{M\in \boD|
d(\boO,M)\le\rho_0\}$. We then can use the polar coordinates
for $v$ a function defined on $\boN$.

\begin{prop}{\label{major}}
Let $\epsilon$ be a positive number. There exists $d>0$ such that for 
every $\alpha\in\R$ and every $v$ solution of \eqref{MSE} on $\boN
\cap \{\alpha-\pi<\theta<\alpha+\pi\}$:
\begin{gather*}
\sup_{0<\rho\le \rho_0}v(\rho, \alpha) < \sup_{ [d,\rho_0]\times
  [\alpha-3\pi/4,\alpha+3\pi/4]} v(\rho, \theta) +\epsilon\\
\inf_{0<\rho\le \rho_0}v(\rho, \alpha) > \inf_{ [d,\rho_0]\times
  [\alpha-3\pi/4,\alpha+3\pi/4]} v(\rho, \theta) -\epsilon
\end{gather*}
\end{prop}

The proof is based on the idea used by R.~Osserman to prove Theorem 10.3 in
\cite{Os}. 

\begin{proof}
First let us assume that $\alpha =\pi/2$. Then $\boN \cap
\{\alpha-\pi<\theta<\alpha+\pi\}$ is isometric to the disk of center
the origin and radius $\rho_0$ minus the segment joining the origin to
$(0,-\rho_0)$. Then we only prove the first inequality since the second is
a consequence of the first by substituting $v$ for $-v$.

Let $\epsilon>0$ and $d$ be a positive number, we note $A$ the point
of coordinates $(0,-d)$ and $D$ the domain in $D(0,\rho_0)$ delimited
by the segment joining $(-\sqrt{\rho_0^2-d^2},-d)$ to $(-d,-d)$, the
half circle of center $A$  and radius $d$: $\theta\mapsto
(d\cos\theta, d(\sin\theta -1))$ for $\theta\in [0,\pi]$ and  the
segment joining $(d,-d)$ to $(\sqrt{\rho_0^2-d^2},-d)$ and containing
$(0,\rho_0/2)$. On the domain $D$ we consider the function $c: M
\mapsto d \left(-\argch(\frac{|AM|}{d})+
\argch(\frac{\rho_0+d}{d})\right)$. The graph of $c$ is a piece of a
catenoid and $c$ is a positive function upper bounded by
$d\argch(\frac{\rho_0+d}{d})$. Let $d$ be small enough such that this
upper-bound is less than $\epsilon$. 

\begin{figure}[h]
\begin{center}
\resizebox{0.8\linewidth}{!}{\input{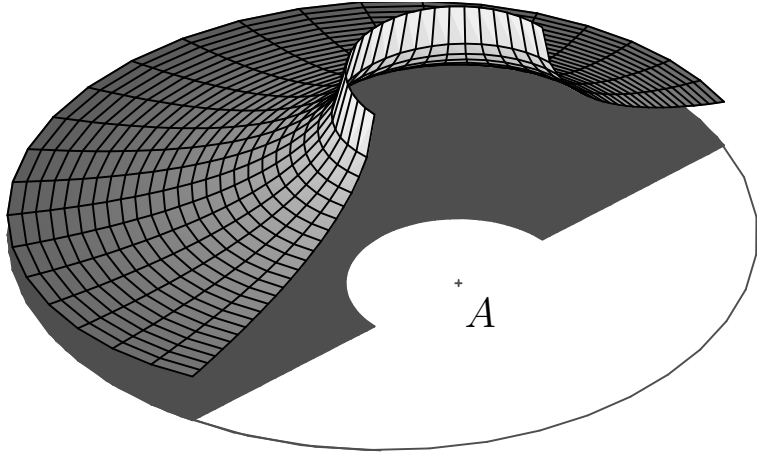}}
\caption{the graph of the function $c$}
\end{center}
\end{figure}

On the part of the boundary of $D$ which is a half circle of center
$A$, the derivatives $\nabla c\cdot n$, where $n$ is the outward
pointing normal, takes the value $+\infty$. This implies, by Lemma $10.2$ in
\cite{Os}, 
that $c+K\ge v$ where $K$ is the supremum of $v$ on the
boundary of $D$ minus the half circle of center $A$. Since this part
of the boundary is included in the set of points of polar coordinates
included in $[d,\rho_0]\times [-\pi/4,5\pi/4]$ and $D$ contains the
segment joining the origin to $(0,1)$ the proposition is established.     
\end{proof}

This proposition is an important fact since it implies that if $u$ is
a solution of the dirichlet problem asked in Theorem \ref{graph}, $u$
is bounded on $\boN\cap \Ome_k^l$ for every $k$ and $l$.

\subsection{Proof of the uniqueness} 

Let us now prove the uniqueness part of Theorem \ref{graph}. We
consider $u$ and $v$ two solutions of the Dirichlet problem asked in
Theorem \ref{graph} and we assume that these two solutions are
different: $u-v$ is not constant. By changing $v$ by $v+c$ where $c$
is a real constant we can suppose:
\begin{enumerate}
\item $\{u-v>0\}\cap\Ome_0^{+\infty}$ and
  $\{u-v<0\}\cap\Omega_0^{+\infty}$ are non-empty ($\Omega_0^{+\infty}$ is
  $\dis\bigcup_{n>0}\Ome_0^n$) and
\item the segment $[\boQ,\boP_1(0)]$ is included in $\{u-v>0\}$.
\end{enumerate}
The first assertion is due to the fact that $u-v$ is non-constant and
that we can exchange $u$ and $v$, the second assertion is a
consequence of the fact that $u$ and $v$ are bounded on the segment
$[\boQ,\boP_1(0)]$ by Proposition \ref{major} and Lemma $2$ in \cite{Ma1}.

We note $\Delta=\{u-v<0\}\cap\Omega_0^{+\infty}$
We then note $\Delta^n$ the intersection of $\Delta$ and $\Omega_0^n$.
The boundary
$\partial\Delta^n$ is composed of a part included in the boundary of
$\Omega$, a part included in the interior of $\Omega_0^n$ and a part
included in the segment $[\boQ,\boP_1(n)]$. Along the first part
$\dd\widetilde{\Psi}= \dd\Psi_u-\dd\Psi_v=0$ and along the second part
$\dd\widetilde{\Psi}$ is positive by Lemma $2$ in \cite{CK}; the union of
this first 
two parts will be noted $\partial^n\Delta$.

We note $\Delta_n^l$ the
intersection of $\Delta_n$ with the part of $\Ome_0^n$ such that a point $A$
is in this part if either $A$ is not in one half-strip included in $\Ome_0^n$
or $A=(x,y)$ is in one of these half-strips which is isometricaly parametrized
by 
$[0,a]\times\R_+$ and $y\le l$. The boundary of $\Delta_n^l$ is composed of 
three parts: the first is $\partial\Delta_n\cap [\boQ,\boP_1(n)]$, the
second is $\partial^n\Delta\cap \overline{\Delta_n^l}$ and the third is
included in the union of the segments parametrized by $[0,a]\times\{l
\}$ in each half-strip; we note $\Gamma^l$ this third part. Since $\dd
\widetilde{\Psi}$ is closed we have:
\begin{equation}{\label{fermeture}}
\int_{\partial\Delta_n\cap [\boQ,\boP_1(n)]} \dd \widetilde{\Psi}+
  \int_{\partial^n\Delta\cap \overline{\Delta_n^l}} \dd \widetilde{\Psi}+
    \int_{\Gamma^l} \dd \widetilde{\Psi}=0 
\end{equation}
We then have 
$$\int_{\partial^n\Delta\cap \overline{\Delta_n^l}} \dd
\widetilde{\Psi} \le |\boQ\boP_1(n)|+\left|\int_{\Gamma^l} \dd
  \widetilde{\Psi} \right|$$
Since $\dd \widetilde{\Psi}\ge 0$ along $\partial^n\Delta$,
$\int_{\partial^n\Delta\cap \overline{\Delta_n^l}} \dd
\widetilde{\Psi}$ increases as $l$ increases. By Lemma \ref{jenkins},
$\int_{\Gamma^l} \dd \widetilde{\Psi}\rightarrow 0$ as $l$ goes to
$+\infty$. Then the above equation implies that $\dd \widetilde{
  \Psi}$ is integrable on $\partial^n\Delta$ and by passing to the
limit in \eqref{fermeture}, we obtain
\begin{equation}
\int_{\partial\Delta_n\cap [\boQ,\boP_1(n)]} \dd \widetilde{\Psi}+
  \int_{\partial^n\Delta} \dd \widetilde{\Psi}=0
\end{equation}
This equation implies that
$$\int_{\partial^n\Delta} \dd \widetilde{\Psi}\le
|\boQ\boP_1(n)|=|\boQ\boP_1(0)|$$
Then  $\dd \widetilde{\Psi}$ is integrable on $\partial\Delta$ and
$\int_{\partial\Delta} \dd \widetilde{\Psi}>0$ since the part of
$\partial\Delta$ included in the interior of $\Ome$ is non empty. We have
$$\int_{\partial\Delta} \dd \widetilde{\Psi}=
\sum_{n=0}^{+\infty}\int_{\partial\Delta \cap \Omega_n^{n+1}} \dd
\widetilde{\Psi}$$
Then $\dis\int_{\partial\Delta \cap \Omega_n^{n+1}} \dd \widetilde{\Psi}
\longrightarrow 0$. We then want to understand what happens on
$\Omega_n^{n+1}$ for big $n$.

Let us consider, for $n\in \N$, $u_n=u\circ f^{-n}$ and $v_n=v\circ
f^{-n}$; the restriction of $u_n$ to $\Ome_0^1$ is equal to the
restriction of $u$ to $\Ome_n^{n+1}$ and the same is true for
$v$. With the same arguements as in the proof of existence, we can
prove that there exist two real sequences $(a_n)$ and $(b_n)$ such that
$u_{n'}-a_{n'}\rightarrow \tilde{u}$ and $v_{n'}
-b_{n'}\rightarrow \tilde{v}$ on $\Ome$ (in fact $a_n=u_n(P)$ and $b_n=v_n(P)$
for a point $P$ in $\Ome$). The functions
$\tilde{u}$ and $\tilde{v}$ are two 
solutions of the Dirichlet problem asked in Theorem \ref{graph} and
the convergence is uniform for every derivative on compact parts. By
changing our subsequence if necessary, we can assume that
$b_{n'}-a_{n'}\rightarrow \pm \infty$ or
$b_{n'}-a_{n'}\rightarrow c\in\R$. We are interested in
$\{u_n-v_n<0\}\cap \Omega_0^1$ and, in a certain sense, these sets
``converge'' to $\{\tilde{u}- \tilde{v}<\lim b_{n'}-
a_{n'}\} \cap \Omega_0^1$.   

Let $\gamma:[0,|\boQ\boP_1(1)|]\rightarrow [\boQ,\boP_1(1)]$ the
parametrization by arc-length of the segment $[\boQ,\boP_1(1)]$, with
$\gamma(0)=\boQ$.

First we assume that $b_{n'}-a_{n'}\rightarrow \pm
\infty$. Let $\epsilon>0$ such that $\epsilon$ is less than
$\dis\frac{1}{3} \int_{\partial^n\Delta} \dd \widetilde{\Psi}$ and
$\dis\frac{1}{2}|\boQ,\boP_1(1)|$. Since $\tilde{u}$ and $\tilde{v}$
are bounded on $[\boQ,\boP_1(1)]$ and that we have uniform convergence
on $\gamma([\epsilon,|\boQ\boP_1(1)|-\epsilon])$, we can ensure that,
for $n'$ big enough, $\gamma ([\epsilon,|\boQ\boP_1(1)| -\epsilon])$ is
included in $\{u_{n'}-v_{n'}<0\}$ or does not intersect
$\{u_{n'}-v_{n'}<0\}$ following the sign of the limit of
$b_{n'}-a_{n'}$. For every $n$,
$\dis\int_{[\boQ,\boP_1(n)]} \dd \widetilde{\Psi}=0$ because of our
hypotheses on $\Psi_u$ and $\Psi_v$, we then have for big $n'$
\begin{equation}{\label{2epsilon}}
\left| \int_{\partial\Delta_{n'+1}\cap [\boQ,
    \boP_1(n'+1)]} \dd\widetilde{ \Psi}\right|\le 2\epsilon
\end{equation}
This equality implies that $\dis\int_{ \partial^{n'+1}\Delta}
\dd\widetilde{\Psi} \le 2\epsilon$ for big $n'$. Then by passing to the
limit we obtain a contradiction because of our hypothesis on $\epsilon$.

We assume now that $b_{n'}-a_{n'}\rightarrow c$. Let
$\epsilon$ be as above. If the segment $\gamma([\epsilon,
|\boQ\boP_1(1)|-\epsilon])$ is included in $\{ \tilde{u}-
\tilde{v}<c\}$ or does not intersect $\{\tilde{u}- \tilde{v}<c\}$ then
the same property is true for $\{u_{n'}-v_{n'}<0\}$ for
big $n$. Then same arguements as above give us a contradiction. We
then can ensure that there is a point in $\gamma([\epsilon,
|\boQ\boP_1(1)|-\epsilon])$ where $\tilde{u}-\tilde{v}=c$. If
$\tilde{u}=\tilde{v}+c$ on $\Omega$, we have $\dd
\Psi_{u_{n'}}-\dd\Psi_{v_{n'}}\rightarrow 0$ uniformly on
$\gamma([\epsilon, |\boQ\boP_1(1)|-\epsilon])$. Then for $n'$ big
enough, the inequality \eqref{2epsilon} is true, this gives us a
countradiction. Then $\tilde{u}\neq \tilde{v}+c$ on $\Omega$ and there
is a compact part of the boundary of $\{\tilde{u}- \tilde{v}<c\}$ that
crosses the segment $\gamma([\epsilon, |\boQ\boP_1(1)|-\epsilon])$,
let us note $\Gamma$ this part of the boundary ($\Gamma$ is oriented as
the boundary of this set). We have
$\dis\int_\Gamma \dd\Psi_{\tilde{u}}-\dd \Psi_{\tilde{v}}>0$ by Lemma
$2$ in \cite{CK}. The curve $\Gamma$ is included in a compact part
of $\Omega_0^2$ 
and since $(u_{n'}-v_{n'})$ converges to $\tilde{u}-
\tilde{v}-c$ (the convergence is uniform for each derivative on every compact
part) we can ensure that for $n'$ big enough
$$\int_{\partial\{u_{n'}-v_{n'}<0\}\cap \Ome_0^2} \dd
\Psi_{u_{n'}}- \dd \Psi_{v_{n'}} \ge \frac{1}{2}
\int_\Gamma \dd\Psi_{\tilde{u}}-\dd \Psi_{\tilde{v}}$$

This implies that
$$\int_{\partial\Delta \cap \Omega_{n'}^{n'+2}} \dd
\widetilde{\Psi} \ge \frac{1}{2} \int_\Gamma \dd\Psi_{\tilde{u}}-\dd
\Psi_{\tilde{v}}>0$$ 
We have then a contradiction with $\dis\int_{\partial\Delta \cap
  \Omega_n^{n+1}} \dd \widetilde{\Psi} \longrightarrow 0$. We have
then proved that our two solutions $u$ and $v$ differ only by a real
additive constant. 

\subsection{Property of the solution 2}

From Corollary \ref{period}, if $u$ is a solution of the Dirichlet
problem asked in Theorem \ref{graph}, there exists a constant $c\in\R$
such that $u\circ f=u+c$. We have the following result for such a
situation.  

\begin{prop}{\label{cno=0}}
Let $(\Ome,\boQ,\phi)$ be a periodic multi-domain with logarithmic singularity,
we note $f$ the isometry associated to the periodicity. Let $v$ be a
solution of the minimal surface equation on $\Ome$ such that:
\begin{enumerate}
\item there exists a constant $c\in \R$ such that $v\circ f=v+c$ and 
\item if $\Psi_v$ is the conjugate function to $v$ with $\Psi_v(\boQ)=0$,
  $\Psi_v$ is non-negative.
\end{enumerate}
Under these hypotheses, the constant $c$ is non zero.
\end{prop}

\begin{proof}
Let us suppose that the constant $c$ is zero and consider $v$ only on
the neighborhood $\boN$ of the singularity point $\boQ$. Since $v\circ f=v$, by
taking the quotient with respect to $f$, the function $v$ can be seen
as a function $\widetilde{v}$ defined on $\boC$ a flat disk of radius $\rho_0$
with a cone 
singularity at its center of angle $2q\pi$ minus the singularity point
(i.e. $\boC$ is $\{(\rho,\theta),\rho\in (0,\rho_0), \theta\in [0,
2q\pi]\}$ where $(\rho,0)$ and $(\rho,2q\pi)$ are identified, the
metric on $\boC$ is the polar metric). The graph of $\widetilde{v}$ is
a minimal surface of $\boC \times \R$. This surface is topologicaly an
annulus then it is conformally parametrized by a Riemann surface $R$
which is an annulus: there are two harmonic maps
$X:R\rightarrow \boC$ and $x_3:R\rightarrow \R$ such that
$x_3=\widetilde{v} \circ X$. We can suppose that $R$ is either
$\{\zeta\in\C|1<|\zeta|<a\}$ for $1<a\le +\infty$ or
$\{\zeta\in\C|0<|\zeta|<a\}$ for $0<a\le +\infty$. On $\boC$, the
function $\Psi_{\widetilde{v}}$ is well defined and satisfies
$\Psi_{\widetilde{v}} \ge 0$ and $\Psi_{\widetilde{v}}=0$ at the cone
singularity by hypothesis. Since $x_3$ is harmonic, on $R$ we can
define 
its harmonic conjugate $x_3^*$; \emph{a priori} $x_3^*$ is
multi-valuated but, for a good choice of $x_3^*$, we have
$x_3^*=\Psi_{\widetilde{v}}\circ X$ then $x_3^*$ is well defined on $R$. As
$\Psi_{\widetilde{v}}=0$ at the cone singularity, $x_3^*(\zeta)$ converges to
zero as either $|\zeta|$ tends to $1$ when $R=\{\zeta\in\C|1<|\zeta|<a\}$ or
$|\zeta|$ tends to zero when $R=\{\zeta\in\C|0<|\zeta|<a\}$. 

In the second case ($R=\{\zeta\in\C|0<|\zeta|<a\}$), $x_3^*$ is a
harmonic function on a pointed disk, which has a continuous extension
to the whole disk. But, an harmonic function can not have an isolated
singularity, then the continuous extension is harmonic and have an
extremum at the origin, since $x_3^*\ge 0$ on $R$; this gives a
contradiction.

In the other case ($R=\{\zeta\in\C|1<|\zeta|<a\}$), since $x_3^*=0$ on
the unit circle, the function $x_3^*$ extends to an harmonic function on
$\{\zeta\in\C|\frac{1}{a}<|\zeta|<a\}$ by Schwartz reflection
principle (the extension is defined by $x_3^*(\zeta)=-x_3^*(
\frac{1}{\overline{\zeta}})$; see \cite{ABR}). By taking the harmonic
conjugate, $x_3$ 
also extends to $\{\zeta\in\C|\frac{1}{a}<|\zeta|<a\}$. Along the unit
circle $\S^1$, $\nabla x_3^*\cdot n=0$, with $n$ the unit tangent
vector to $\S^1$, since $x_3^*=0$ along this circle. Besides by Rolle's Theorem
there is a point on $\S^1$ where 
$\nabla x_3\cdot n=0$. At this point, $\nabla x_3^*=0$ and the
local structure of a critical point of a non constant harmonic function implies
that $x_3^*$ must be negative on $R$ in the neighborhood of this
point, this contradicts $x_3^*\ge 0$ on $R$.  
\end{proof}

\begin{rem}{quotient1}
This proposition proves that, for a solution $u$ of the Dirichlet problem
asked in Theorem \ref{graph}, the constant $c$ such that $u\circ f=u+c$ is
non-zero, then $u$ can not pass to the quotient by $f$. But
since $u\circ f=u+c$ the derivatives of $u$ are invariant by $f$ then
the derivatives of $u$ are well defined on the quotient. Besides, this
implies that the function $\Psi_u$ is also invariant by $f$.
\end{rem}

An other property of such a solution $u$ is that $u$ can be build by
the way used in the proof of the existence. As a consequence, the
conjugate function $\Psi_u$ is non-negative on $\Ome$.


\section{The regularity over the singularity point{\label{regularite}}}
As in the preceding section we consider $V$ a polygon that bounds a
multi-domain with cone singularity $(D,Q,\psi)$; using Construction \ref{C1}
and \ref{C2}, we get a multi-domain with logarithmic singularity
$(\Ome,\boQ,\phi)$; we suppose that $D$ satisfies the hypothesis H.
Theorem \ref{graph} allows us to construct a minimal surface in
$\R^3$ which is a graph over $\Ome$; we want to use this surface to build a
$r$-noid with genus $1$ and horizontal ends of type I with $V$ as flux polygon
therefore the graph needs to be regular up to its boundary. In this section, we
understand the behaviour of this minimal surface over the singularity
point $\boQ$: we shall prove the following result.
\begin{thm}{\label{reg1}}
Let $u$ be a solution of the Dirichlet problem asked in Theorem
\ref{graph}. Then the surface $\{\phi_{\psi(Q)}(x),u(x)\}$ for $x\in \boN$ is
a minimal surface with boundary and its boundary is
$\{\phi_{\psi(Q)}(x),u(x)\}$ for $x\in \partial \boN$ and the vertical
straight line passing through $\psi(Q)$. 
\end{thm}

In fact, this surface can also be seen as a minimal surface in
$\boD\times\R$; with this point of view, Theorem \ref{reg1} can be
stated as follows. 
\begin{thm}{\label{reg2}}
Let $u$ be a solution of the Dirichlet problem asked in Theorem
\ref{graph}. Then the surface $\{\phi(x),u(x)\}_{x\in\boN}\subset
\boD\times\R$ is a minimal surface with boundary and its boundary is
$\{\phi(x),u(x)\}$ for $x\in \partial \boN$ and the vertical
straight line $\{\boO\}\times\R$.
\end{thm}

It is easy to see that Theorem \ref{reg2} implies Theorem \ref{reg1},
then we shall prove Theorem \ref{reg2}. The problem in this theorem is the
behaviour near the singularity point not along $\partial \boN$.

\subsection{A compactness result}

In the proof of Theorem \ref{reg2}, we shall need to make converge a
sequence of minimal surfaces. The following theorem will be an
important tool. 

\begin{thm}{\label{comp}}
Let $\zeta_1$, $\zeta_2$ and $\zeta_3$ be three different points in $\S^1$. Let
$\gamma_n:\S^1\rightarrow\boD\times\R^{k_1}$ be a sequence of Jordan curves
in $\boD\times\R^{k_1}$ and $f_n:\S^1\rightarrow\R^{k_2}$ be a sequence of
continuous maps; we note $g_n=(\gamma_n,f_n)$. We note
$\Gamma_n=\gamma_n(\S^1)$ and put, for 
$1\le i\le 3$, $p_i^n=\gamma_n(\zeta_i)$. We suppose that:
\begin{enumerate}
\item $\dis I=\inf_{n\in\N} \{d(p_1^n,p_2^n), d(p_1^n,p_3^n),
  d(p_2^n,p_3^n)\}>0$, 
\item for every $m>0$ there exists $\epsilon>0$ such that, for every
  $n\in\N$, if $\zeta,\zeta'\in \S^1$ and
  $d(g_n(\zeta),g_n(\zeta'))<\epsilon$, one component of
  $\Gamma_n\backslash \{\gamma_n(\zeta),\gamma_n(\zeta')\}$ is of
  diameter less than $m$ and
\item there exist $M>0$ and a sequence $X_n:\Delta\rightarrow
  \boD\times\R^{k_1+k_2}$ ($\Delta$ is the unit disk) such that
  $g_n=X_n|_{\S^1}$ and 
$$\int_\Delta (|{X_n}_x|^2+|{X_n}_y|^2)\dd x\dd y<M$$ 
\end{enumerate}

The family $\{\gamma_n\}_{n\in\N}$ is then equicontinuous.
\end{thm}

By Arzela's theorem, this implies that if, for example, the sequence
$(p_1^n)$ converges there exists a subsequence $\gamma_{n'}$ that
converges for uniform convergence. Theorem \ref{comp} is in fact very similar
to a classical result used in the resolution of the classical Plateau problem
(see for example, Lemma 3.2 in \cite{Cou} or \cite{Hi}).

\begin{proof}
The proof is based on the following lemma.
\begin{lem}[Courant-Lebesgue]{\label{CL}}
Let $X$ be of class $C^0(\Delta,\boD\times\R^k)\cap
C^1(\inter{\Delta},\boD\times\R^k)$ ($\inter{\Delta}$ is the interior of
$\Delta$) and satisfy 
$$\int_\Delta (|X_x|^2+|X_y|^2)\dd x\dd y<M$$
for some $M\in\R_+$. Then, for every $\zeta_0\in\S^1$ and for each
$\delta\in(0,1)$, there exists a number
$\rho\in(\delta,\sqrt{\delta})$ such that the distance between the images
$X(\zeta)$, $X(\zeta')$ of the two intersection points $\zeta$ and $\zeta'$ of 
$\S^1$ 
with the circle $\partial D_\rho(\zeta_0)$ can be estimated by
$$d(X(\zeta),X(\zeta'))\le \left(\frac{4M\pi}{\ln(1/\delta)}\right)^{1/2}$$
\end{lem}

The proof of this lemma can be found in \cite{Hi}.

Let $e$ be a positive number, we suppose that $e< I$. By the second
hypothesis, there exists $\epsilon$ such that, for every $n\in\N$, if
$\zeta,\zeta'\in \S^1$ and $d(g_n(\zeta),g_n(\zeta'))<\epsilon$, one
component of $\Gamma_n\backslash \{\gamma_n(\zeta),\gamma_n(\zeta')\}$ 
is of diameter less than $e$. Let $\delta>0$ be such that 
$$\left(\frac{4M\pi}{\ln(1/\delta)}\right)^{1/2}\le \epsilon$$
and for every $\zeta\in\S^1$ we have $|\zeta-\zeta_i|>\delta$ for at
least two of the points $\zeta_1,\zeta_2,\zeta_3$. 

Let $n\in\N$ and $\zeta_0\in\S^1$, by the third hypothesis and the
Courant-Lebesgue Lemma, there exists $\delta<\rho<\sqrt{\delta}$ such
that, if $\zeta$ and $\zeta'$ are the two intersections points of
$\S^1$ with 
the circle $\partial D_\rho(\zeta_0)$, $d(X_n(\zeta),X_n(\zeta'))\le
\epsilon$. $\S^1\backslash \{\zeta,\zeta'\}$ is composed of two arcs:
one, $A'$, 
contains $\zeta_0$ and all the points that are at a distance less than
$\delta$ from $\zeta_0$, the second arc, $A''$, contains two of the
three 
points $\zeta_1,\zeta_2,\zeta_3$. Since $d(X_n(\zeta),X_n(\zeta'))\le
\epsilon$, one of 
the two arcs $\gamma_n(A')$, $\gamma_n(A'')$ is of diameter less than
$e$; but $\gamma_n(A'')$ contains two of the points $p_1^n$, $p_2^n$,
$p_3^n$ and $e<I$, then it is $\gamma_n(A')$ that is of diameter less
than $e$. We then have proved that if $|\zeta-\zeta_0|<\delta$,
$d(\gamma_n(\zeta),\gamma_n(\zeta_0))\le e$; this proves that the
family 
$\{\gamma_n\}$ is equicontinuous.
\end{proof}

\subsection{Proof of Theorem \ref{reg2}}
\subsubsection{Preliminaries}
Let $u$ be a solution of the Dirichlet problem asked in Theorem
\ref{graph}. We use the notations introduced in Section
\ref{exist}. The function $u$ is contructed as the limit of a sequence 
$(u_n)$ where $u_n$ is a solution of a Dirichlet problem on
$\widetilde{\Omega}_{-n}^n$. Since we are only interested in the
behaviour on $\boN$ we shall use the polar coordinates given by $\phi$; $u$
is then defined on $[0,\rho_0] \times \R$ and $u_n$ is defined on
$[0,\rho_0] \times [-2nq\pi-\frac{\pi}{2},
  2nq\pi+\frac{\pi}{2}]$. Because of the periodicity of $u$, to prove
Theorem \ref{reg2}, it is enough to make the proof on a period
$[0,\rho_0]\times [0,2q\pi]$. 

By construction, $u_n$ takes the value $+\infty$ on $(0,\rho_0)\times
\{-2nq\pi-\frac{\pi}{2}\}$ and the value $-\infty$ on
$(0,\rho_0)\times 
\{2nq\pi+\frac{\pi}{2}\}$, besides, by taking  $\Psi_{u_n}(\boQ)=0$,
we have $\Psi_{u_n}\ge 0$ on $\widetilde{\Omega}_{-n}^n$, this proves,
by Theorem 3 in \cite{Ma1}, that the minimal surface
$\{\phi(x),u_n(x)\}_{x\in\boN 
\cap \widetilde{\Omega}_{-n}^n} \subset \boD\times\R$ has the vertical
straight-line passing by $\boO$ as boundary. The idea of the proof of Theorem
\ref{reg2} is to follow the behaviour of the graph near this vertical
straight-line when $n$ goes to $+\infty$.

We now need a result on the behaviour of graph bounded by vertical
line. 
\begin{lem}{\label{bord}}
Let $v$ be a solution of \eqref{MSE} on a sector of $\boD$
$\{(r,\theta)\in\boD|\ r\le r_0,\ \alpha_1\le \theta\le \alpha_2\}$
($\alpha_1<0<\alpha_2$). Suppose that the graph of $v$ in $\R^3$ is a complete
minimal
surface with boundary and the part of the boundary over the origin is
an interval of the vertical straight-line passing by the origin. Then
if $v$ is bounded on $\{\theta=0\}$, $\dis\lim_{r\rightarrow 0}v(r,0)$
exists. Besides the normal to the graph $v$ at the points $(0,0,
\dis\lim_{r\rightarrow 0}v(r,0))$ is $\pm(0,1,0)$.
\end{lem}
\begin{proof}
We note $\Sigma$ the graph of $v$. All the cluster points of
$(r,0,v(r,0))$ as $r$ goes to $0$ are in the boundary of the graph and
more precisely in the part of the boundary consisting in the interval of the
vertical straight-line passing 
by the origin. Besides, the curve $r\mapsto (r,0,v(r,0))$ is in the
intersection of the vertical plane $\{y=0\}$ and the graph of $v$,
this curve is then in the intersection of two minimal surfaces. Let
$(0,0,a)$ be a cluster point of $(r,0,v(r,0))$ as $r$ goes to $0$, since
the boundary of $\Sigma$ is a vertical straight-line near $(0,0,a)$
the normal to $\Sigma$ at this point is horizontal and is
$(\cos\alpha,\sin\alpha,0)$ for some $-\pi\le\alpha<\pi$. Near
$(0,0,a)$, $\Sigma$ is then a graph over the vertical plane
$\{x\cos\alpha+y\sin\alpha= 0\}$ and is tangent to this plane at
$(0,0,a)$. Then, if $|\alpha|\neq\pi/2$, the intersection of $\Sigma$
and $\{y=0\}$ is only the vertical straight-line near $(0,0,a)$ then
no $(r_n,0,v(r_n,0))$ can converge to $(0,0,a)$. This implies that the
normal at $(0,0,a)$ is $\pm(0,1,0)$ and then since $\Sigma$ is a graph
over $\{y=0\}$ the intersection of $\Sigma$ with this plane near
$(0,0,a)$ is the vertical straight-line and some smooth curves passing
by $(0,0,a)$ such that $x\neq 0$ along them. One of this curve is then
$(r,0,v(r,0))$, by continuity; this prove that $\dis\lim_{r\rightarrow
0}v(r,0)=a$.  
\end{proof}

We want to understand the convergence of the graph of $u_n$ near the
singularity point $\boO$. let $\alpha\in\R$, to simplify our notation
in the following we suppose that $\alpha\in[0,2q\pi]$, let $k\in\N$
and note $\alpha_1=\alpha-2kq\pi$ and $\alpha_2=\alpha+2kq\pi$. We
know 
that $u_n(r,\alpha_1)$ and $u_n(r,\alpha_2)$ are bounded by Lemma $2$ in
\cite{Ma1}. Then 
since the graph of $u_n$ has a vertical straight-line as boundary,
$u_n(r,\alpha_1)\rightarrow a_n$, $u_n(r,\alpha_2)\rightarrow b_n$ and
$u_n(r,\alpha)\rightarrow c_n$ as $r$ goes to $0$, for some $c_n\in
(b_n,a_n)$. The graph  
$\Sigma_n\subset\boD\times\R$  
of $u_n$ over $U= [0,\rho_0] \times [\alpha_1,\alpha_2]$ is a minimal
surface bounded by $\rho \mapsto (\rho,\alpha_1,u_n(\rho,\alpha_1))$
for $0\le \rho\le \rho_0$, $\theta\mapsto
(\rho_0,\theta,u_n(\rho_0,\theta))$ for $\alpha_1\le \theta\le
\alpha_2$, $\rho\mapsto (\rho,\alpha_2,u_n(\rho,\alpha_2))$ for
$\rho_0\ge \rho\ge 0$ and the segment $\{\boO\} \times 
      [b_n,a_n]$. $\Sigma_n$ is a minimal surface of $\boD \times \R$ 
      which is of the type of the disk then we can parametrize
      conformally it 
by the upper half-disk $\Delta^+=\{\zeta\in\Delta| \Im(\zeta)\ge
0\}$: there exist an harmonic map $X_n:\Delta^+\rightarrow \boD$ and
an harmonic function $z_n:\Delta^+\rightarrow \R$ such that the
surface $\Sigma_n$ is $\{(X_n(\zeta),z_n(\zeta)),\zeta\in
\Delta^+\}$. We choose $X_n$ and $z_n$ such that
$(X_n(-1),z_n(-1))=(0,0,a_n)$, $(X_n(0),z_n(0))=(0,0,c_n)$ and
$(X_n(1),z_n(1))=(0,0,b_n)$. We want to apply Theorem \ref{comp} to
the sequence $(X_n,z_n)|_{\partial\Delta^+}$, but we can not do this
directly. 

\begin{prop}{\label{convergence}}
Let $(X_n,z_n)$ be as above then , if $k$ is big enough, a subsequence
$(X_{n'})$
converges uniformly to an harmonic map $X: \Delta^+\rightarrow \boD$
and $(z_{n'})$ converges uniformly on each compact subset to an harmonic
function 
$z:\inter{\Delta^+} \rightarrow \R$ (where $\inter{\Delta^+}$ is the
interior of $\Delta^+$). On $\inter{\Delta^+}$, $X$ and $z$ satisfy
$z=u\circ X$.  
\end{prop}

\begin{proof}
In fact the function $z_n$ is not the good function to consider. On
$U=[0,\rho_0] \times [\alpha_1,\alpha_2]$, the function $\Psi_{u_n}$ is
defined with $\Psi_{u_n}(0,0)=0$. $\Psi_{u_n}$ corresponds to the harmonic
conjugate 
of $z_n$. More precisely, if $\zeta_0\in\inter{\Delta^+}$ and
$z_n^*$ is the conjugate function to $z_n$ such that
$z_n^*(\zeta_0)=\Psi_{u_n}(X_n(\zeta_0))$, we have $z_n^*=\Psi_{u_n}
\circ X_n$. \emph{A priori}, the preceding equality is true only in the
interior of $\Delta^+$, but, since $\Psi_{u_n}$ is Lipschitz
continuous, $z_n^*$ can be extended to the boundary such that the
equality is true everywhere. 

By a result of J.C.C.~Nitsche \cite{Ni}, the area of $\Sigma_n$ is bounded by
$Area(U)+ \int_{\partial U}|u_n|$. Since $u_n$ converges to $u$ uniformly on
each compact subset of $\Ome\backslash \{\boQ\}$, the 
$u_n$ are uniformly bounded functions on $U$ by Proposition \ref{major}. This
proves that the areas of 
$\Sigma_n$  are uniformly bounded by a constant $M$. Since $(X_n,z_n)$
are conformal:
$$\int_{\Delta^+}|{X_n}_x|^2+ |{X_n}_y|^2+|\nabla z_n|^2\dd x\dd y=
Area(\Sigma_n)^2< M^2$$

Let $F_n:\Delta^+\rightarrow \boD\times\R^3$ be the map defined by
$F_n(\zeta)=(X_n(\zeta),z_n^*(\zeta),x,y)$ with $\zeta=x+iy$. We shall
apply Theorem \ref{comp} with $\gamma_n$ the restriction of $F_n$ to the
boundary of 
$\Delta^+$, $f_n=z_n$ and with $\zeta_1=-1$, $\zeta_2=0$ and
$\zeta_3=1$. We note $G_n=(F_n,z_n)$. Since $z_n^*$ is conjugate to
$z_n$, $|\nabla z_n^*|=|\nabla z_n|$; then 
$$\int_{\Delta^+}\left(|{G_n}_x|^2+|{G_n}_y|^2\right) \dd x\dd y\le
2M^2+2 Area(\Delta^+)$$ 
We have the third hypothesis of Theorem \ref{comp}.

Because of the function $x$ in $F_n$, $d(F_n(\zeta_i),F_n(\zeta_j))\ge
1$ for $i\neq j$, this is the first hypothesis of Theorem \ref{comp}.

To prove that $(\gamma_n,f_n)$ satisfies the second hypothesis of Theorem
\ref{comp}, we need a lemma.
\begin{lem}{\label{techni}}
If $k$ is big enough, for big $n$ we have:
\begin{gather}
1\le \inf_{(0,\rho_0]\times[-2kq\pi,-2(k-1)q\pi]}u_n-
  \sup_{(0,\rho_0]\times [0,2q\pi]}u_n \label{eq1}\\  
1\le \inf_{(0,\rho_0]\times[0,2q\pi]}u_n- \sup_{(0,\rho_0]\times 
  [2kq\pi,2(k+1)q\pi]}u_n\label{eq2}  
\end{gather}
\end{lem}
\begin{proof}
There exists a constant $c\in\R$ such that $u\circ f= u+c$, by
Proposition \ref{cno=0} $c\neq 0$ and by construction $c<0$, this is due to
the value $+$ and $-\infty$ on $\widetilde{\boL}^+$ and $\widetilde{\boL}^-$
for the function $u_n$. This implies
that 
\begin{gather*}
\sup_{(0,\rho_0]\times[2lq\pi,2(l+1)q\pi]}u- \sup_{(0,\rho_0] \times
  [2mq\pi,2(m+1)q\pi]}u= c(l-m)\\
\inf_{(0,\rho_0]\times[2lq\pi,2(l+1)q\pi]}u- \inf_{(0,\rho_0] \times
  [2mq\pi,2(m+1)q\pi]}u= c(l-m)
\end{gather*}

This implies that for a $k$ big enough:
\begin{gather*}
2\le \inf_{(0,\rho_0] \times [-2(k+1)q\pi,-2(k-2)q\pi]}u-
  \sup_{(0,\rho_0] \times [-2q\pi,4q\pi]}u\\  
2\le \inf_{(0,\rho_0] \times [-2q\pi,4q\pi]}u- \sup_{(0,\rho_0]\times 
  [2(k-1)q\pi,2(k+2)q\pi]} u
\end{gather*}

Let us apply Proposition \ref{major} with $\epsilon=1/4$, there exists then
$d$ such that for every $l\in \Z$ and every $n>l+2$:
\begin{gather*}
\sup_{(0,\rho_0] \times [2lq\pi,2(l+1)q\pi]}u_n \le \sup_{[d,\rho_0]
    \times [2(l-1)q\pi,2(l+2)q\pi]}u_n +\frac{1}{4}\\
\inf_{(0,\rho_0] \times [2lq\pi,2(l+1)q\pi]}u_n \ge \inf_{[d,\rho_0]
    \times [2(l-1)q\pi,2(l+2)q\pi]}u_n -\frac{1}{4}\\
\end{gather*}

But since $u_n\rightarrow u$ uniformly on every compact subset of
$\Ome\backslash \{\boQ\}$ we have for big $n$:
\begin{gather*}
\sup_{[d,\rho_0] \times [2(l-1)q\pi,2(l+2)q\pi]}u_n \le \sup_{[d,\rho_0]
    \times [2(l-1)q\pi,2(l+2)q\pi]}u +\frac{1}{4}\\
\inf_{[d,\rho_0] \times [2(l-1)q\pi,2(l+2)q\pi]}u_n \ge \inf_{[d,\rho_0]
    \times [2(l-1)q\pi,2(l+2)q\pi]}u -\frac{1}{4}\\
\end{gather*}

Then in using these inequalities with $l=0$ and $l=-k$, we have for
$n$ big enough:
\begin{equation*}
\begin{split}
\inf_{(0,\rho_0]\times[-2kq\pi,-2(k-1)q\pi]}u_n &\ge 
\inf_{[d,\rho_0] \times [-2(k+1)q\pi,-2(k-2)q\pi]}u_n -\frac{1}{4}\\ 
&\ge \inf_{[d,\rho_0] \times [-2(k+1)q\pi,-2(k-2)q\pi]}u -\frac{1}{2}\\
&\ge \sup_{[d,\rho_0] \times [-2q\pi,4q\pi]}u +\frac{3}{2}\\
&\ge \sup_{[d,\rho_0] \times [-2q\pi,4q\pi]}u_n +\frac{5}{4}\\
&\ge \sup_{(0,\rho_0] \times [0,2q\pi]}u_n +1
\end{split}
\end{equation*}
This shows \eqref{eq1}. With $l=0$ and $l=k$, we get \eqref{eq2}:
$$\inf_{(0,\rho_0]\times[0,2q\pi]}u_n \ge \sup_{(0,\rho_0]\times 
  [2kq\pi,2(k+1)q\pi]}u_n +1$$
\end{proof}

Let $k\in\N$ given by Lemma \ref{techni}. In $\boD$, the curves $\zeta\mapsto
X_n(\zeta)$ for $\zeta\in\partial \Delta^+$ have the same image
$\Gamma$ for every $n$. $\Gamma$ is 
a Jordan curve in $\boD$ so for every $m>0$ there exists $\epsilon>0$ such
that, if $p',p''\in \Gamma$ and $d(p',p'')<\epsilon$, one of the components
of $\Gamma \backslash \{p',p''\}$ is of diameter less than $m$. Let
$\delta$ be $\min\{\rho_0/2,1/2\}$ if $\boA$ is the point of coordinates
$(\rho_0,\alpha_1)$ then $2\delta$ is less than the distance between
$\boO$ and 
$\boA$ and if $m<\delta$, in the above property, there is only one
component of $\Gamma \backslash \{p',p''\}$ with diameter less than
$m$. 

Let $0<m<\delta$, there exists $\epsilon>0$ that satisfies the above
property, there also exists $\eta$ such that if $\zeta',\zeta''\in\partial
\Delta^+$ and $|\zeta'-\zeta''|<\eta$, one of the components of
$\partial \Delta^+\backslash \{\zeta',\zeta''\}$ is of diameter less
than $m$. Since $m<\delta\le 1/2$ this component is unique. Let
$n\in\N$ be big enough such that \eqref{eq1} and \eqref{eq2} are
satisfied. Let $\zeta'$ and 
$\zeta''$ in $\partial \Delta^+$ such that the distance between 
$G_n(\zeta')$ and $G_n(\zeta'')$ is less than $\min\{\epsilon,\eta,
1/2\}$. Let us note $F_n(\zeta')=(p',\Psi_{u_n}(p'),x',y')$ and
$F_n(\zeta'')=(p'',\Psi_{u_n}(p''),x'',y'')$. We have $d(p',p'')\le 
\epsilon$, then there exists one component of $\Gamma \backslash
\{p',p''\}$ with diameter less than $m$. Let $I$ be the part of
$\partial \Delta^+$ that parametrize this component; the end points
of $I$ are $\zeta'=x'+iy'$ and $\zeta''=x''+iy''$. Since
$|\zeta'-\zeta''|<\eta$, $I$ or its complemantary in $\partial \Delta^+$
is of diameter less than $m$; let us prove that it is $I$. If $\boO
\notin X_n(I)$ then $\{\zeta\in \partial\Delta^+|\zeta\in\R\}\cap I=
\emptyset$, but $\{\zeta\in \partial\Delta^+|\zeta\in\R\}$ is of
diameter $2$ then $I$ is of diameter less than $m$. If $\boO\in
X_n(I)$ then $\boA\notin X_n(I)$ and the points $p'$ and $p''$ are at a
distance less than $m$ from $\boO$ then the point $p'$ and $p''$ can not
be on the part of $\Gamma$ : $\theta\mapsto (\rho_0,\theta)$. There are
different possible cases. First, $I$ can be included in $\{\zeta\in
\partial\Delta^+|\zeta\in\R\}$ then since $|\zeta'-\zeta''|<\eta$, $I$
is of diameter less than $m$. The second case is when $\{\zeta\in
\partial\Delta^+|\zeta\in\R\}$ is included in $I$ then we can suppose
that $p'\in\{(\rho,\alpha_1)\}_{0\le \rho\le \rho_0}$ and
$p''\in\{(\rho,\alpha_2)\}_{0\le \rho\le \rho_0}$  this implies, by Lemma
\ref{techni}, 
that $|z_n(\zeta')-z_n(\zeta'')|\ge 2$ which is impossible since the
distance between $G_n(\zeta')$ and $G_n(\zeta'')$ is less than
$\min\{\epsilon,\eta, 1/2\}<2$. For the third case, we can suppose that
$\zeta'\in\{\zeta\in \partial\Delta^+|\zeta\in\R\}$ and
$p''\in\{(\rho,\alpha_2)\}_{0\le \rho\le \rho_0}$ (all the other cases are
given by permutating $\zeta'$ and $\zeta''$ and $\alpha_1$ and
$\alpha_2$). Since $\boA\notin X_n(I)$, we have $1\in I$. Besides
$z_n$ is decreasing along $\{\zeta\in \partial\Delta^+|\zeta\in\R\}$,
then, since $c_n-z_n(\zeta'')\ge 1$, by Lemma \ref{techni}, and
$|z_n(\zeta')-z_n(\zeta'')|\le 
1/2$, $\zeta'\in[0,1]$. Since $|\zeta'-\zeta''|\le 1/2$, $\zeta''$ is in
the part of $\partial\Delta^+$ with $\Re(\zeta)\ge 0$ then the
shortest component of $\partial\Delta^+\backslash\{\zeta',\zeta''\}$ is
the one that contains $1$ then it is $I$.

So it was proved that $I$ and $X_n(I)$ are of diameter less than
$m$. Since $\Psi_{u_n}$ is $1$-Lipschitz continuous and
$z_n^*=\Psi_{u_n}\circ X_n$ the set $z_n^*(I)$ is of diameter less
than $m$. Then, if $\zeta',\zeta''\in\partial\Delta^+$ are such that the
distance between $G_n(\zeta')$ and $G_n(\zeta'')$ is less than
$\min\{\epsilon,\eta, 1/2\}$, one component of $F_n(\partial\Delta^+)
\backslash \{F_n(\zeta'),F_n(\zeta'')\}$ is of diameter less than
$\sqrt{3} m$; this is the second hypothesis.

Let us apply the result of Theorem \ref{comp}, the family $\{
F_n|_{\partial\Delta^+} \}$ is equicontinuous then there exists a
subsequence $(F_{n'}|_{\partial\Delta^+})$ that converges
uniformly. Since the $F_n$ are harmonic maps, the sequence $(F_{n'})$
converges uniformly on $\Delta^+$. Then $X_{n'}\rightarrow X$ uniformly
on $\Delta^+$ and $z_{n'}^*\rightarrow z^*$ with $X$ and $z^*$ harmonic
map. The fact that $z_n$ and $z_n^*$ are harmonic conjugates implies that
$(z_{n'})$ 
converges to the harmonic conjugate $z$ of $z^*$ and the convergence is
uniform on each compact subset of $\inter{\Delta^+}$. The equality $z=u\circ
X$ is just the limit of the equality $z_n=u_n\circ X_n$. This ends the proof
of Proposition \ref{convergence}.
\end{proof}

\begin{cons}{alpha}
Let us summarize what we have done above. We take $\alpha\in \R$ and
$k\in \N^*$. For every $n\in \N$, we consider $(X_n,z_n)$ a conformal
parametrization of the graph of $u_n$ over $[0,\rho_0]\times
[\alpha-2kq\pi, \alpha+2kq\pi]$ by the upper-half disk $\Delta^+$ such
that the part parametrized by $\Delta^+ \cap \R$ is the vertical
segment in the boundary which is above $\boO$. Then Proposition
\ref{convergence} says us that if we take $k$ big enough we have a
subsequence $(X_{n'},z_{n'})$ that converges to $(X,z)$ as described
in Proposition \ref{convergence}. We also have the sequence $(z_n^*)$ and
$z_{n'}^*\rightarrow z^*$ uniformly on $\Delta^+$. $X$ satisties t$\widetilde{\boL}^+$he following
property. 
\end{cons}

\begin{cor}
Let $X$ be constructed as in Construction \ref{alpha}. Then
$$X(\Delta^+)=[0,r_0] \times [\alpha_1,\alpha_2]$$ 
\end{cor}
\begin{proof}
Let $\boA\in [0,r_0] \times [\alpha_1,\alpha_2]$ then for every $n$
there exists $\zeta_n$ such that $X_n(\zeta_n)=\boA$. We have
$X_{n'}\rightarrow X$ uniformly on $\Delta^+$ and, since $\Delta^+$ is
compact, a subsequence $n''$ of $n'$ is such that
$\zeta_{n''}\rightarrow\zeta\in \Delta^+$. Then $\boA=X_{n''}
(\zeta_{n''}) \rightarrow X(\zeta)$ and $\boA=X(\zeta)$.
\end{proof}

\subsubsection{The proof}
Before proving Theorem \ref{reg2}, we need a remark
\begin{prop}
Let $\alpha\in\R$ then $u(\rho,\alpha)$ converges as $\rho$ goes to $0$.  
\end{prop}
\begin{proof}
Let $\alpha\in\R$ and $k\in\N$ be big enough such that Construction
\ref{alpha} can be done. By construction, we have
$X_n$, $X$, $z_n^*$ and $z^*$ defined on $\Delta^+$. Let $U$ be a
compact neighborhood of $0$ in $\Delta^+$ such that $X(U)\subset
[0,\rho_0]\times 
(\alpha-\pi/2,\alpha+\pi/2)$. Using the projection $H:(\rho,\theta)
\mapsto (\rho\cos\theta,\rho\sin\theta)$ we can consider $X$ as a
harmonic map from $\Delta^+$ to $\R^2$. For $n$ big enough,
$X_n(U)\subset [0,r_0]\times (\alpha-\pi/2,\alpha+\pi/2)$, the maps
$X_n$, $X$, $z_n^*$ and $z^*$ can be extended to $U'$ where $U'$ is
the union of $U$ and $\{\zeta\in\C|\ \bar{\zeta}\in U\}$. $X_n$ and
$z_n^*$ are extended by Schwarz reflection principle to harmonic
map. Since $X_{n'}\rightarrow X$ and $z_{n'}^*\rightarrow z^*$ on $U$
we have the same convergence for their extensions to $U'$. In the same
way, the functions $z_n$ and $z$ can be extended to $U'$ and the
convergence $z_{n'}\rightarrow z$ is uniform on each compact of
$U'$. The map $(X,z)$ on $U'$ gives a minimal surface of $\R^3$. Since
the point $(X_n(0),z_n(0))$ is for every $n\in\N$ the limit point of
$(\rho, \alpha, u_n(\rho,\alpha))$ as $\rho$ goes to $0$. The normal to
the minimal surface $(X_n,z_n)$ at the origin is
$(\sin\alpha,-\cos\alpha,0)$. Then, the normal to the minimal surface
$(X,z)$ is also $(\sin\alpha,-\cos\alpha,0)$ at the origin. Then the
intersection of the minimal surface $\{X(\zeta),z(\zeta)\}_{\zeta\in
  U}$ with the vertical plane of equation $x\sin\alpha-y\cos\alpha=0$
in the neighborhood of the point $(X(0),z(0))$ is composed of a piece
of the vertical straight line passing by this point and a curve which
is above $\{(\rho,\alpha)\}_{\rho>0}$ with $(X(0),z(0))$ as end
point. Then this curve is, in fact, $\{\rho,\alpha,
u(\rho,\alpha)\}_{\rho>0}$, this implies that $\dis\lim_{\rho\rightarrow 0}
u(\rho,\alpha)=z(0)$.
\end{proof}

We can then make the proof of Theorem \ref{reg2}
\begin{proof}[Proof of Theorem \ref{reg2}]. 
Let $k\in\N$ be big enough such that we can make Construction
\ref{alpha} with $\alpha=0$. We then get a map $X:\Delta^+\rightarrow
\boN$ and two functions $z:\inter{\Delta^+} \rightarrow \R$ and
$z^*:\Delta^+\rightarrow \R$. Let $I$ be 
the part of $\partial\Delta^+$ such that $\zeta\in I
\Leftrightarrow X(\zeta)=\boO$. Since $X$ is monotone, $I$ is connected
and by construction $\{\zeta\in \partial\Delta^+|\zeta\in\R\} \subset
I$. Let $\zeta_1$ and $\zeta_2$ denote the two end-points of $I$ and
consider the biholomorphic map $h:\Delta^+\rightarrow \Delta^+$  such
that $h(-1)=\zeta_1$, $h(0)=(0)$ and $h(1)=\zeta_2$. If
$\widetilde{X}=X\circ h$, $\widetilde{z}=z\circ h$ and
$\widetilde{z}^*=z^*\circ h$ (it is obvious that the conjugate
harmonic function to $\widetilde{z}$ is $\widetilde{z}^*$), we make
only a reparametrization of the minimal surface $(X,z)$. For
$\zeta\in \partial\Delta^+$, $\widetilde{X}(\zeta)=\boO \Leftrightarrow
\zeta\in\R$. Then let consider $\zeta_0\in \partial\Delta^+ \backslash
\R$. We have $\widetilde{X}(\zeta_0)\neq\boO$ then since
$\widetilde{z}(\zeta) =u\circ\widetilde{ X} (\zeta)$ we can define
$\widetilde{z} (\zeta_0)$ by making $\zeta$ converging to
$\zeta_0$. $\widetilde{z}$ can be also defined at $-1$ and $1$: as
$\zeta\in\partial\Delta^+ \backslash \R$ goes to $-1$ or $1$,
$X(\zeta)$ goes to $\boO$ along $\{(\rho,\alpha_1)\}$ or $\{(\rho,
\alpha_2)\}$ and then $\widetilde{z}(\zeta)$ goes to
$\dis\lim_{\rho\rightarrow 0} u(\rho,\alpha_1)=z(-1)$ or
$\dis\lim_{\rho\rightarrow 0} u(\rho,\alpha_2)=z(1)$. Since
$\widetilde{X}=\boO$ on $\R\cap \Delta^+$, $\widetilde{z}^*=0$ on the
same set. then we can extend $\widetilde{z}^*$ to the whole disk by
Schwartz reflection principle, this prove that $\widetilde{z}$ can be
also extended to the interior of the whole disk by reflection. Since
we have define $\widetilde{z}$ on the circular part of the boundary of
$\Delta^+$, $\widetilde{z}$ is then defined also on the boundary of
the disk and is continuous on the boundary then $\widetilde{z}$ is the
harmonic extention to the disk of this countinuous function on the
circle. This proves that $\widetilde{z}$ is a continuous function on
$\Delta^+$. Then we have $(\widetilde{X},\widetilde{z})_{\Delta^+}$
which is a parametrization of the graph of $u$ above $[0,\rho_0] \times
[\alpha_1,\alpha_2]$ (because $\widetilde{z}= u\circ\widetilde{X}$) and
$(\widetilde{X},\widetilde{z})_{\Delta^+}$ has a boundary such that the
part of this boundary which is above $\boO$ is a vertical
segment. This proves Theorem \ref{reg2} because $[0,\rho_0] \times
[\alpha_1,\alpha_2]$ contains several period of the domain $\boN$.  
\end{proof}

\subsection{Property of the solution 3}

Let $u$ be a solution of the Dirichlet problem asked in Theorem
\ref{graph}. We then can understand the boundary behaviour of the
graph of $u$ in $\R^3$.

From Theorem \ref{reg1}, we know that, over the neighborhood $\boN$
of $\boQ$, the graph is bounded by the vertical straight-line passing
by the point $\psi(Q)$.

The other points where there are boundary components for the graph of
$u$ are the vertices of $\Ome$. This points statisfy the hypotheses of
Theorem $3$ in  \cite{Ma1}. Then the graph of $u$ is bounded by vertical
straight-lines near the vertices.

The last remark we can make is the following. We have $\Psi_u(\boQ)=
0= \Psi_u(\boV)$; this implies that the conjugate surface of the graph
of $u$, which is bounded by the conjugate of the boundary of the graph
of $u$, has its boundary included in the plane $\{z=0\}$. 

This remark implies that, if $\Sigma$ is the graph of $u$ over
$\Ome_0^1$, a period of $\Ome$, its conjugate surface can be extend by
symmetry with respect to the plane $\{z=0\}$, we note $\Sigma^*$ this
symmetric surface. We then have the following result.

\begin{lem}{\label{totcurv}}
$\Sigma^*$ is of finite total curvature and its total curvature is
$4\pi r$.
\end{lem}

\begin{proof}
Using arguements given in \cite{Ma1}, we can prove that $\Sigma$ is of
finite total curvature, this implies that the same is true for
$\Sigma^*$. Then, as in Proposition 2.2 in \cite{HK}, each catenoidal end
gives a contribution of $2\pi$ to the total curvature (see \cite{JM}, for the
original arguements) and, using Gauss-Bonnet Theorem, we compute the value of
the total curvature and get $4\pi r$. 
\end{proof}


\section{The period problem{\label{peri}}}
\subsection{The general case}{\label{periodproblem}}

We now try to build a $r$-noid with genus $1$ and horizontal ends of type I
for a given polygon of flux.

Let $V=(v_1,\dots,v_r)$ be a polygon that bounds a multi-domain with
cone singularity $(D,Q,\psi)$. Using Construction \ref{C1} and
Construction \ref{C2} as in Section \ref{exist}, we get a multi-domain
with logarithmic singularity $(\Ome,\boQ,\phi)$ and a solution $u$ on
$\Ome$ of the Dirichlet problem asked in Theorem \ref{graph}; as in
Section \ref{exist}, we assume that the period $2q\pi$ of $\Ome$ is
the angle at the cone singularity of $D$. Let us
consider $\Ome_0^1$ one period of the multi-domain $\Ome$ and $\Sigma$ the
graph of $u$ over $\Ome_0^1$. We know, because of the result of the
preceding sections, that $\Sigma$ is a minimal surface bounded by $r-1$
vertical lines passing by the vertices $P_i$ ($i\neq 1$) of the polygon $V$,
two vertical half-lines over $P_1$, a
vertical segment over $\psi(Q)$ and two curves over the segment
$[\psi(Q),\psi(P_1)]$. The conjugate surface of $\Sigma$ is included in
$\{z\ge 0\}$ and the conjugates of the $r-1$ vertical lines, the two vertical
falf-lines and the 
vertical segment are exactly the intersection of the conjugate surface
with the plane $\{z=0\}$. We then can extend the conjugate surface
by symmetry with respect to this plane, we get a new surface that we note
$\Sigma^*$. $\Sigma^*$ is a solution for the Plateau problem at
infinity for the data $V$ if the two components of boundary of
$\Sigma^*$, coming, by conjugation, from the two curves which are over
$[\psi(Q),\psi(P_1)]$, glue together such a way that $\Sigma^*$ has no
boundary. 

In fact this two components of boundary differ from a translation, how
can we compute the vector of the translation? 
From the function $u$
we can derive three closed 1-forms $\dd X_1^*$, $\dd X_2^*$ and $\dd X_3^*$
on $\Ome$ which are the differential of the three coordinate
functions of the conjugate surface to the graph of $u$ (These
1-forms depend only on the first derivatives of $u$). For example,
we have $\dd X_3^*=\dd \Psi_{u}$. In $\boN$, the neighborhood of
$\boQ$ in $\Ome$, we consider the path $\Gamma:
\theta\mapsto(\rho,\theta)$ for some $\rho<\rho_0$ and
$\theta\in[0,2q\pi]$; $\Gamma$ is a lift of a generator of
$\pi_1(D \backslash\{Q\})$. Then the two components of boundary of
$\Sigma^*$ differ from the following vector, called the \emph{period
vector}:
$$\left(\int_\Gamma \dd X_1^*,\ \int_\Gamma \dd X_2^*,\ \int_\Gamma
\dd X_3^*\right)$$
Since $\dd X_3^*=\dd \Psi_{u}$ and as $\Psi_{u}$ is invariant by
$f$, $\dis \int_\Gamma \dd X_3^*=0$. Obviously, the value of the
integrals is the same for every $\Gamma$ which is the lift of a
generator of $\pi_1(D \backslash\{Q\})$; in fact, from Remark \ref{quotient1}, 
since $\dd X_i^*$ depends only on the derivatives of $u$, $\dd X_i^*$
is well defined on $D \backslash Q$

Then the question of the existence of a solution to the Plateau
problem at infinity for the data $V$ becomes: knowing if there exists
$(D,Q,\psi)$ bounded by $V$ such that the associated period vector is
zero; this is the period problem. 

\begin{rem}{doublement}
We now give some explanations on Remark \ref{angle} and the hypothesis on $D$
saying that its angle at the cone singularity is the period of $\Ome$. Let
$V=(v_1,...,v_r)$, $(D,Q,\psi)$, $(\Ome,\boQ,\phi)$ and $u$ be as above ($D$
sastisfies the hypothesis H). We also suppose the the period vector associated
to $D$ is zero. Then $\Sigma^*$ which is the conjugate surface to $\Sigma$, the
graph of $u$ over $\Ome_0^1$, extended by symmetry is a $r$-noid with genus
$1$ and horizontal ends of type $I$ having $V$ as flux polygon. Let $f$ be the
isometry of $\Ome$ associated to its periodicity. Let $a\in\N^*$, then the
quotient of $\boW$ by the group $\{f^{an}\}_{n\in\Z}$ ($\boW$ is given by
Construction \ref{C1} applied to $D$) is a multi-domain with cone singularity
that bounds the polygon
$$V_a=(\underbrace{ v_1,\dots,v_r,\dots \dots,v_1,\dots,v_r}_{a\
  \textrm{times}})$$ 

Besides, if $\Sigma_a$ is the graph of $u$ on $\Ome_0^a$, the conjugate
surface $\Sigma_a^*$ of $\Sigma_a$ extended by symmetry is a $ar$-noid with
genus $1$ and horizontal ends of type $I$ having $V_a$ as flux polygon. In fact
$\Sigma_a^*$ is just $\Sigma^*$ that we cover $a$ times. Then we also can find
solution to the Plateau problem at infinity for $D$ that does not satisfy the
hypothesis H.
\end{rem}

\subsection{The period map and the proof of Theorem 
\ref{main}{\label{periodmap}}}

In this subsection we explain how we shall prove Theorem \ref{main}.

Let $V=(v_1,\dots,v_r)$ be a polygon bounded by an immersed polygonal
disk $(\boP,\psi)$. For each $A$ in the interior of $\boP$, Construction
\ref{C1} gives 
us 
$(\boW_A,\boA,\phi_A)$ a multidomain with a logarithmic singularity
. Then by applying Construction \ref{C2} and Theorem
\ref{graph}, we get a periodic multi-domain with logarithmic singularity
$(\Ome_A,\boA,\phi_A)$ and a function $u_A$ 
defined on $\Ome_A$. We then have the three closed $1$-forms $\dd
X_i^*(A)$. To prove Theorem \ref{main}, we need to find a 
point $A\in\inter{\boP}$ such that the period vector associated to the above
constuction for $A$ is $0$. 

In fact, we have a map from the interior $\boP$ to $\R^2$ which associates to
every point $A\in\inter{\boP}$ 
the vector $\dis \left(\int_\Gamma \dd X_1^*(A),\ \int_\Gamma \dd X_2^*(A) 
\right)$, this map is the period map and will be noted $Per$. Then to
prove Theorem \ref{main}, the problem is to prove that this map
vanishes at one point of $\inter{\boP}$. The period map satisfies the following
proposition.

\begin{prop}
The period map $Per$ is continuous on the interrior of $\boP$.
\end{prop}

\begin{proof}
Let us consider a sequence $(A_n)$ of points in $\inter{\boP}$ that converges
to 
$A\in\inter{\boP}$. Let $\Gamma$ be a closed path in $\boP\backslash\left\{A,
  A_0, 
A_1, \dots, A_n, \dots\right\}$ such that for every $n$, $\Gamma$ is a
generator of $\pi_1(\boP \backslash\{A_n\})$ and $\Gamma$ is a
generator of $\pi_1(\boP \backslash\{A\})$. Since for $A$ (or every
$A_n$) we have $u_A\circ f=u_A+c$ the derivatives of $u_A$ are well
defined on $\boP\backslash \{A\}$. Then the two closed $1$-forms $\dd
X_1^*(A)$ and $\dd X_2^*(A)$ are well defined on $\boP\backslash
\{A\}$. The same is true for $A_n$. We then have:
$$Per(A_n)=\left(\int_\Gamma \dd X_1^*(A_n),\ \int_\Gamma \dd X_2^*(A_n)
\right)$$

On $\boP$, we have the sequence of the derivatives of $u_{A_n}$ and
these derivatives converge to the derivatives of $u_A$ if there is no
line of divergence (a line of divergence is a phenomenon linked to the
behaviour of the first 
derivatives so we can use this arguement in this case). Since the
arguements used in the proof of the existence part of Theorem
\ref{graph} are always true, there is no line of divergence and $\dd
X_1^*(A_n) \rightarrow \dd X_1^*(A)$ and  $\dd X_2^*(A_n) \rightarrow
\dd X_2^*(A)$, the convergence is uniform along $\Gamma$. This proves that
$Per(A_n) \rightarrow Per (A)$ by integration. 
\end{proof}

The idea to prove that the period map vanishes at one point in the interior of
$\boP$
is then to extend continuously $Per$ to the boundary of $\boP$ and
show that the degree of the period map along the boundary of $\boP$ is
non zero. In fact, we shall use a modified boundary of $\boP$. Then, using
Proposition 3.20 in \cite{Fu}, this proves that there 
exists a point $A \in \inter{\boP}$ 
where $Per(A)=0$. The following section is devoted to the extension
of $Per$ to the boundary and to the proof of Theorem \ref{degree}
that establishes that the degree of the period map is non-zero along the
boundary.  

To extend the period map on the boundary we make a renormalization of the
map $Per$: let $A$ be a point of $\inter{\boP}$, if $||Per(A)||\le
1$ then we do not change the value of $Per(A)$ but if $||Per(A)||\ge
1$ the new value of $Per(A)$ is $\dis\frac{Per(A)}{||Per(A)||}$. The new
period map is always continuous and for every point the norm of
the period at this point is less than one.


\section{The period map on the boundary of $\boP${\label{extention}}}
We use the notations introduced in Subsection \ref{periodmap}.

\subsection{The behaviour on the edges}

\begin{prop}{\label{edge}}
Let $(A_n)$ be a sequence in $\inter{\boP}$ such that $A_n\rightarrow A$ where
$A$ is a point in the interior of one edge of the boundary of
$\boP$. Then $Per(A_n)$ converges to $\dd\psi|_A(N)$ where $N$ is the
outer unit normal to the edge at $A$, we recall that $\psi$ is the
developping map of $\boP$.
\end{prop}

\begin{proof}
As in \cite{Ma1}, we note $\Ome(\boP)$ the multi-domain obtained when we
glue to every edge $[P_i,P_{i+1}]$ a half strip isometric to
$[P_i,P_{i+1}]\times \R_+$. Then, for every $A_n$ the covering map
$\pi:\boW_{A_n} \rightarrow \boP$ extends to a covering map $\pi:
\Ome_{A_n}\rightarrow \Ome(\boP)$ and the derivatives of $u_{A_n}$ are
then well defined on $\Ome(\boP)\backslash \{A_n\}$, by Remark
\ref{quotient1}. Suppose that the point $A$ is in 
the interior of the edge $[P_1,P_2]$  and that $|P_1P_2|=2$. By
choosing a good chart, we can suppose that $[-1,1]\times \R_+$ is
the half-strip glued to this edge. We note $D_r$ the domain in
$\R^2$ which is the intersection of the domain $y\le 0$ and the disk
of center $(0,r)$ and radius $\sqrt{r^2+1}$. By choosing a $r>0$ big
enough, $D_r$ is a neighborhood of the edge $[P_1,P_2]$ in $\boP$. 

We suppose that $A$ is the point $(a,0)$ ($-1<a<1$) and the points
$A_n$ lie in $D_r$ and have coordinates $(a_n,b_n)$ ($b_n<0$). We
have $a_n\rightarrow a$ and $b_n\rightarrow 0$. For every $n$, we note
$L_n$ the half straight-line $\{(a_n,b_n+t)\}_{t\ge 0}$ and $L$ the
half straight-line $\{(a,t)\}_{t\ge 0}$. By using the
covering map $\pi$ and $u_{A_n}$, we can define on $\Ome(\boP)
\backslash L_n$ a function $u_n$ which has the same derivatives as
$u_{A_n}$; $u_n$ is solution of the minimal surface equation and has
the value $+\infty$ (resp. $-\infty$) on $\pi(\boL_i^+)$
(resp. $\pi(\boL_i^-)$). We then want to understand the convergence of
$u_n$. By the same arguement as in the proof of Theorem \ref{graph},
there is no line of divergence in $\Ome(\boP) \backslash L$; this
proves that, for a subsequence, $(u_n)$ converges to a solution $u$ of
\eqref{MSE} on $\Ome(\boP) \backslash L$. By Lemma \ref{lemdiv2}, $u$
takes 
the value $+\infty$ (resp. $-\infty$) on $\pi(\boL_i^+)$
(resp. $\pi(\boL_i^-)$). We then need to understand its behaviour near
$L$ to know the function $u$.

Since by convention $\Psi_{u_n}(A_n)=0$, we fix $\Psi_u(A)=0$
\begin{lem}{\label{psi=t}}
With this convention, $\Psi_u(a,t)=t$ for $t\ge 0$.
\end{lem}

\begin{proof}
Since $\Psi_u$ is $1$-Lipschitz continuous, $\Psi_u(a,t)\le t$; let us
suppose that for some $t_0$ we have $\Psi_u(a,t_0)= t_0-\epsilon$ with
$\epsilon >0$, then for $t>t_0$ $\Psi_u(a,t)\le t-\epsilon$. We have
$\Psi_u(1,t)=t$ for every $t\ge 0$ because $u$ takes the value
$+\infty$ along $\pi(\boL_0^+)$. We have:
$$\Psi_u(1,t)- \Psi_u(a,t)= \lim_{n\rightarrow+\infty} \Psi_{u_n}(1,t)-
\Psi_{u_n}(a,t) =\lim_{n\rightarrow+\infty} \int_{[(a,t),(1,t)]}
\dd\Psi_{u_n}$$ 
By Lemma \ref{jenkins}, the integral is always less than
$2\sqrt{2}\frac{1-a}{t}$ for big $t$. So for $t$ big enough, this
upper-bound is less than $\epsilon$: this give us a contradiction.
\end{proof}

The result of Lemma \ref{psi=t} says us, by Lemma \ref{lemdiv2}, that $u$
takes the value $+\infty$ 
on one side of $L$ and $-\infty$ on the other side; more pricisely,
$u(a+\eta,t)$ tends to $+\infty$ (resp. $-\infty$) if $\eta$ tends to
$0$ by negative value (resp. positive value). There is only one solution for
the Dirichlet problem for such 
boundary condition (we apply Theorem $7$ in \cite{Ma1} with the polygon
$(\overrightarrow{P_1A},\overrightarrow{AP_2},v_2,\cdots,v_r)$), this proves
that the limit for the subsequences of 
$(u_n)$ is unique then the sequence $(u_n)$ must converge to the function $u$.

In $[-1,1]\times \R_+\cup D_r$, the $1$-form $\dd X_1^*(A_n)$ and $\dd
X_2^*(A_n)$ are given by:
\begin{gather*}
\dd X_1^*(A_n)= \frac{q_np_n}{W_n}\dd x+ \frac{1+q_n^2}{W_n}\dd y\\
\dd X_2^*(A_n)= -\frac{1+p_n^2}{W_n}\dd x- \frac{p_nq_n}{W_n}\dd y 
\end{gather*}
with $p_n$ and $q_n$ the first derivatives of $u_n$ (see \cite{Os}). Using
this expressions, we can also define the $1$-forms $\dd X_1^*(A)$ and $\dd
X_2^*(A)$. 

Let $\eta_1$ be a small positive number, $\eta_2<\eta_1$ and $l$ 
positive numbers. Let $\Gamma$ be the closed path which consists in the
segment $[(a+\eta_1,l), (a-\eta_1,l)]$, the segment
$[(a-\eta_1,l),(a-\eta_1,0)]$, the half circle in $D_r$ of center $A$
and radius $\eta_1$ and the segment
$[(a+\eta_1,0),(a+\eta_1,l)]$. Besides we call $\Gamma_1$ the part of
$\Gamma$ consisting in the two vertical segments and the half circle,
$\Gamma_2$ the union of the two segments $[(a+\eta_1,l),(a+\eta_2,l)]$
and $[(a-\eta_2,l), (a-\eta_1,l)]$ and $\Gamma_3$ the segment
$[(a+\eta_2,l), (a-\eta_2,l)]$ (see Figure \ref{edge2}). There exists
$n_0\in\N$ such that for 
$n\ge n_0$ and every $\eta_1$, $\eta_2$ and $l$, the period for the
point $A_n$ is computed by
$$\left(\int_\Gamma \dd X_1^*(A_n), \int_\Gamma \dd X_2^*(A_n)
\right)$$ 

\begin{figure}[h]
\begin{center}
\resizebox{0.8\linewidth}{!}{\input{figgenr15.pstex_t}}
\caption{\label{edge2}}
\end{center}
\end{figure}

Let $0<\alpha<1$. By Lemma \ref{jenkins}, we have $\dis\frac{|q_n|}{|p_n|}\le
\sqrt{2}\frac{2}{l} \frac{1}{1-\frac{4}{l^2}}$ on $\Gamma_2\cup
\Gamma_3$. Then in choosing $l$ big enough, we can ensure that:
$$ \left|\int_{\Gamma_2+\Gamma_3} \dd X_1^*(A_n)\right| <\frac{\alpha}{2}
\int_{\Gamma_2+ \Gamma_3}\dd X_2^*(A_n)$$ 

We have:
\begin{gather*}
\lim_{n\rightarrow +\infty} \int_{\Gamma_1}\dd X_1^*(A_n)=
\int_{\Gamma_1} \dd X_1^*(A)\\
\lim_{n\rightarrow +\infty} \int_{\Gamma_1}\dd X_2^*(A_n)=
\int_{\Gamma_1} \dd X_2^*(A)\\
\end{gather*}
the same is true on $\Gamma_2$ and as $\eta_2$ tends to $0$, we have
$\dis\int_{\Gamma_2} \dd X_2^*(A)\rightarrow+\infty$. The last
assertion is due to Lemma $1$ in \cite{JS} which implies that, as $\eta$ goes
to zero $\dis\frac{p}{W}(a+\eta,l)\longrightarrow 1$ and
$\dis p(a+\eta,l)\ge \frac{C}\eta$ for some constant $C$. We then can
choose $\eta_2$ such that for big $n$, we have
\begin{equation*}
\frac{\left| \int_{\Gamma_1} \dd X_1^*(A_n) \right|}{\int_{\Gamma_2}
  \dd X_2^*(A_n)}< \frac{\alpha}{2}\quad\textrm{and}\quad \frac{\left|
    \int_{\Gamma_1} \dd X_2^*(A_n) \right|}{\int_{\Gamma_2} \dd
  X_2^*(A_n)}< \frac{\alpha}{8}
\end{equation*}

This implies first that $\dis\lim_{n\rightarrow +\infty} \int_\Gamma
\dd X_2^*(A_n)=+\infty$, since $\dis \int_{\Gamma_3}
\dd X_2^*(A_n)\ge 0$; then the period at $A_n$, for big $n$, is
renormalized and must have a non negative second coordinate. Secondly,
for big $n$, we have:
\begin{equation*}
\begin{split}
\frac{\left| \int_\Gamma \dd X_1^*(A_n) \right|}{\left| \int_\Gamma
    \dd X_2^*(A_n) \right|}&\le \frac{\left| \int_{\Gamma_1} \dd
    X_1^*(A_n) \right|+ \left| \int_{\Gamma_2} \dd X_1^*(A_n) \right|+
  \left| \int_{\Gamma_3} \dd X_1^*(A_n) \right|}{-\left|
    \int_{\Gamma_1} \dd X_2^*(A_n) \right|+ \int_{\Gamma_2} \dd
  X_2^*(A_n)+ \int_{\Gamma_3} \dd X_2^*(A_n)}\\
&\le \frac{\frac{\alpha}{2}\int_{\Gamma_2} \dd X_2^*(A_n)+
  \frac{\alpha}{2} \int_{\Gamma_2} \dd X_2^*(A_n)+
  \alpha\int_{\Gamma_3} \dd
  X_2^*(A_n)}{-\frac{\alpha}{8}\int_{\Gamma_2} \dd X_2^*(A_n) +
  \int_{\Gamma_2} \dd X_2^*(A_n)+ (1-\frac{\alpha}{8})\int_{\Gamma_3}
  \dd X_2^*(A_n)}\\
&\le \frac{\alpha \int_{\Gamma_2+\Gamma_3} \dd
  X_2^*(A_n)}{(1-\frac{\alpha}{8}) \int_{\Gamma_2+\Gamma_3} \dd X_2^*(A_n)}\\
&\le \frac{\alpha}{1-\frac{\alpha}{8}} 
\end{split}
\end{equation*}

This proves that the renormalized period converges to the vector
$(0,1)$; this is what we want to prove.
\end{proof}

\subsection{The behaviour at the vertices}

Because of Proposition \ref{edge}, it is clear that we can not extend the
period map to the vertices and obtain a continuous map on the
boundary. In fact the idea to solve this problem is to make a blowing-up
of $\boP$ at its vertices.

Let $P_i$ be a vertex of $\boP$, there exists $\alpha>0$ such that a
neighborhood of $P_i$ is isometric to $\{(\rho,\theta),\ 0\le \rho<\mu,
0\le \theta\le \alpha\}$ with the polar metric $\dd \rho^2+\rho^2\dd
\theta^2$. A blowing-up at $P_i$ consists in remplacing the point $P_i$
with the segment of all the points $(P_i,\theta)_{0\le \theta\le
  \alpha}$ and if $(A_n)$ is a sequence of points of $\inter{\boP}$ converging
to $P_i$ in the original topology we shall say that $(A_n)$ converges to
$(P_i,\theta)$ if $\theta_n\rightarrow \theta$ where
$(\rho_n,\theta_n)$ are the coordinates of $A_n$ near $P_i$. If we make this
blowing-up at all the vertices, we get a new topological space that we
note $\widetilde{\boP}$; $\widetilde{\boP}$ is always a
topological space homeomorphic to the closed unit disk and its
interior is equal to the interior of $\boP$. Then the
question is to understand what is the limit of $Per(A_n)$ when $A_n$
tend to some $(P_i,\theta)$.

Let us consider the case where $i=1$; a neighborhood of $P_1$ in $\boP$
is $\{(\rho,\theta),\ 0\le \rho<\mu,0\le \theta\le \alpha\}$. We know
that there is a bijection between the Alexandrov-embedded $r$-noid with
genus $0$ and horizontal ends and the polygonal immersed disk (see
\cite{CR}). In this bijection, the corresponding $r$-noid to $\boP$
will be noted $\Sigma(\boP)$ and $\Sigma(\boP)^+$ is the conjugate of
a graph over the multi-domain $\Ome(\boP)$ which contains $\boP$ and has
as vertices the vertices of $\boP$. Besides this graph is bounded by
$r$ vertical straight-lines passing by the vertices of $\boP$ (see
\cite{CR} and \cite{Ma1}). Let us
consider $\boC$ the conjugate of the straight-line passing by
$P_1$. $\boC$ is a strictly convex curve and there exists
$\gamma:(-\pi/2,\alpha+\pi/2) \rightarrow \{z=0\}$ a parametrization
of $\boC$ by its normal. We then have the following result.

\begin{prop}{\label{vertex}}
Let $(A_n)$ be a sequence of points in the interior of $\boP$ converging to
the point 
$(P_1,\theta)$ of $\partial \widetilde{ \boP}$. Then $Per(A_n)$ converges to: 
\begin{itemize}
\item $(0,-1)$ if $\theta=0$,
\item $(-\sin\alpha, \cos\alpha)$ if $\theta=\alpha$,
\item $\overrightarrow {\gamma(\theta-\pi/2)\gamma (\theta +\pi/2)}$ or
  $\dis \frac{\overrightarrow {\gamma(\theta-\pi/2)\gamma (\theta
      +\pi/2)}}{|| \overrightarrow {\gamma(\theta-\pi/2)\gamma (\theta
      +\pi/2)}||}$, following the sign of $|| \overrightarrow
  {\gamma(\theta-\pi/2)\gamma (\theta +\pi/2)}||-1$, if $\theta\in
  (0,\alpha)$. 
\end{itemize}
\end{prop}

\subsubsection{Preliminaries}
Let $(A_n)$ be a sequence of points in $\inter{\boP}$ converging to
$(P_1,\theta$) and we suppose that a neighborhood of $P_1$ in $\boP$
is $\{(\rho,\theta),\ 0\le \rho<\mu,0\le \theta\le \alpha\}$. In
$\Ome(\boP)$, a neighborhood of $P_1$ is then $T(-\frac{\pi}{2},
\alpha+\frac{\pi}{2}, \mu)=\{(\rho,\theta),\ 0\le
\rho<\mu, -\frac{\pi}{2}\le \theta\le \alpha+\frac{\pi}{2} \}$. Let
$u_n$ be the restriction of the solution $u_{A_n}$ to the period
${\Ome_A}_0^1$; we can remark that the period ${\Ome_A}_0^1$ can be
identified with $\Ome(\boP) \backslash [A_n,P_1]$ in using the
covering map $\pi$.

We note $\Sigma_n$ the graph in $\R^3$ of $u_n$ and $\Sigma_n^*$ the
minimal surface consisting in the union of the conjugate surface of
$\Sigma_n$ with its symetric with respect to $\{z=0\}$ (the conjugate
surface to $\Sigma_n$ is normalized such that the conjugates of the
vertical lines satisfy $z=0$). We also note $\widetilde{\Sigma}_n^*$ the
periodic minimal surface consisting in the union of the conjugate surface of
the graph of $u_{A_n}$ with its symetric with respect to $\{z=0\}$. In
a certain way, $\Sigma_n^*$ is a period of $\widetilde{\Sigma}_n^*$,
and the vector that lets $\widetilde{\Sigma}_n^*$ invariant is the non
renormalized $Per(A_n)$. 

From Lemma \ref{totcurv}, we remark that the total curvature of
$\Sigma_n^*$ does not depend on $n$ and is $4\pi r$.

\subsubsection{The convergence of the graphs{\label{convgraphe}}}
To understand the behaviour of the surface $\Sigma_n^*$ when $n$ goes
to $+\infty$ we need to know the behaviour of the sequence $(u_n)$. We
consider $u_n$ as a function on $\Ome(\boP) \backslash [A_n,P_1]$;
then the study of the convergence is on the limit multi-domain
$\Ome(\boP)$. Using the arguements of the proof of Theorem
\ref{graph}, we see that there is no line of divergence so a
subsequence $(u_{n'})$ converges to a function $u$ solution of
\eqref{MSE} on $\Ome(\boP)$ and taking the value $+\infty$ on
$\pi(\boL_i^+)$ and $-\infty$ on $\pi(\boL_i^-)$ (more precisely,
if $B$ is a point in $\boP$ we have $u_{n'}-u_{n'}(B)\rightarrow
u$). By Theorem $7$ in \cite{Ma1}, 
such a solution $u$ is unique; then, in fact, the sequence $(u_n)$
converges to $u$. The graph of this function is the conjugate surface
to $\Sigma(\boP)^+$.

\medskip

We also need to study the convergence of $(u_n)$ near the point $P_1$,
and to do this we shall renormalize a neighborhood of $P_1$. $u_n$ is
defined on the neighborhood $T(-\frac{\pi}{2}, \alpha+\frac{\pi}{2},
\mu) \backslash [A_n,P_1]$ of $P_1$; more precisely, the derivatives of
$u_n$ are well defined on $T(-\frac{\pi}{2}, \alpha+\frac{\pi}{2}, \mu)
\backslash \{A_n\}$. We have $A_n=(\rho_n,\theta_n)$ with $\theta_n
\rightarrow \theta$ and $\rho_n \rightarrow 0$. We then renormalized
by $\frac{1}{\rho_n}$: we get a function $v_n$ defined on
$T(-\frac{\pi}{2}, \alpha+\frac{\pi}{2}, \frac{\mu}{\rho_n}) \backslash
\{(\rho,\theta_n),\rho\in[0,1]\}$ by
$v_n(\rho,\beta)=\frac{1}{\rho_n}u_n(\rho_n\rho,\beta)$. $v_n$ is a
solution of the minimal surfaces equation. We want to understand the
asymptotic behaviour of $v_n$. Since $\frac{\mu}{\rho_n}\rightarrow
+\infty$ the limit multi-domain is $T(-\frac{\pi}{2},
\alpha+\frac{\pi}{2}, +\infty)\backslash \{(1,\theta)\}$ for the
derivatives of $v_n$. In the following we note $B(\beta)$ the point of
polar coordinates $(1,\beta)$ and $L(\beta)$ the half straight-line
$\{(\rho,\beta)\}_{\rho>0}$. 

First we must study the lines of divergence. We know that $v_n$ takes
the value $+\infty$ on $L(-\frac{\pi}{2})$ and the value $-\infty$ on
$L(\alpha+\frac{\pi}{2})$; we have $\Psi_{v_n}(P_1)=0$,
$\Psi_{v_n}(B(\theta_n))=0$ and $\Psi_{v_n}\ge 0$. Let $L$ be a line of
divergence, $L$ must have an end-point, otherwise we can apply the
arguement of the proof of Theorem \ref{graph} with the point
$P_1$. This end point can not be on $L(\alpha+\frac{\pi}{2})$ or
$L(-\frac{\pi}{2})$ because of Lemma \ref{lemdiv1}. Then the end point must be
$P_1$ or $B(\theta)$. If the line of divergence has two end-points,
it is the segment $[P_0,B(\theta)]$ then we have
\begin{equation*}
\begin{split}
0=\left|\Psi_{v_n}(B(\theta_n))-\Psi_{v_n}(P_1)\right|
=\left|\int_{[P_1,B(\theta_n)]}\dd\Psi_{v_n} \right|
&\longrightarrow \left|\int_{[P_1,B(\theta)]}\lim\dd\Psi_{v_n}
\right|\\ 
&\longrightarrow 1
\end{split}
\end{equation*}
This is a contradiction.

We then can ensure that $L$ is a half straight-line with $P_1$ or
$B(\theta)$ as end-point. Suppose that the end-point is $P_1$ then
$L$ is some $L(\beta)$. 
\begin{lem}{\label{div1}}
Let $L$ be a line of divergence with $P_1$ as end-point, $L$ is some
$L(\beta)$. Then $\beta \notin (\theta-\frac{\pi}{2},\theta+ \frac{\pi}{2})$  
\end{lem}

\begin{proof}
Let us suppose that $\beta \in (\theta-\frac{\pi}{2}, \theta+
\frac{\pi}{2})$. We note $C_\rho$ the point of $L(\beta)$ with coordinates
$(\rho,\beta)$. Since $\dd\Psi_{v_n}$ is closed, 
$$\int_{[P_1,B(\theta_n)]}\dd \Psi_{v_n}+ \int_{[B(\theta_n),
  C(\rho)]}\dd \Psi_{v_n}+ \int_{[C(\rho), P_1]}\dd \Psi_{v_n}= 0$$
The first integral is always zero, then 
$$ \left| \int_{[C(\rho), P_1]}\dd \Psi_{v_n}\right| \le
|B(\theta_n)C(\rho)|$$
Then by passing to the limit for a subsequence making $L$ appears we get,
$\rho\le 
|B(\theta)C(\rho)|$, but $|B(\theta)C(\rho)|=\sqrt{\rho^2+1-2\rho\cos
  (\beta-\theta)}$ then for big $\rho$ the inequality is not true.
\end{proof}

We suppose now that the end-point of $L$ is $B(\theta)$. We note
$(\rho',\gamma')$ the polar coordinates on $T(-\frac{\pi}{2},\alpha+
\frac{\pi}{2},+\infty)$ with $B(\theta)$ as origin point; $\gamma'$ is
chosen such that the coordinates of $P_0$ in this new coordinates are
$(1,\pi)$. In this polar coordinates, the line of divergence $L$ is some
$L'(\beta)=\{\gamma'=\beta,\ \rho'>0\}$. 

\begin{lem}{\label{div2}}
Let $L$ be a line of divergence with $B(\theta)$ as end point, $L$ is some
$L'(\beta)$. Then $\beta\notin (-\pi,-\frac{\pi}{2})\cup (\frac{\pi}{2},\pi]$.
\end{lem} 

\begin{proof}
The proof is the same as the one of Lemma \ref{div1} in exchanging $P_1$
and $B(\theta)$.
\end{proof}

\begin{figure}[h]
\begin{center}
\resizebox{0.8\linewidth}{!}{\input{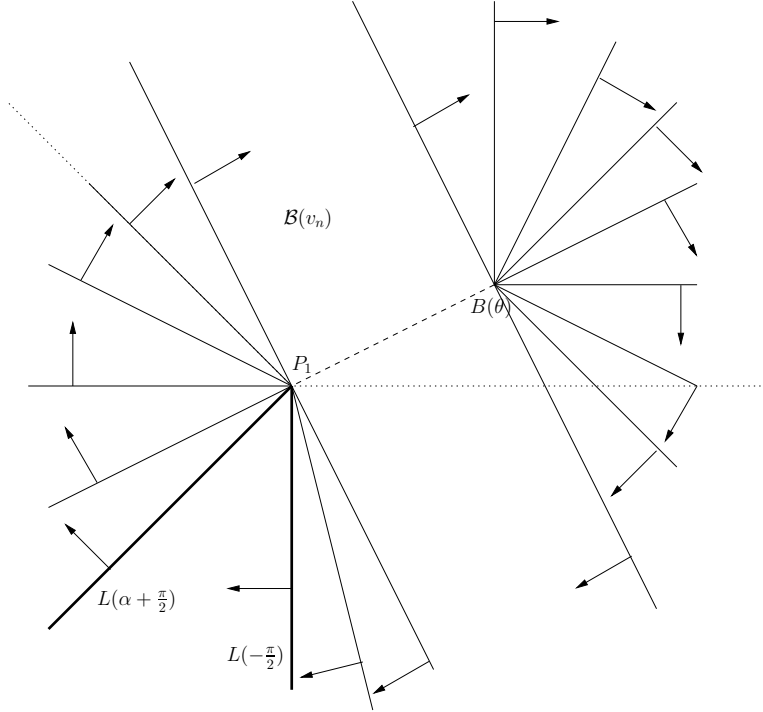}}
\caption{the asymptotic behaviour of $v_n$ \label{strip}}
\end{center}
\end{figure}

In the following, we shall prove that, in fact, all the lines of
divergence that we have not excluded by Lemma \ref{div1} and \ref{div2} yet
appear. We first observe that the 
allowed lines of divergence do not intersect themselves.  Since, for every $n$,
$\Psi_{v_n}(P_1)= 0= \Psi_{v_n}(B(\theta_n))$ and $\Psi_{v_n}\ge 0$,
there is only one possibility for the limit normal on each line of
divergence. Besides, we know that $\boB(v_n)=\{P\in T(-\frac{\pi}{2},
\alpha+ \frac{\pi}{2},+\infty)|\ |\nabla v_n(P)|\textrm{ is
  bounded}\}$ contains a strip which is delimited by the two half
straight-lines with $P_1$ as end-point $L(\theta-\frac{\pi}{2})$ and
$L(\theta+\frac{\pi}{2})$ and the two half straight-lines with
$B(\theta)$ as end-point $L'(-\frac{\pi}{2})$ and $L'(\frac{\pi}{2})$
for the polar coordinates centred on $B(\theta)$. If all the lines of
divergence appear, $\boB(v_n)$ is exactly this strip (see Figure
\ref{strip}). 

Let us suppose that one allowed line of divergence $L$ do not appear,
\emph{i.e.} $L\subset\boB(v_n)$. Then the connected component $\boB$ of
$\boB(v_n)$ that contains $L$ is then a multi-domain which is such
that there exists a subset $K$ such that $\boB \backslash K$ is
isometric to an angular sector (let us observe that the angle at the
vertex can be greater than
$2\pi$) minus the set of the points at a distance less than $d$ from the
vertex of the angular sector ($d$
is a positive number). On $\boB$ we have a subsequence $v_{n'}$ that
converges to some function $v$. Since $\boB$ is bounded by lines of
divergence or by the boundary of $T(-\frac{\pi}{2}, \alpha+
\frac{\pi}{2},+\infty)$ the value of $v$ is $+\infty$ on one side of
the angular sector and $-\infty$ on the other side, by Lemma
\ref{lemdiv2}.
Besides $\Psi_v \ge 0$ since $\Psi_{v_n}\ge 0$ for every $n$. Then, the
function $v$ satisfies many  
conditions that contradict Theorem $2$ in \cite{Ma2}; this proves there is no
sub-sequence such that one of the allowed lines of divergence does not
appear. We then know the limit of the normal to the graph of $v_n$ for
all the points outside the strip.

On the strip, there is a sub-sequence $(v_{n'})$ that converges to some
function $v$ (in fact $v_{n'}$ is not well defined on the strip since
it is defined only outside $[P_1,B(\theta_n)]$, but the
derivatives are well defined and converge to the derivatives of some
function $v$). The function $v$ takes the value $-\infty$ on
$L(\theta+ \frac{\pi}{2})$ and $L'(-\frac{\pi}{2})$ and the value
$+\infty$ on $L(\theta- \frac{\pi}{2})$ and $L'(\frac{\pi}{2})$ by
Lemma \ref{lemdiv2}; such a solution $v$ is unique and is a peace of
helicoid. More precisely, if the strip is isometricaly
parametrized by $\R\times[-1/2,1/2]$ with $P_1=(0,1/2)$ and
$B(\theta)=(0,-1/2)$ then we have $v(x,y)=x\tan(\pi y)$.

\begin{rem}{limK}
Let $P$ be a point in the strip, the curvature of the helicoid at this point
is non zero, then the curvature of the graph of $v_n$ over the 
point $P$ goes to the curvature of the helicoid at this point. This
implies that there exists a sequence of point $p_n\in \Sigma_n$  such
that the curvature of $\Sigma_n$ at $p_n$ goes to $+\infty$. 
\end{rem}

\begin{rem}{limcn}
Besides, we know that there exists, for each $n$, a constant $c_n$ such
that $u_{A_n}\circ f_n=u_{A_n}+c_n$, then the result above proves that
$\dis \frac{c_n}{\rho_n}\rightarrow+\infty$.
\end{rem}

\medskip

We also want to know what occurs when we homothetically expand the 
sequence
$u_n$ in a general way. Let $(M_n)$ be a sequence of point in $\Ome(\boP)$ such that, for every
$n$, $M_n\neq A_n$ and $(\lambda_n)$ a sequence of positive number such
that $\lambda_n$ goes to $+\infty$. Let $h_n$ be the homothety of center
$M_n$ and ratio $\lambda_n$, by applying the dilatation $h_n$ we get a
function $v_n$ on $h_n(\Ome(\boP)) \backslash \{A_n\}$ defined by $v_n(M)=
\lambda_n 
u_n(h_n^{-1}(M))$. The question we ask is: what asymptotic behaviours
are possible? By taking a subsequence, if it is necessary, we can
suppose that we are in one of the following cases.

\begin{case}{c1}
For every $i$, $d(M_n,h_n(P_i))\rightarrow+\infty$ and $d(M_n,h_n(A_n))
\rightarrow+\infty$. In this case the limit multi-domain for the
sequence $v_n$ is $\R^2$ then if there is no line of divergence a
sub-sequence $(v_{n'})$ must converge to a linear function by Bernstein
Theorem. If there is a line of divergence $L$, each connected
component of the domain of convergence of a subsequence $(v_{n'})$is a strip
or a half-plane and $(v_{n'})$ converges to a function $v$ with the value
$+\infty$ on one side and $-\infty$ on the other side; but such
solution of \eqref{MSE} does not exist by Proposition $1$ in \cite{Ma2} so the
domain of convergence is 
empty and we have only lines of divergence which are all parallel to
$L$ and the limit normal is constant on $\R^2$. 
\end{case}

\begin{case}{c2}
For every $i$, $d(M_n,h_n(P_i))\rightarrow+\infty$ and $d(M_n,h_n(A_n))
\rightarrow d\ge 0$. In this case the limit multi-domain is $\R^2$
minus the limit point $A'$ of $h_n(A_n)$. Since
$d(M_n,h_n(P_1))\rightarrow+\infty$, we have
$d(h_n(A_n),h_n(P_1))\rightarrow+\infty$ then $\lambda_n\rho_n
\rightarrow +\infty$; by Remark \ref{limcn}, this implies that
$\lambda_n c_n$, which is the 
vertical jump over $[h_n(A_n),h_n(P_1)]$ for the function $v_n$, goes
to $+\infty$. Then the derivatives of $v_n$ can not
converge on $\R^2 \backslash \{A'\}$ and there are lines of
divergence. Since $\Psi_{v_n}(h_n(A_n))=0$ and $\Psi_{v_n}\ge 0$, using
arguements that we have already seen, we can ensure that a line of
divergence must be a half straight-line with $A'$ as end-point and
that there is only one possibility for the limit normal along the line
of divergence. If the domain of convergence $\boB(v_n)$ is non empty,
each connected component of it is an angular sector of $\R^2$; then on
one component, a subsequence $(v_{n'})$ converges to a solution $v$ of
$\eqref{MSE}$ with 
the value $+\infty$ on one side and $-\infty$ on the other side.
By Proposition $2$ in \cite{Ma2}, such a solution does not exist and
we have only the lines of divergence as asymptotic behaviour.
\end{case}

\begin{case}{c3}
There exists $j\neq 1$ such that $d(M_n,h_n(P_j))\rightarrow
c\ge0$. This implies that, for $i\neq j$,
$d(M_n,h_n(P_i))\rightarrow+\infty$ and $d(M_n,h_n(A_n))\rightarrow
+\infty$; then the limit multi-domain is an angular sector isometric
to some $T(0, \beta,+\infty)$ with $P_j'=\lim h_n(P_j)$ as vertex. As
above, the lines of divergence must be half straight-lines with $P_j'$
as end-point and the limit normal on a line of divergence is given by
the condition $\Psi_{v_n}\ge 0$. As in Case \ref{c2}, $\boB(v_n)$ is
empty and the asymptotic behaviour is given by the lines of
divergence.  
\end{case}

\begin{case}{c4}
$d(M_n,h_n(P_1))\rightarrow c\ge0$ and $d(h_n(P_1),h_n(A_n))\rightarrow 0
\textrm{ or }+\infty$. In this case, the limit multi-domain is an
angular sector with $P_1'=\lim h_n(P_1)$ as vertex and the asymptotic
behaviour is the same as in Case \ref{c3}.
\end{case}

\begin{case}{c5}
$d(M_n,h_n(P_1))\rightarrow c\ge 0$ and $d(h_n(P_1),h_n(A_n))
\rightarrow c>0$. The limit-multi-domain is then an angular sector,
with $P_1'$ as vertex, minus the point $A'=\lim h_n(A_n)$. We are in
the situation studied above and we know that the asymptotic behaviour
is a domain of convergence which is a strip where $(v_n)$ converges to a
piece of helicoid and, outside the strip, lines of divergence with
$A'$ or $P_1'$ as end-point (see Figure \ref{strip}).
\end{case}


\subsubsection{The convergence of $\Sigma_n^*$}

We first fix the notations. Let $Q$ be a point in the interior of $\boP$ and
$q_n$ 
the corresponding point in $\Sigma_n$. For every $n$, the value
$\Psi_{u_n}(Q)$ is well defined and we normalize $\Sigma_n^*$ such
that $q_n^*$, the corresponding point to $q_n$, has coordinates
$(0,0,\Psi_{u_n}(Q))$; $\Sigma_n^*$ is then still symetric with respect to
the plane $\{z=0\}$. We want to determine the limit of the sequence of
minimal surfaces $(\Sigma_n^*)$.

Let $M$ be a surface in $\R^3$. In the following, when we shall
talk about a geodesical disk $D(m,\mu)$ of center $m\in M$ and radius
$\mu$, we shall consider the disk of radius $\mu$ in the tangent
plane to $M$ at $m$ with the exponential map $\exp_m: T_mM\rightarrow
\R^3$. Besides, we shall sometimes identify a point in the tangent
plane with its image by $\exp_m$. 

Let $m_n^*$ be a sequence of point of $\R^3$ such that $m_n^*$ lies in
the part of $\Sigma_n^*$ which is the conjugate of $\Sigma_n$ (the
third coordinate of $m_n^*$ is non-negative). We make the following
assumption: $K_{\widetilde{\Sigma}_n^*}(m_n^*)\rightarrow
+\infty$. Because of 
Remark \ref{limK}, we know that such a sequence exists. Let $\mu$ be a
positive 
number, and, for each $n$, we consider $D(m_n^*,\mu)$ the closed
geodesical disk of center $m_n^*$ and radius $\mu$ in
$\widetilde{\Sigma}_n^*$. On $D(m_n^*,\mu)$, the function $a\mapsto
|K_{\widetilde{\Sigma}_n^*}(a)|(\mu-d(a,m_n^*))^2$ admits a maximum in
the interior 
of the disk (the distance $d$ is the distance in the tangent space)
and let 
$p_n^*$ be a point where the maximum is reached. Let us note
$\lambda_n=\sqrt{|K_{\widetilde{\Sigma}_n^*}(p_n^*)|}$ and $\mu_n=
(\mu-d(p_n^*,m_n^*))$. The geodesical disk $D(p_n^*,\mu_n)$ is
included in 
$D(m_n^*,\mu)$ (in fact the image by the exponential map in
$\widetilde{\Sigma}_n^*$ of the geodesical disk $D(p_n^*,\mu_n)$ is included
in 
the image of $D(m_n^*,\mu)$); besides since $p_n^*$ realizes the
maximum we have 
$\lambda_n^2 \mu_n^2\ge K_{\widetilde{\Sigma}_n^*}(m_n^*)\mu^2$, then
$\lambda_n \mu_n
\rightarrow +\infty$ and $\lambda_n\rightarrow +\infty$. By
translating $p_n^*$ to the origin and homothetically expanding the disk
$D(p_n^*,\mu_n)$ by the factor $\lambda_n$, we obtain a new geodesical
disk $D'_n=D(0,\lambda_n \mu_n)$, $(D'_n)$ is a sequence of minimal
surfaces. We have $|K_{D'_n}(0)|=1$. Let $R$ be a positive number and
$\tilde{a}$ a point in the geodesical disk $D(0,R)$ of $D'_n$; we then
note $a$ the point in $D(p_n^*,\mu_n)$ corresponding to
$\tilde{a}$. $\exp_{p_n^*}(a)$ is included in the image of the disk
$D(m_n^*,\mu)$. Then there exists $a'\in D(m_n^*,\mu)$ such that
$\exp_{m_n^*}(a')=\exp_{p_n^*}(a)$ and $d(m_n^*,a')\le
d(m_n^*,p_n^*)+d(p_n^*,a)$; we have $K_{\widetilde{\Sigma}_n^*}(a)=
K_{\widetilde{\Sigma}_n^*}(a')$. With this notation we then have:

\begin{equation*}
\begin{split}
|K_{D'_n}(\tilde{a})|(\lambda_n \mu_n-R)^2\le
|K_{D'_n}(\tilde{a})|(\lambda_n \mu_n-d(\tilde{a},0))^2
&=|K_{\widetilde{\Sigma}_n^*}(a)|(\mu_n-d(a,p_n^*))^2\\
&\le |K_{\widetilde{\Sigma}_n^*}(a')|(\mu-d(a',m_n^*))^2\\
&\le \lambda_n^2 \mu_n^2
\end{split}
\end{equation*}
the equality is due to the fact that the function
$|K_{\widetilde{\Sigma}_n^*}(\cdot)|(\mu_n-d(\cdot,p_n^*))^2$ is
invariant 
under rescaling. Thus the curvature on $D'_n$ is uniformly bounded on
the sequence of geodesical disks $D(0,R)$ of $D'_n$. Then there
exists a subsequence $(D'_{n'})$ that converges to a complete minimal
surface that we denote $D'_\infty$; this surface is complete since
$\lambda_n \mu_n\rightarrow +\infty$. Besides $D'_\infty$ is non flat
since at the origin its curvature is $-1$. Since $D'_\infty$ is non
flat there is a point $\tilde{a}$ where the normal has a negative third
coordinate; there exists a neighborhood $U$ of $\tilde{a}$ in the
tangent 
plane to $D'_\infty$ at $\tilde{a}$ such that $D'_\infty$ and
$D'_{n'}$, for big $n'$, are graphs over $U$ and $D'_{n'} \rightarrow
D'_\infty$ as graphs. Let $\tilde{a}_{n'}$ be the sequence of points
in $D'_{n'}$ over $\tilde{a}$ as graphs. The normal at
$\tilde{a}_{n'}$ to $D'_{n'}$ has a negative third coordinate then
before the rescaling $\tilde{a}_{n'}$ correponds to a point $a_{n'}$
which lies in the conjugate of $\Sigma_{n'}$. Let $b_{n'}$ be the
point in $\Sigma_{n'}$ corresponding to $a_{n'}$ and $B_{n'}$ the
projection on $\Ome(\boP)$ of $b_{n'}$. The convergence of $D_{n'}$ to
$D'_\infty$ near $\tilde{a}$ says us that $\boB(v_{n'})$ is non empty 
where $v_{n'}$ is the rescaled function of $u_{n'}$ with the factor
$\lambda_{n'}$ and $B_{n'}$ as origin points. We are then in Case \ref{c1}
or Case \ref{c5} of the preceding subsection, but if it is Case \ref{c1} the
limit graph would be 
a plane which is impossible since $D'_\infty$ is non flat. Then it is
Case \ref{c5} and the limit graph is a piece of helicoid then this
implies that $D'_\infty$ is a catenoid. This also implies that
$\lambda_{n'}\sim \frac{c}{\rho_{n'}}$ with $c$ some real constant (we
recall that $\rho_n=|P_1A_n|$). 

\begin{rem}{maxK}
If for example, we take $m_n^*$ such that $K_{\Sigma_n}(m_n^*)=\max
K_{\Sigma_n}(\cdot)$ the above arguements show that $\max
K_{\Sigma_n}(\cdot)=O(\left(\frac{1}{\rho_n}\right)^2)$.  
\end{rem}

\medskip

The boundary of $\Sigma_n^*$ is composed of two closed paths
$\Gamma_n^1$ and $\Gamma_n^2$, we want to understand the behaviour of
this boundary when $n$ goes to $+\infty$. The boundary of the graph
$\Sigma_n$ is composed of $r-1$ straight-lines, they are over the
points $P_i$ for $i\neq 1$, and a curve which consists in a half
straight-line over $P_1$ that goes down from the infinity to some
point called $t_n^1$, a curve which is a graph over the segment
$[P_1,A_n]$ joining $t_n^1$ to some point $t_n^2$, a vertical segment
$[t_n^2,t_n^3]$ over $A_n$,  a curve which is a graph over the segment
$[P_1,A_n]$ joining $t_n^3$ to some point $t_n^4$ (by a vertical
translation, this curve is the same as the curve joining $t_n^1$ to
$t_n^2$) and a half straight-line over $P_1$ with $t_n^4$ as end-point
and going down to the infinity. Then the path $\Gamma_n^1$
(resp. $\Gamma_n^2$) consists in the conjugate of the curve joining
$t_n^1$ to $t_n^2$ (resp. $t_n^3$ to $t_n^4$) with its symetric with
respect to the plane $\{z=0\}$. 

As in Subsection \ref{convgraphe}, we consider $v_n$ the rescaled
function with factor $\frac{1}{\rho_n}$ on
$T(-\frac{\pi}{2},\alpha+\frac{\pi}{2}, \frac{r}{\rho_n}) \backslash
[P_1,B(\theta_n)]$. By Remark \ref{maxK}, there exists a constant which
bounds the curvature on the graph of $v_n$ for every $n$. Let $p_n$ be
a point in the graph of $v_n$ which is above the middle of the
segment $[P_1,B(\theta_n)]$ and $\boY_n$ the
conjugate of the graph of $v_n$ that we extend by symmetry and
periodicity such that the conjugate point of $p_n$ is the origin of
$\R^3$. Let $R$ be 
a positive number and consider the sequence of geodesical disks
$D(0,R)$ in $\boY_n$. Since the curvature is uniformly bounded on the
disk there exists a subsequence such that the sequence of disks
converges, this implies that, using a Cantor diagonal process, there
exists a subsequence $(\boY_{n'})$ which converges to some minimal
surface $\boY$. Since the graphs of $v_n$ converge near $p_n$ to a
piece of an helicoid, $\boY$ is a catenoid whose flux is the
vector $2\overrightarrow{P_1B(\theta)}$. There is only one possible
limit for subsequences of $(\boY_n)$ so the sequence $(\boY_n)$ converges to
$\boY$. Let us consider the catenoid given in cylindrical coordinates by
$(u,v) \mapsto (\frac{1}{\pi}u,v, \frac{1}{\pi}\argch(u))$. Thus $\boY$ is the
translated by $(0,0,-\frac{1}{\pi})$ of the image by a rotation of axis
$\{z=0, x\cos\theta+y\sin\theta=0\}$ and angle $\frac{\pi}{2}$ of this
catenoid. Suppose that $p_n$ is the point in the graph of $v_n$ which is
the limit of the points $(1/2,\beta,v_n(1/2,\beta))$ when $\beta\rightarrow
\theta_n$ with $\beta<\theta_n$. We note $\boY_n'$ the conjugate of the graph
of $v_n$ that we extend by symmetry but not by periodicity such that the
conjugate of the point $p_n$ is the origin. What we have proved just above
implies that $(\boY_n')$ converges to the part of the catenoid $\boY$
included in $\{x\cos\theta+y\sin\theta\le 0\}$; this part is noted $\boY_-$. If
$p_n$ is build with 
$\beta>\theta_n$ we get the half catenoid included in
$\{x\cos\theta+y\sin\theta\ge 0\}$.  

This proves that the 
rescaled paths $\frac{1}{\rho_n}\Gamma_n^1$ and $\frac{1}{\rho_n}
\Gamma_n^2$ converge to two circles of the same radius. Let us call
$s_n^i$ the point in $\Sigma_n^*$ which is the conjugate of $t_n^i$ for
$1\le i\le 4$. The above proof shows that for every $\epsilon>0$ and
for big $n$ the ball of center $s_n^1$ (resp. $s_n^4$) and radius
$\epsilon$ contains $\Gamma_n^1$ (resp. $\Gamma_n^2$) and even a
bigger and bigger part of a surface near a half catenoid. We consider the 
case where $p_n$ is build with $\beta<\theta_n$. If $B(R)$ is the ball of
center the origin and radius $R$, we have $B(R)\cap \boY_n'$ is near from
$B(R)\cap \boY_-$ for big $n$. Besides for big $n$, $\rho_n R<\epsilon$ then
the ball $B(s_n^1,\epsilon)$ contains the homothetic by $\rho_n$ of a surface
near $B(R)\cap \boY_-$. The same is true for $p_n$ build with
$\beta>\theta_n$.  
This implies that the sequence of
total curvatures of the part of $\Sigma_n^*$ included in this ball has
an lower limit bigger than $2\pi$, since the total curvature of an half
catenoid is $2\pi$. In fact, we have to remember that near the points
$s_n^1$ and $s_n^4$ the surface $\Sigma_n^*$ behaves like small
half-catenoid for big $n$. 

\medskip

We know that $(\Sigma_n^*)$ is a sequence of symetric minimal
surface with finite total curvature $4\pi r$. Each surface
$\Sigma_n^*$ has two connected components of boundary $\Gamma_n^1$ and
$\Gamma_n^2$. We know that the diameter of each component goes to zero
and if $B_n$ is a sequence of balls of diameter $\epsilon$ centred at
$s_n^2$ then the sequence of total curvatures of 
$\Sigma_n^*\cap B_n$ has an lower limit bigger than $2\pi$. By
taking a subsequence we can suppose that $(s_n^2)$ and $(s_n^3)$ diverge
to the infinity or converge to $s_\infty^2$ and $s_\infty^3$. Then by
results explained in \cite{CR}, there exist a finite number of distinct
properly and simply immersed branched minimal surfaces
$M_1,\dots,M_k\subset\R^3$ with finite total curvature, a finite subset
$X\in\R^3$ contained in $M=M_1\cup\cdots\cup M_k$ and a subsequence of
$(\Sigma_n^*)$, that we still call $(\Sigma_n^*)$, such that:
\begin{enumerate}
\item $(\Sigma_n^*)$ converges to $M$ (with finite multiplicity)
  on compact subsets of $\R^3 \backslash (X\cup\{s_\infty^2,
  s_\infty^3\})$ in the $C^m$-topology for any positive integer $m$;
\item on each $M_i$ the multiplicity $m_i$ is well defined and is such
  that 
$$m_1C(M_1)+\cdots+m_k C(M_k) \le C(\Sigma_n^*);$$
\item $X$ is the singular set of the limit $\Sigma_n^*\rightarrow
  M$. Given a point $p\in X$, the amount of total curvature of the
  sequence $(\Sigma_n^*)$ which disappears through the point $p$ is a 
  positive multiple of $4\pi$.
\end{enumerate}

Since at $s_\infty^2$ and $s_\infty^3$, there is $2\pi$ of total
curvature that disappear, we have $m_1C(M_1)+\cdots+m_k C(M_k) \le
C(\Sigma_n^*)-4\pi=4\pi(r-1)$ even if $s_n^2$ or $s_n^3$ diverges. We
know that the sequence of functions 
$u_n$ on $\Ome(\boP)\backslash [P_1,A_n]$ converge on $\Ome(\boP)$ to
the solution $u$ of the Dirichlet problem asked in Theorem $7$ of
\cite{Ma1}. Then $q_n$ converges to $q_\infty$ the point in the graph of $u$
which is above $Q$. The conjugate surface to the graph of $u$, after
being extended by symetry with respect to the plane $\{z=0\}$, is the
solution $\Sigma(\boP)$ to the plateau problem at infinity with genus
$0$ for the data $\boP$. Let $q_\infty^*$ be the point in
$\Sigma(\boP)$ that correspond to $q_\infty$. We then have
$q_n^*\rightarrow q_\infty^*$ and, in a neighborhood of $q_\infty^*$,
$\Sigma_n^*$ converges to $\Sigma(\boP)$, then $\Sigma(\boP)$ is one $M_i$. We
then can suppose that 
$M_1=\Sigma(\boP)$, since $C(\Sigma(\boP))=4\pi (r-1)$, $m_1=1$ and for,
$i\neq 1$, $M_i$ is a plane and besides $X$ is empty. In fact, in
\cite{CR}, C.~Cos\'\i n and A.~Ros proves that, in such a convergence,
no plane can appear. Then finally, $\Sigma_n ^*\rightarrow
\Sigma(\boP)$.

Since the problems of convergence appear only near the points
$s_\infty^2$ and $s_\infty^3$, the curves in $\Sigma(\boP)$ which are
the conjugates of the $r-1$ straight-lines that are over the points
$P_i$, for $i\neq 1$, are the respective limits of the curves in $\Sigma_n^*$
which are the conjugates of the $r-1$ straight-lines in $\Sigma_n$ that
are over the points $P_i$, for $i\neq 1$.

\subsubsection{The convergence of $Per(A_n)$}

$Per(A_n)$ corresponds to the vector that defines the translation under
which $\widetilde{\Sigma}_n^*$ is invariant. Then $Per(A_n)$ is
$\overrightarrow{s_n^2s_n^3}$ or $\frac{\overrightarrow{s_n^2s_n^3}}{
  ||\overrightarrow{ s_n^2s_n^3}||}$, following the value of
$||\overrightarrow{ s_n^2s_n^3}||$. 

As above, we consider $\boC$ the curve in $\Sigma(\boP)$ which is the
conjugate of the 
vertical straight-line which is over $P_1$ in the graph of $u$. This
curve is strictly convex then we can parametrized $\boC$ by its
normal.  Then there exists a parametrization $\gamma:
(-\frac{\pi}{2}, \alpha+\frac{\pi}{2})\rightarrow \{z=0\}$ of $\boC$
such that the normal to $\gamma(\beta)$ is $(\sin\beta,-\cos\beta,0)$;
we cover $\boC$ as we cover the straight-line over $P_1$ when we go
down. We know that $\boC$, outside $s_\infty^2$ and $s_\infty^3$, is the
limit, in the Hausdorff topology, of 
a part of the conjugate of the boundary of the graph $\Sigma_n$ of
$u_n$ which is included in $\{z=0\}$. For each $n$ this set is
composed of three strictly convex arc, we note $\boC_n^1$, $\boC_n^2$
and 
$\boC_n^3$: $\boC_n^1$ is the conjugate of the vertical half
straight-line that have $t_n^1$ as end-point ($\boC_n^1$ has $s_n^1$
as end-point), $\boC_n^2$ is the conjugate of the vertical segment
$[t_n^2,t_n^3]$ ($\boC_n^2$ is joining $s_n^2$ to $s_n^3$) and
$\boC_n^3$ is the conjugate of the vertical half straight-line that
have $t_n^4$ as end-point ($\boC_n^3$ has $s_n^4$ as end-point). As
$\boC$, this three strictly convex arcs can be parametrized by their normal,
then 
there exist $\gamma_n^1:(-\frac{\pi}{2},\theta_n] \rightarrow  \{z=0\}$,
$\gamma_n^2:[\theta_n-\pi,\theta_n+\pi] \rightarrow  \{z=0\}$ and
$\gamma_n^3:[\theta_n,\alpha+\frac{\pi}{2}) \rightarrow  \{z=0\}$ such
that $\gamma_n^i$ parametrized $\boC_n^i$ and the normal at the point
$\gamma_n^i(\beta)$ is $(\sin\beta,-\cos\beta,0)$. We have
$\gamma_n^1(\theta_n)=s_n^1$, $\gamma_n^2(\theta_n-\pi)=s_n^2$,
$\gamma_n^2(\theta_n+\pi)=s_n^3$ and $\gamma_n^3(\theta_n)=s_n^4$. We
note $I_n^i$ the definition set of $\gamma_n^i$. We then have
$I_n^i\rightarrow I^i$ where $I^1=(\frac{\pi}{2},\theta]$, $I^2=
[\theta-\pi, \theta+\pi]$ and $I^3=[\theta,\alpha+\frac{\pi}{2})$. The
question is: is $\gamma_n^i$ converging on $I^i$?

Since there is a half catenoid that appears near the points $s_n^i$ for
big $n$, if $\beta\in(\theta-\frac{\pi}{2},\theta]\subset I^1$ the
curvature at $\gamma_n^1(\beta)$ becomes infinite, the same is true
for the sequence of point $(\gamma_n^2(\beta))$ for
$\beta\in[\theta-\pi, \theta-\frac{\pi}{2}) \cup
(\theta+\frac{\pi}{2}, \theta+\pi]\subset I^2$ and the sequence of
points $(\gamma_n^3(\beta))$  for $\beta\in
[\theta,\theta+\frac{\pi}{2})\subset I^3$.

The set $\boC\backslash\{s_\infty^2,s_\infty^3\}$ is composed of a
finite number of convex arcs. Let $a$ be a point in $\boC \backslash
\{s_\infty^2,s_\infty^3\}$, at this point the convergence of
$\Sigma_n^*$ to $\Sigma(\boP)$ well behaves so there exists a
neighborhood $U$ of $a$ in the tangent plane to $\Sigma(\boP)$ at $a$
(this plane is vertical) such that, over $U$, $\Sigma(\boP)$ and
$\Sigma_n^*$, for big $n$, are graphs and, as graphs, $\Sigma_n^*$
converge to $\Sigma(\boP)$. Over $U$, there is a neighborhood of $a$
in $\boC \backslash \{s_\infty^2,s_\infty^3\}$ and this neighborhood
is the limit of the part of $\Sigma_n^*$ which is included in
$\{z=0\}$. Since each unit vector is reached a finite number of times on
$\boC_n^1\cup\boC_n^2 \cup \boC_n^3$, by taking a subsequence $n'$, there
exist $\beta$, $\epsilon$ and $i\in\{1,2,3\}$ such that, for every
$n'$, the part of $\Sigma_{n'}^*\cap \{z=0\}$ which is over $U$ contains
$\gamma_{n'}^i(\beta-\epsilon,\beta+\epsilon)$. Then the convergence
as graphs of $\Sigma_n^*$ implies that on $(\beta- \epsilon,\beta+
\epsilon)$ the sequence $(\gamma_{n'}^i)$ converges to some map
$\widetilde{\gamma}$ that parametrized a neighborhood of $a$ in $\boC$
such that the normal to the point $\widetilde{\gamma}(\omega)$ is
$(\sin\omega, -\cos\omega,0)$ (the convergence is in the $C^m$
topology for every $m$). Since $\gamma_{n'}^i\rightarrow
\widetilde{\gamma}$, for every $\omega\in
(\beta-\epsilon,\beta+\epsilon)$ the curvature at the point
$\gamma_{n'}^i(\omega)$ remains bounded. Then, if $i=1$, $(\beta-
\epsilon, \beta+\epsilon)\cap [\theta-\frac{\pi}{2},\theta]=
\emptyset$, if $i=2$, $(\beta-\epsilon, \beta+\epsilon)\cap
\left([\theta-\pi,\theta- \frac{\pi}{2}]\cup [\theta+ \frac{\pi}{2},
  \theta+\pi]\right)= \emptyset$ and, if $i=3$, $(\beta-\epsilon,
\beta+\epsilon)\cap [\theta, \theta+ \frac{\pi}{2}]= \emptyset$. 

Then by applying the above arguement to a countable number of points and
constructing a subsequence by diagonal Cantor process, there exist
\begin{enumerate}
\item $k_3$ open intervals in $I=(-\frac{\pi}{2},\alpha+
  \frac{\pi}{2})$: $J_1,\dots,J_{k_1}\subset I^1\backslash
  [\theta-\frac{\pi}{2},\theta]$, $J_{k_1+1},\dots,J_{k_2}\subset
  I^2\backslash \left([\theta-\pi,\theta- \frac{\pi}{2}]\cup [\theta+
    \frac{\pi}{2}, \theta+\pi]\right)$ and
  $J_{k_2+1},\dots,J_{k_3}\subset I^3 \backslash [\theta, \theta+
  \frac{\pi}{2}]$, these intervals satisfy $J_j\cap J_l=\emptyset$ if
  $j\neq l$, and
\item a map $\widetilde{\gamma}$ with value in $\{z=0\}$, defined on
  the union of the $k_3$ intervals,  $\widetilde{\gamma}$ parametrizes
  $\boC\backslash \{s_\infty^2,s_\infty^3\}$ by its normal,
\item a subsequence $n'$,
\end{enumerate} 
such that, on $J_j$, the sequence $(\gamma_{n'}^i)$ (for
the corresponding $i$) converges to $\widetilde{\gamma}$; besides, the
$J_j$ are maximal in the sense that $\widetilde{\gamma}$ restricted to
$J_j$ parametrized one of the convex arcs that composed $\boC\backslash
\{s_\infty^2,s_\infty^3\}$. The curve $\boC$ has total curvature
$\alpha+\pi$ and since $\widetilde{\gamma}$ parametrizes by the normal
$\boC\backslash \{s_\infty^2,s_\infty^3\}$, which has the same total
curvature, we have
\begin{equation}{\label{eqcourb}}
\alpha+\pi=\sum_{j=1}^{k_3}l(J_j)\le l(I)=\alpha+\pi
\end{equation}
where $l(J_j)$ is the length of the interval $J_j$.

\medskip

$\bullet$ We now assume that $\theta\notin \{0,\alpha\}$. Under this
hypothesis,the computation \eqref{eqcourb} implies that there exist points 
\begin{multline*}
-\frac{\pi}{2}=\beta_0<\beta_1<\cdots<\beta_{k_1}=\theta-
\frac{\pi}{2} <\beta_{k_1+1}<\cdots\\
\cdots<\beta_{k_2}= \theta+ \frac{\pi}{2}
<\beta_{k_2+1}<\cdots<\beta_{k_3}=\alpha+\frac{\pi}{2}
\end{multline*}
such that, for every $j$, $J_j=(\beta_{j-1},\beta_j)$. 

\begin{lem}
for $j\notin\{ 0,k_1,k_2,k_3\}$, the sequence $(\gamma_{n'}^i)$ (for the
corresponding $i$) converges in a neighborhood of $\beta_j$ in the
$C^1$ topology, then we can extend the definition of
$\widetilde{\gamma}$ at $\beta_j$.
\end{lem}

\begin{proof}
We apply Proposition \ref{convcomp}, if there is no $C^1$ convergence near
$\beta_j$, since $\Sigma_n^*\rightarrow \Sigma(\boP)$ in the Hausdorff
topology, $\Sigma(\boP)$ contains a segment. But this is not true
then we can extend the definition of $\widetilde{\gamma}$ at
$\beta_j$.  
\end{proof}

We then have $\gamma_{n'}^1\rightarrow \widetilde{\gamma}$ on
$(-\frac{\pi}{2},\theta- \frac{\pi}{2})$, $\gamma_{n'}^2\rightarrow
\widetilde{\gamma}$ on $(\theta- \frac{\pi}{2},\theta+ \frac{\pi}{2})$
and $\gamma_{n'}^3\rightarrow \widetilde{\gamma}$ on $(\theta+
\frac{\pi}{2},\alpha+ \frac{\pi}{2})$.

We shall then study the behaviour near $\theta-\frac{\pi}{2}=\theta'$.

\begin{figure}[h]
\begin{center}
\resizebox{1.0\linewidth}{!}{\input{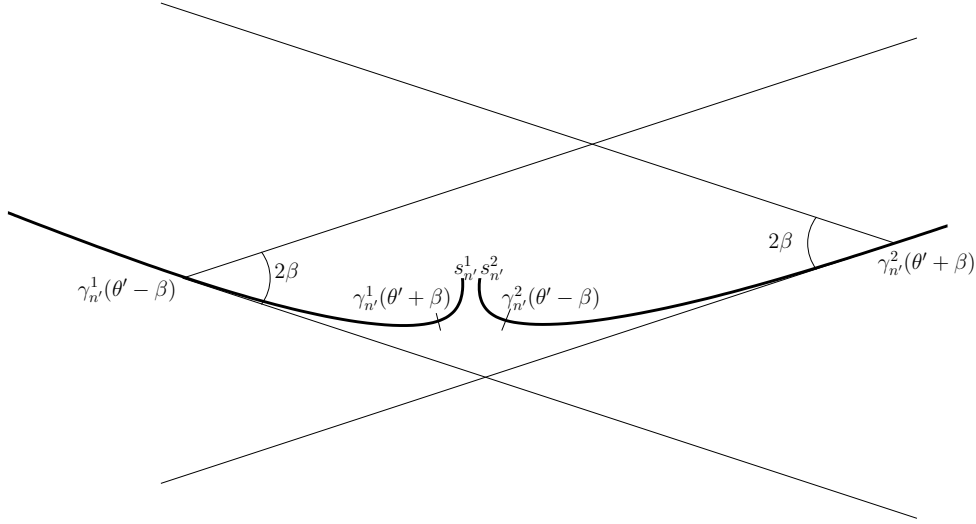}}
\caption{the local behaviour near $s_{n'}^2$ \label{figsecto}}
\end{center}
\end{figure}

Let $\beta<\frac{\pi}{4}$ be a small
positive angle and $\epsilon>0$. We apply Proposition \ref{convencadr} to
$\gamma_{n'}^1$ on $(\theta'-\beta,\theta'+\beta)$, then the set
$\gamma_{n'}^1 (\theta'-\beta,\theta'+\beta)$ is included in an
angular sector, we note $S$, of vertex $\gamma_{n'}^1(\theta'-\beta)$
and angle $2\beta$. We have $\gamma_{n'}^1(\theta'+\beta)\in S$. Then,
for big $n'$, because of the behaviour of the sequence $(\Sigma_{n'}^*)$
near the point $s_{n'}^2=\gamma_{n'}^2(\theta)$, we have the distance
between $\gamma_{n'}^1(\theta'+\beta) $ and $s_{n'}^2$ less than
$\epsilon$ and the distance between $\gamma_{n'}^1(\theta'+\beta) $ and
$\gamma_{n'}^2(\theta'-\beta)$ less than $\epsilon$. Besides by
applying Proposition \ref{convencadr} to $\gamma_{n'}^2$ on
$(\theta'-\beta,\theta'+\beta)$, we have that $\gamma_{n'}^2(\theta'-
\beta,\theta'+\beta)$ is included in an angular sector of vertex
$\gamma_{n'}^2(\theta'-\beta)$ and angle $2\beta$. Since at
$\gamma_{n'}^1(\theta'-\beta)$ and $\gamma_{n'}^2(\theta'-\beta)$, the
two strictly convex curves $\boC_{n'}^1$ and $\boC_{n'}^2$ have the
same normal and their curvature have the same sign, the angular sector
of vertex $\gamma_{n'}^2(\theta'-\beta)$ is just the translation of
$S$ by the vector $\overrightarrow{ \gamma_{n'}^1(\theta'-\beta)
  \gamma_{n'}^2(\theta'-\beta)}$. This second angular sector is then
included in the set $S_\epsilon$ of points of $\R^2$ which are at a
distance less than $\epsilon$ from $S$. Then by passing to the limit
$n'\rightarrow +\infty$, we get that $\widetilde{\gamma}\left(
  (\theta'-\beta, \theta'+\beta)\backslash\{\theta'\} \right)$ is
included in the set of points in $\R^2$ which are at a distance less
than $\epsilon$ from an angular sector which have $\widetilde{\gamma}
(\theta'-\beta) $ as vertex and $2\beta$ as angle. The points
$s_{n'}^2$ are also in this set for big $n'$ (see Figure \ref{figsecto}).

We can do the same work in starting from the point $\gamma_{n'}^2
(\theta'+\beta)$ and in covering the curves $\boC_{n'}^1$ and
$\boC_{n'}^2$ in the opposite sense. We then get that
$\widetilde{\gamma} \left( (\theta'-\beta,
  \theta'+\beta)\backslash\{\theta'\} \right)$ and $s_{n'}^2$, for
  big $n'$, are
included in the set of points in $\R^2$ which are at distance less
than $\epsilon$ from an angular sector which have $\widetilde{\gamma}
(\theta'+\beta)$ as vertex and $2\beta$ as angle. Since $\beta$ is
less than $\frac{\pi}{4}$, the intersection of the two sets we have
just build is a compact subset of $\R^2$ which contains $s_{n'}^2$, then
$s_\infty^2$ exists. Then
$\widetilde{\gamma}(t)$ admits 
two limits: one when $t\rightarrow \theta'$, $t<\theta'$ and one when
$t\rightarrow \theta'$, $t>\theta'$. This two limit points are in
$\boC$ and $\widetilde{\gamma}$ parametrized $\boC$ except for two
points where the normal is different. Since the normals at this two
limit points are the same, the two limit points are, in fact, a unique
point and then the definition of $\widetilde{\gamma}$ extend at
$\theta'$. Now, letting $\beta$ and $\epsilon$ goes to zero, we see
that the limit compact is just $\widetilde{\gamma}(\theta')$ then
$\widetilde{\gamma}(\theta')=s_\infty^2$.

In the same way we see that $s_\infty^3$ exists and
$\widetilde{\gamma}(\theta+\frac{\pi}{2})= s_\infty^3$. In fact this prove
that $\widetilde{\gamma}$ parametrizes $\boC$ by its normal on $I$ and then 
$\widetilde{\gamma}=\gamma$. Then the sequence $(Per(A_n))$ has only one
possible cluster point then the sequence converges to the limit given in
Proposition \ref{vertex}.

\medskip

$\bullet$ We now study the case $\theta=0$ (the case $\theta=\alpha$ can be
done 
in the same way). In this case $k_1=0$ and using the same arguements
as above we can prove that: $s_\infty^3$ exists,
$\gamma=\widetilde{\gamma}$ and
$\gamma(\frac{\pi}{2})=s_\infty^3$. Because of the asymptotic
behaviour of the curve $\boC$, we know that there exists
$0<\beta<\frac{\pi}{2}$ such that for every $t\in (-\frac{\pi}{2},
-\frac{\pi}{2}+\beta)$, $|\gamma(t)s_\infty^3|>2$ and $\frac{
  \overrightarrow{\gamma(t)s_\infty^3}}{|\gamma(t)s_\infty^3|}$ is at
a distance less than $\epsilon$ from the vector $(0,-1)$. Let
$0<\beta'<\beta$, we consider the angular sector $S$ of vertex
$\gamma(-\frac{\pi}{2}+\beta')$ and angle $2\beta$ which contains
$\gamma(-\frac{\pi}{2}, -\frac{\pi}{2}+\beta')$ and has, as a part of
its boundary, the half straight-line tangent to $\boC$ at
$\gamma(-\frac{\pi}{2}+\beta')$ and with
$\gamma(-\frac{\pi}{2}+\beta')$ as end-points; $S$ exists because of
Proposition \ref{convencadr}. Then, because of the asymptotic behaviour of
$\boC$, there exist $0<\beta'<\beta$ and $\epsilon'$ such that, for
every point $s$ at a distance less than $\epsilon'$ from $S$, we have
$|ss_\infty^3|>2$ and $\frac{\overrightarrow {ss_\infty^3}}{|s
  s_\infty^3|}$ is at a distance less than $\epsilon$ from $(0,-1)$. 

We apply now Proposition \ref{convencadr} to $\gamma_{n'}^2$ on
$(-\frac{\pi}{2} 
-\beta, -\frac{\pi}{2}+\beta)$ that we cover in the opposite sense. We
then obtain that $\gamma_{n'}^2 (-\frac{\pi}{2} -\beta, -\frac{\pi}{2}
+\beta)$ is included in an angular sector with vertex $\gamma_{n'}^2
(-\frac{\pi}{2}+\beta')$ and angle $2\beta$. Since $\gamma_{n'}^2
\rightarrow \gamma$ on $(-\frac{\pi}{2}, -\frac{\pi}{2}+\beta)$, these
angular sector converges to $S$, in fact, the angular sectors we have
build are just the translations of $S$ by the vector $\overrightarrow
{\gamma(-\frac{\pi}{2}+\beta') \gamma_{n'}^2(-\frac{\pi}{2}+\beta')}$
which goes to zero. Besides we know that the distance
between $\gamma_{n'}^2 (-\frac{\pi}{2}-\beta')$ and $s_{n'}^2$ goes to
zero. Then for big $n'$, $s_{n'}^2$ is at a distance less than
$\epsilon'$ from $S$. Then, using that $s_{n'}^3$ goes to
$s_\infty^3$, this proves that for big $n'$, $|s_{n'}^2s_{n'}^3|>2$
and $\frac{ \overrightarrow {s_{n'}^2s_{n'}^3}}{|s_{n'}^2s_{n'}^3|}
\rightarrow (0,-1)$.

Since in each case, there is only one possible limit for $(Per(A_{n'}))$
this proves that $(Per(A_n))$ converges to this limit. Proposition \ref{vertex}
is then established.

\subsection{Conclusion}

We use the preceding subsections to extend $Per$ to $\widetilde{\boP}$ then we
get a continuous map on $\widetilde{\boP}$. First we have the following remark.
\begin{prop}
The period map does not vanish on the boundary of $\widetilde{\boP}$ 
\end{prop}
 
\begin{proof}
The only points where $Per$ can vanish are points in the vertices and
if $Per(A)=0$ we have $\overrightarrow{ \gamma(\theta-\frac{\pi}{2})
  \gamma(\theta+\frac{\pi}{2})}=0$ for some strictly convex curve
$\gamma$. But $\gamma$ on $(\theta-\frac{\pi}{2},
\theta+\frac{\pi}{2})$ is a graph over a straight-line then the above
vector can not be $0$.
\end{proof}

This proves that we can compute the degree of the period map along the
boundary of $\widetilde{\boP}$.

\begin{thm}{\label{degree}}
The degree of the period map along the boundary of $\widetilde{\boP}$
is $-(r-1)$
\end{thm}

\begin{proof}
The edges of $\boP$ does not contribute toward the degree. So only the
behaviour at the vertices is interesting for the degree. Let us compute the
contribution of the vertex $P_1$; we use the
notation of the preceding section. Since the curve $\gamma$ is
strictly convex, the map $\theta \mapsto \frac{ \overrightarrow
  {\gamma(\theta-\frac{\pi}{2}) \gamma(\theta +\frac{\pi}{2})}}{|
  \gamma(\theta-\frac{\pi}{2}) \gamma(\theta +\frac{\pi}{2})|}$ is a
monotone map. Let $\theta$ be in $[0,\alpha]$, the unit vector 
tangent to $\gamma$ at $\gamma(\theta)$ is $(\cos\theta,\sin\theta)$ which
turns in the 
clockwise sense, when $\theta$ increases. Besides, for
$0<\theta<\alpha$, $\gamma([\theta-\pi/2,\theta+\pi/2])$ is a graph
over a straight-line generated by $(\cos\theta,\sin\theta)$. This
implies that $Per(P_1,\theta)\cdot (\cos\theta,\sin\theta)$ never
vanishes ($(P_1,\theta)$ is a point of the boundary of
$\widetilde{\boP}$); by looking at the behaviour for small $\theta$ we
have 
$Per(P_1,\theta)\cdot (\cos\theta,\sin\theta)\ge 0$. Besides, for
$\theta=0$, the basis composed of $(1,0)$ and $Per(P_1,0)$ is
an indirect one and, for $\theta=\alpha$, the basis
($(\cos\alpha,\sin\alpha),Per(P_1,\alpha)$ is a direct 
one. Then, when $\theta$ increases from $0$ to $\alpha$,
$Per(P_1,\theta)$ describes an angle  $\alpha$ (since the unit tangent
describes this angle) plus $\pi$ (since the basis composed of the unit
tangent at $\gamma(\theta)$ and $Per(P_1,\theta)$ is an orthonormal indirect
one for $\theta=0$ 
and an orthonormal direct one for $\theta= \alpha$). Since, when we
describe $\partial\widetilde{\boP}$ in the clockwise sense, $\theta$
decreases, the contribution of the vertex $P_1$ towards the degree is
$-\frac{1}{2\pi}(\alpha+\pi)$. 

The degree is then $\dis-\frac{1}{2\pi}(r\pi +\sum_{i=0}^{r-1}\alpha_i)$ 
where $\alpha_i$ is the inner angle at the vertex $P_i$. Then by
applying Gauss-Bonnet Theorem to $\boP$, $\dis \sum_{i=0}^{r-1}
\alpha_i=r\pi-2\pi$ and then the degree is $-(r-1)$. 
\end{proof}

Theorem \ref{degree} then proves that the degree of the period map is
non zero 
and then there exists $A\in\inter{\boP}$ such that $Per(A)=0$ and then Theorem
\ref{main} is proved.

\section{An other example of solution the Plateau problem at 
infinity{\label{rotation}}}
Theorem \ref{main} gives a wide class of solutions of the Plateau problem at
infinity with genus $1$. But it gives no example of polygon which is the the
flux polygon of an $r$-noid of genus $1$ but not the flux polygon of an
$r$-noid of genus $0$. Corollary \ref{regpoly} gives such polygons. In fact we
study the case where the multi-domain with cone singularity bounded by the
polygon is invariant under a "rotation".

\begin{thm}{\label{rot}}
Let $V$ be a polygon that bounds a multi-domain with cone singularity
$(D,Q,\psi)$. We suppose that $D$ satisfies the hypothesis H and that there
exists an isometry $h$ of $D$ such that 
$\psi \circ h =R \circ \psi$ where $R$ is the rotation in $\R^2$ with center
$\psi(Q)$ and angle $\alpha\in(0, 2\pi)$. Then the period vector associated to 
$D$ vanishes and there exists an Alexandrov-embedded $r$-noid with genus $1$
and horizontal ends having $V$ as flux polygon.
\end{thm}

\begin{proof}
From Construction \ref{C1} and \ref{C2}, we have a multi-domain with
logarithmic singularity $(\Ome,\boQ,\phi)$ and Theorem \ref{graph} gives us a
function $u$ on $\Ome$ which is unique up to an additive constant. By
construction, the map $h$ can be lifted to $\Ome$ to a map $\tilde{h}$
which is an isometry of $\Ome$ such that $\phi_{\psi(Q)} \circ \tilde{h}= R
\circ \phi_{\psi(Q)}$. Since $\tilde{h}$ is an isometry of $\Ome$, $u\circ
\tilde{h}$ is a solution of the same Dirichlet problem as $u$ then $u
\circ \tilde{h}= u+c$ with $c\in \R$. Then the two equations $\phi_{\psi(Q)}
\circ \tilde{h}= R \circ \phi_{\psi(Q)}$ and $u \circ \tilde{h}= u+c$ imply
that:
$$\tilde{h}^* (\dd X_1^*, \dd X_2^*)= R(\dd X_1^*, \dd X_2^*)$$
with $\dd X_1^*$ and $\dd X_2^*$ the $1$-forms associated to $u$ as in
Subsection \ref{periodproblem}. Then, if $\Gamma$ is a lift of a generator of
$\pi_1(D\backslash \{Q\})$, we have:
\begin{equation*}
\begin{split}
\int_\Gamma (\dd X_1^*, \dd X_2^*)&= \int_{\tilde{h}(\Gamma)} (\dd X_1^*, \dd
X_2^*) \\ 
&=\int_\Gamma \tilde{h}^*(\dd X_1^*, \dd X_2^*) \\
&=\int_\Gamma R(\dd X_1^*, \dd X_2^*) \\
&=R\left(\int_\Gamma (\dd X_1^*, \dd X_2^*)\right)
\end{split}
\end{equation*} 
The first equality is due to the fact that $\Gamma$ and $\tilde{h}(\Gamma)$
are lifts of two closed pathes that give the same generator of
$\pi_1(D\backslash \{Q\})$. Besides $R$ has a unique fixed point, since
$\alpha\in(0,2\pi)$, which is $0$. Then the period vector vanishes.

\end{proof}

We then can give some examples of polygons $V$ that satify this condition; we
consider the case where $V$ is a regular convex polygon or a regular star
polygon (see \cite {Cox}). 

\begin{figure}[h]
\begin{center}
\resizebox{1\linewidth}{!}{\input{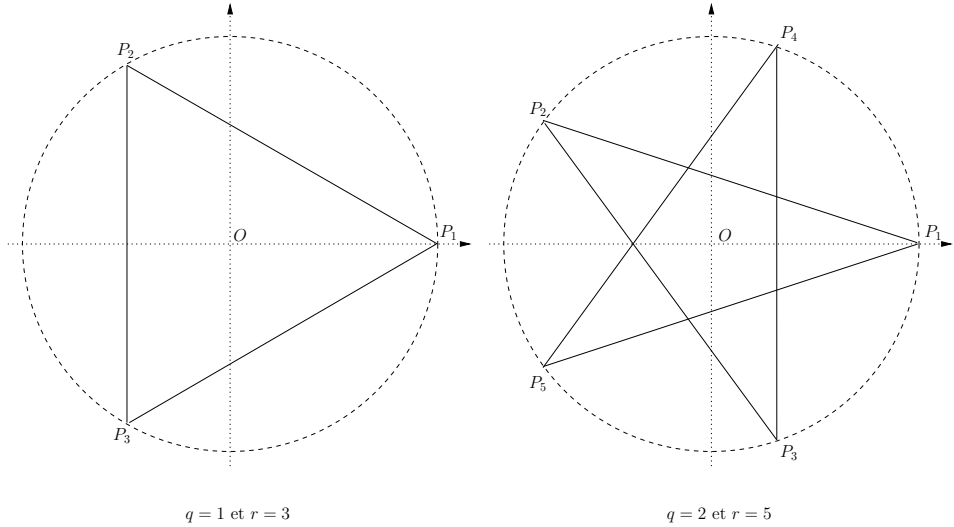}}
\caption{Examples of polygons studied in Corollary \ref{regpoly}
  {\label{fig8}}} 
\end{center}
\end{figure}

\begin{cor}{\label{regpoly}}
Let $q$ and $r$ be in $\N^*$ such that $\gcd(q, r)=1$ and $2q<r$. We note, for
$i=1,\dots,r+1$, $\dis P_i=e^{2(i-1)\sqrt{-1}\frac{q}{r}\pi}\in \C= \R^2$
($P_1=P_{r+1}$) (see Figure \ref{fig8}). Then there exists an
Alexandrov-embedded $r$-noid of genus $1$ and horizontal ends with
$(\overrightarrow{P_1P_2}, \overrightarrow{P_2P_3}, \dots,
\overrightarrow{P_rP_{r+1}})$ as flux polygon. 
\end{cor}

\begin{proof}
The idea of the proof is to build a multi-domain with cone singularity bounded
by the above polygon and satisfying the hypothesis of Theorem \ref{rot}. Let
$\R^2$ be identified with $\C$. We note $T$ the set of
$(\rho,\theta)\in\R_+\times[0,\pi]$ such that $\rho e^{i\theta}$ is in the
triangle $P_1P_2O$, where $O$ is the origin. We note $D$ the set of
$(\rho,\theta)$ in $\R^+\times [0,2q\pi]$ such that $(\rho,\theta)$ is in $D$
if $(\rho,\theta')$ is in $T$ where $\theta=n
\left(2\frac{q}{r}\pi\right)+\theta'$ is the only writing with $n\in\Z$ and
$\theta'\in [0,2\frac{q}{r}\pi)$ (like an euclidean division of $\theta$ by
$2\frac{q}{r}\pi$). We remark $\{0\}\times [0,2q\pi]\subset D$ and $(\rho,0)\in
D$ iff $(\rho,2q\pi)\in D$ iff $\rho\in[0,1]$. Because of these remarks, we
can consider the polar metric on $D$ (we identify all the points $(0,\theta)$
and call the new point $Q$) and indentify the point $(\rho,0)$ with $(\rho,
2q\pi)$ for $0\le \rho\le 1$. Then if we consider on $D$ the map
$\psi:(\rho,\theta)\mapsto (\rho\cos\theta,\rho \sin \theta)$, the triplet
$(D,Q,\psi)$ is a multi-domain with cone singularity which is bounded by the
polygon $(\overrightarrow{P_1P_2}, \overrightarrow{P_2P_3}, \dots,
\overrightarrow{P_rP_{r+1}})$. The angle of $D$ is $2q\pi$. Besides on $D$ the
map 
$$
h: (\rho,\theta)\longmapsto
\left\{\begin{array}{cc}
(\rho,\theta+2\frac{q}{r}\pi)&\textrm{if } \theta\le
2q\pi-2\frac{q}{r}\pi\\ 
(\rho,\theta-2\frac{r-1}{r}q\pi)&\textrm{if } \theta\ge
2q\pi-2\frac{q}{r}\pi 
\end{array}
\right.$$
is well defined and is an isometry of $D$. Besides $\psi\circ h=R\circ \psi$
with $R$ the rotation of center $O$ and angle $2\frac{q}{r}\pi<\pi$. Then we
can apply Theorem \ref{rot}.
\end{proof}

When $q=1$, the polygon is a convex regular polygon and the multi-domain $D$
that we build is in fact an immerssed polygonal disk, then Theorem \ref{main}
gives also the result but what we now know is that in the proof of Theorem
\ref{main} we can choose $A$ as the isobarycenter of the polygon.

When $q>1$, the polygon $(\overrightarrow{P_1P_2}, \overrightarrow{P_2P_3},
\dots, \overrightarrow{P_rP_{r+1}})$ does not bound an immersed polygonal
disk. Then we get new examples of polygons which are flux polygons of an
Alexandrov-embedded $r$-noid with genus $1$.

\begin{rem}{histoire}
In \cite{JM}, L.~P.~Jorge and W.~H.~Meeks give Weierstrass data for
$r$-noid with genus $0$ and horizontal ends having as flux polygon the polygon
studied in Corollary \ref{regpoly} with $q=1$. Then Corollary \ref{regpoly}
gives $r$-noids similar to the Jorge-Meeks examples for the genus $1$ case. In
fact this examples are known by H.~Karcher which gave in \cite{Ka} Weierstrass
data for some $r$-noids which correspond to the ones we have just built. These
Weierstrass data are expressed in terms of the Weierstrass
$\mathfrak{p}$-function. 
\end{rem}


\appendix

\section{Convergence of sequences of solution of \eqref{MSE} and line
  of divergence{\label{lineofdiv}}} 

The aim of this appendix is to explain some results on the convergence
of sequences of solutions of the minimal surface equation.

Let $(u_n)$ be a sequence of solution of \eqref{MSE} on a multi-domain
$D$ ($D$ has possibly a cone or logarithmic singularity), the first
result is then

\begin{Prop}
If $(u_n)$ is an uniformly bounded sequence on $D$, there exists a
subsequence which converges to a solution of \eqref{MSE} on $D$. The
convergence is uniform on every compact subset of $D$.
\end{Prop}

This result proves that if the sequence $(||\nabla u_n||)$ is
uniformly bounded on every compact subset of $D$ there exists a
subsequence $(u_{n'}-u_{n'}(P))$, where $P$ is a point in $D$, which
converges to a solution of \eqref{MSE} on $D$ ($D$ is connected). We
then can define the domain of convergence of the sequence $(v_n)$ as
the set of the points $P\in D$ where the sequence $(||\nabla u_n (P)||)$
is bounded, we note this set $\boB(u_n)$. We know that $\boB(v_n)$ is
an open subset of $D$ such that $(||\nabla u_n||)$ is uniformly bounded
on every compact subset of $\boB(u_n)$ (see \cite{Ma1}). Then, on every
connected component of the domain of convergence, we can make converge
a subsequence $(u_{n'}-u_{n'}(P))$. In fact, in our paper we often write that
a subsequence $(u_{n'})$ converges instead of $(u_{n'}-u_{n'}(P))$ but, since
the value on the boundary are often infinite, this does not matter; when
it is necessary we use $(u_{n'}-u_{n'}(P))$. The question is to understand
the domain of convergence, in fact we shall understand $D\backslash
\boB(u_n)$. 

Let $P$ be in $D\backslash \boB(u_n)$, then for a subsequence
$W_{n'}\rightarrow +\infty$ and since the normal to the graph at the
point over $P$ is 
$$N_{n'}(P)=\left(\frac{p_{n'}}{W_{n'}}(P), \frac{q_{n'}}{W_{n'}}(P),
-\frac{1}{W_{n'}}(P) \right)$$
(we use euclidean coordinates in a neighborhood of $P$) we can
suppose that $N_{n'}(P)\rightarrow N$ where $N$ is a unit horizontal
vector. We then have the following result.

\begin{Thm}
Let $(D,\psi)$ be a multi-domain. Let $(u_n)$ be a sequence of solutions
of \eqref{MSE} on $D$. Let $P\in D$ and $N$ be a unit horizontal
vector and $L$ the geodesic of $D$ passing by $P$ and normal to
$N$. If the sequence $(N_n(P))$ converges to $N$, then $N_n(A)$
converges to $N$ at every point $A$ of $L$
\end{Thm} 
Since $D$ is locally isometric to $\R^2$ we can see $N$ as a vector in
$\R^2$ and then the vector $N$ is well defined at all the points of
$D$ in fact this definition coincides with the parallel transport of
$N$. The proof of this theorem is in \cite{Ma1}.

Then in our situation $N_{n'}(P)$ is converging to a unit horizontal
vector then $N_{n'}$ converges to this vector along a straight-line
$L$. Such a line $L$ is called a \emph{line of divergence} of the
sequence since $L$ must be included in $D\backslash \boB(u_n)$. Then
the set $D\backslash \boB(u_n)$ is an union of geodesics of $D$ and
the question is what are the possible lines of divergence: the
following lemmas give some tools to answer to this discussion.

First, we observe that the existence of a line of divergence has a
consequence on the sequence of $1$-forms $\dd\Psi_{u_n}$. If
$N_n\rightarrow N$ along a line of divergence and $T$ is a segment
included in this line of divergence then $\dis\int_T\dd\Psi_{u_n}$
converges to $|T|$ the length of $T$ if the orientation of $T$ is such
that $N$ is the right-hand normal. Since $(\Psi_{u_n})$ is a sequence
of $1$-Lipschitz continuous function, we can always take a subsequence
and suppose that it has a limit $\Psi$, then the above remark on
$\dd\Psi_{u_n}$ allows us to make some calculations on $\Psi$.

As in Figure \ref{strip}, when we make a picture to explain the convergence of
a sequence of solutions of \eqref{MSE}, we draw the limit
normal along the lines of divergence to explain the asymptotical behaviour.

We then have the two following results concerning the lines of divergence and
the convergence of sequences of solutions of \eqref{MSE}

\begin{Lem}{\label{lemdiv1}}
Let $(u_n)$ be a sequence of solutions of \eqref{MSE} on
$[0,1]^2$ such that, for all $n$, the function $u_n$ tends to
$+\infty$ on $\{1\}\times ]0,1[$. Then, no line of divergence of the
sequence $(u_n)$ has $(1,\frac{1}{2})$ as end-point. 
\end{Lem}

\begin{proof}
Let us suppose that there exists such a line of divergence $L$. We note
$A=(1,0)$, $B=(1,\frac{1}{2})$ and $C=(1,1)$. Let us consider a point
$M\in L$. We suppose that the limit normal along $L$ is such that the basis
composed of $\overrightarrow{MB}$ and the limit normal is direct. Since
$\dd\Psi_{u_n}$ is closed, we then have:  

$$\int_{[A,B]}\dd\Psi_{u_n}+ \int_{[B,M]}\dd\Psi_{u_n}+
\int_{[M,A]}\dd\Psi_{u_n}=0$$ 
Since $u_n$ takes the value $+\infty$ on $[A,B]$, this equality proves
that:

$$|AB|+\int_{[B,M]}\dd\Psi_{u_n}\le |MA|$$
But, by our choice of limit normal, we have
$\int_{[B,M]}\dd\Psi_{u_n}\rightarrow |BM|$ for a subsequence, then
$|AB| +|BM|\le |MA|$ which contradicts the triangle inequality. If the
limit normal is the opposite of the one we consider, we can do the
same arguement with the triangle $BCM$.
\end{proof}

\begin{Lem}{\label{lemdiv2}}
Let $(u_n)$ be a sequence of solutions of \eqref{MSE} on
$[0,1]^2$ such that $(u_n)$ converges in the interior of the square to
a solution $u$ and such that we are in one of the following two cases
\begin{itemize}
\item for all $n$, the function $u_n$ tends to $+\infty$ on
$\{1\}\times ]0,1[$ or,
\item for all $n$, $u_n$ is the restriction to $[0,1]^2$ of a solution
$v_n$ of \eqref{MSE} defined on $[0,2]\times[0,1]$ and for $y\in
]0,1[$ we have $N_n(1,y)\rightarrow(1,0,0)$.
\end{itemize}
Then, the limit function $u$ tends to $+\infty$ on $\{1\}\times ]0,1[$.
\end{Lem}

\begin{proof}
let $0<\epsilon<\frac{1}{2}$, for $x\in[0,1]$ we note $A_x=(x,\epsilon)$
and $B_x=(x,1-\epsilon)$. Since $\dd\Psi_{u_n}$ is closed, for $x<1$
we have:
$$\int_{[A_x,B_x]}\dd\Psi_{u_n}- \int_{[A_1,B_1]}\dd\Psi_{u_n}=
\int_{[B_x,B_1]}\dd\Psi_{u_n}+ \int_{[A_1,A_x]}\dd\Psi_{u_n}$$
Then, by passing to the limit and using that $\int_{[A_1,B_1]}
\dd\Psi_{u_n} \rightarrow 1-2\epsilon$ in the two cases, we obtain
that:
$$\left|\int_{[A_x,B_x]}\dd\Psi_u-(1-2\epsilon)\right|\le 2(1-x)$$

This proves that $\dd\Psi_u=\dd y$ on $[A_1,B_1]$. Let $x\in]0,1[$ and
$v$ be the solution of \eqref{MSE} on $A_xB_xB_1A_1$ such that $v$
tends to $+\infty$ on $[A_1,B_1]$ and $v=u$ on the rest of the
boundary we shall prove that $u=v$. 

If $u\neq v$ there exist $\eta\neq 0$ such that $\Ome=\{u-v>\eta\}$ is
non-empty. The boundary of $\Ome$ is composed of one part included in
the interior of $A_xB_xB_1A_1$ and one part included in
$[A_1,B_1]$. Since $\dd\Psi_u$ and $\dd\Psi_v$ are closed we have
$$\int_{\partial \Ome} \dd\Psi_u-\dd\Psi_v =0$$
But on the part of the boundary included in $[A_1,B_1]$ the integral
is $0$ since $\dd\Psi_u=\dd y=\dd\Psi_v$ and on the part included in
the interior of $A_xB_xB_1A_1$ the integral is negative by Lemma $2$
in \cite{CK}. Then we have a contradiction and $u=v$ then $u$ goes to
$+\infty$ on $[A_1,B_1]$.
\end{proof}

\section{Convex curves in $\R^2$}

A convex curve in $\R^2$ is a curve such that its geodesical curvature
has always the same sign. A curve is strictly convex if this
geodesical curvature is positive or negative.

If a curve $c:s\mapsto c(s)\in \R^2$ is convex, the map with value in
$\S^2$ that associates to the parameter $s$ the normal to $c$ at $c(s)$
is a monotone map and, if $c$ is strictly convex, the above map is
strictly monotone. This implies that, if $c:I\rightarrow \R^2$ is a
strictly convex curve, there exists $h:J\rightarrow I$ a
diffeomorphism such that for $\beta\in J$ the normal to $c$ at the
point $c\circ h(\beta)$ is $(\sin\beta,-\cos\beta)$. We then say that
$c$ is parametrized by its normal.

\begin{Prop}{\label{convgraph}}
Let $c:I\rightarrow \R^2$ be a strictly convex curve parametrized by
its normal. If $I\subset(\beta_0+\frac{\pi}{2},\beta_0+\frac{3\pi}{2})$,
the curve $c$ is a graph over the straight line
$y\cos\beta_0-x\sin\beta_0=0$.
\end{Prop}

\begin{proof}
We suppose that $\beta_0=0$. If $c(\beta)=(x(\beta),y(\beta))$, the
map $\beta\mapsto x(\beta)$ is a local diffeomorphism since the second
coordinate of the normal is always positive. Besides $\beta\mapsto
x(\beta)$ is injective: if $x(\beta')=x(\beta'')$ there would exist,
by Rolle's Theorem, $\beta \in(\beta',\beta'') $ such that the normal
at $c(\beta)$ is horizontal, it is impossible. This proof that
$\beta\mapsto x(\beta)$ is a global diffeomorphism then $c$ is a graph
over $y=0$.
\end{proof}

\begin{figure}[h]
\begin{center}
\resizebox{0.8\linewidth}{!}{\input{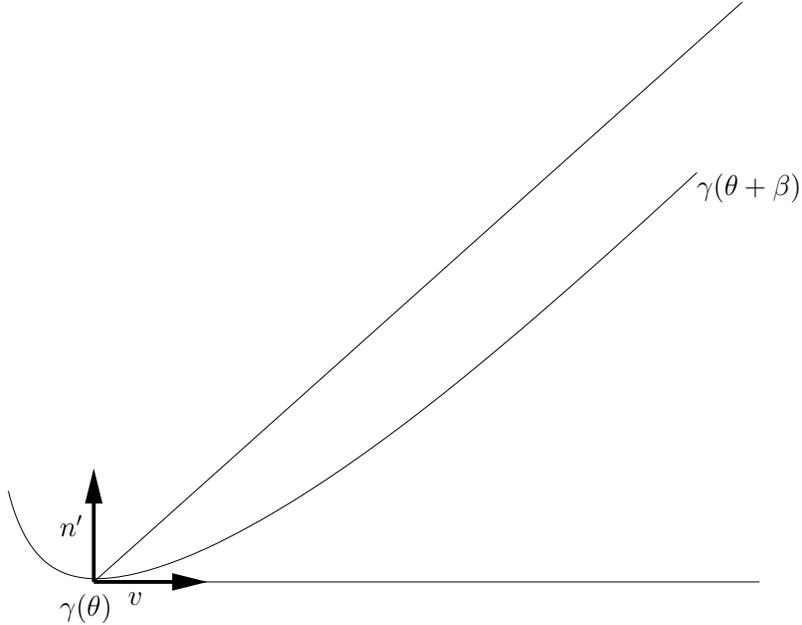}}
\caption{the situation of Proposition \ref{convencadr} \label{fig12}}
\end{center}
\end{figure}

\begin{Prop}{\label{convencadr}}
Let $c$ be a strictly convex curve parametrized by its normal on
$(\theta-\epsilon,\theta+\beta]$ with $\beta\le\frac{\pi}{2}$. Let $v$
be the unit tangent vector to $c$ at $c(\theta)$ and $n'$ the unit
vector which is the normal to $c$ at $c(\theta)$ if the curvature is
positive or the opposite of the normal if the curvature is
negative. Let $w$ be the unit vector $\cos\beta v+ \sin\beta n'$. Then
 the curve $c([\theta,\theta+\beta])$ is in the angular sector
 delimited by the two half straight-lines with $c(\theta)$ as
 end-point and respectively generated by $v$ and $w$ (see Figure
 \ref{fig12}). 
\end{Prop}

\begin{proof}
We can suppose that $\theta=0$, the curvature is negative and
$c(0)=(0,0)$. In this case $c([0,\beta])$ is the graph of a function $f$
over a segment $[0,a]$ by Proposition \ref{convgraph}. By hypotheses,
$f$ is a 
convex function and $f'(0)=0$. This implies first that $f\ge 0$. If
$\beta=\frac{\pi}{2}$ the proposition is proved. If $\beta<
\frac{\pi}{2}$ and $c([0,\beta])$ is not in $\{(x,y)|\ x\ge 0,\ y\ge
0,\ y\le x\tan\beta\}$ there is a parameter $b\in[0,a]$ such that
$f(b)=b\tan\beta$ and then there exists $d<b$ such that
$f'(d)=\tan\beta$. Since the normal map is injective on $[0,\beta]$,
we have $(d,f(d))=c(\beta)$ which is impossible since $d<a$. The
proposition is then proved 
\end{proof}

\begin{Lem}{\label{convlim}}
Let $c$ be a strictly convex curve parametrized by its normal on
$[\beta-\epsilon,\beta)$. We suppose that the curve $c$ is included in
a compact of $\R^2$. Then $c(t)$ converges when $t\rightarrow \beta$.
\end{Lem}

\begin{proof}
Since we are in a compact part of $\R^2$, it is enough to prove that
there is only one cluster point for $c(t)$. Let us suppose that
$c(t_n)\rightarrow p_1$ and $c(s_n)\rightarrow p_2$ and, for every
$n$, $t_n<s_n<t_{n+1}$. We then apply Proposition \ref{convencadr} on
$[t_n,\beta)$ and $[s_n,\beta)$. This proves that $p_1$ is in an
angular sector with vertex $c(s_n)$ and $\beta-s_n$ as angle and $p_2$
is in an angular sector with vertex $c(t_n)$ and $ \beta-t_n$ as
angle. Letting $n$ goes to $+\infty$, we obtain that $p_1$ is in
the half straight-line with $p_2$ as end-point and generated by $v$
(where $v$ is the limit unit tangent vector) and $p_2$ is in
the half straight-line with $p_1$ as end-point and generated by
$v$. This situation is possible only if $p_1=p_2$.
\end{proof}

\begin{Prop}{\label{convcomp}}
Let $(c_n)$ be a sequence of stricly convex curves parametrized by their
normal on $(\theta-\epsilon,\theta+\epsilon)$. We suppose that $(c_n)$
converges to $\widetilde{c}$ on $(\theta-\epsilon,\theta) \cup
(\theta, \theta+\epsilon)$ in the $C^1$ topology. Then: 
\begin{itemize}
\item we have $\widetilde{c}(t)\rightarrow p_1$ when
  $t \rightarrow \theta$, $t<\theta$ and $\widetilde{c}(t)\rightarrow
  p_2$ when $t \rightarrow \theta$, $t>\theta$,
\item as sets, $c_n(\theta-\epsilon,\theta+\epsilon)$ converges to
  $\widetilde{c} (\theta-\epsilon,\theta)\cup \widetilde{c} (\theta,
  \theta+\epsilon) \cup [p_1,p_2]$.
\end{itemize}
If $p_1=p_2$, we have moreover that $(c_n)$ converges to $\widetilde{c}$ (that
we extend by $\widetilde{c}(\theta)=p_1$) on
$(\theta-\epsilon,\theta+\epsilon)$ in the $C^1$ topology.
\end{Prop}

\begin{proof}
$\epsilon$ is supposed to be small and we choose
$\epsilon'<\epsilon$. We then apply Proposition \ref{convencadr} to 
$c_n(\theta-\epsilon',\theta+\epsilon')$ at the points
$c_n(\theta-\epsilon')$ and $c_n(\theta +\epsilon')$ and we get
that, for every $n$, $c_n(\theta-\epsilon',\theta+\epsilon')$ is
included in an angular sector of vertex $c_n(\theta-\epsilon')$ and 
angle $2\epsilon'$ and an angular sector of vertex
$c_n(\theta+\epsilon')$ and angle $2\epsilon'$ (here, we apply
Proposition \ref{convencadr} to the curve $c$ that we cover in the opposite
sense). Letting $n$ goes to 
$+\infty$, we get that $\widetilde{c}((\theta-\epsilon', \theta)\cup
(\theta,\theta+\epsilon'))$ is included in the intersection of two
angular sectors of angle $2\epsilon'$, one has $\widetilde{c}
(\theta-\epsilon')$ as vertex the other has $\widetilde{c}
(\theta+\epsilon')$ as vertex. The intersection of this two angular
sector is a compact; then, by Lemma \ref{convlim}, $p_1$ and $p_2$
exists. We then 
have also proved that the cluster points of the sequence of curves
$(c_n)$ are $\widetilde{c} ((\theta-\epsilon, \theta)\cup
(\theta,\theta+\epsilon))$ and points included in the intersection of
the two angular sectors. If we let $\epsilon'$ goes to $0$ the
intersections of the two sectors converge to the segment
$[p_1,p_2]$ which have $(\sin\theta,-\cos\theta)$ as normal. We then
must show that all the points of the segment $[p_1,p_2]$ are the limit
of a sequence $(c_n(t_n))$.

We suppose now that $\theta=0$. By Proposition \ref{convgraph}, the
curves $c_n$ 
are graphs over $\{y=0\}$. Let $a$ and $b$ be the respective first
coordinates of $\widetilde{c} (-\epsilon/2)$ and $\widetilde{c}
(\epsilon/2)$; we suppose $a<b$. Since $c_n\rightarrow \widetilde{c}$,
we then can ensure that, over the segment $[a,b]$, the 
curves $c_n$ are graphs for big $n$. We have $p_1=(x_1,y_1)$ and
$p_2=(x_2,y_2)$; since the segment $[p_1,p_2]$ is horizontal,
$y_1=y_2$. Besides by convergence of $(c_n)$, $a<x_1\le x_2<b$. Let
$x\in[x_1,x_2]$; since $c_n$ is a graph over $[a,b]$ for big $n$,
there exist one parameter $t_n$ such that $c_n(t_n)$ has $x$ as first
coordinate. Then the only possible cluster point for the sequence
$(c_n(t_n))$ is $(x,y_1)$ since every cluster point must have $x$ as
first coordinate and be in the segment $[p_1,p_2]$ or in the curve
$\widetilde{c}$ but all this set is a graph over $\{y=0\}$. This then
proves that all the points of the segment $[p_1,p_2]$ is in the limit
set.  

If $p_1=p_2$, $(c_n)$ converges in the $C^1$ topology since the normal at
the point $\widetilde{c}(\theta)=p_1$ is $(\sin\theta,-\cos\theta)$ by
continuity. 
\end{proof}


\bigskip

\noindent Laurent Mazet

\noindent Laboratoire Emile Picard (UMR 5580), Universit\'e Paul Sabatier,

\noindent 118, Route de Narbonne, 31062 Toulouse, France.

\noindent E-mail: mazet@picard.ups-tlse.fr

\end{document}